\documentclass[onefignum,onetabnum]{siamart190516}

\usepackage{caption}
\usepackage{subcaption}
\usepackage{float}
\usepackage{tikz}
\usetikzlibrary{matrix}
\usetikzlibrary{shapes,patterns,arrows,snakes,decorations.shapes}
\usetikzlibrary{positioning,fit,calc}
\usetikzlibrary{plotmarks}
\usepackage{pgfplots,pgfplotstable}
\pgfplotsset{width=10cm, compat=newest}

\usepackage{amsmath, amsfonts, amssymb,mathrsfs}
\usepackage{amsmath,amssymb}
\usepackage{bm,color}
\usepackage{diagbox}
\newcommand{\bs}{\boldsymbol}

\definecolor{pblue}{rgb}{0.36, 0.54, 0.66}
\definecolor{pred}{rgb}{0.8, 0.25, 0.33}
\definecolor{pyellow}{rgb}{1.0, 0.75, 0.0}
\definecolor{porange}{rgb}{1.0, 0.49, 0.0}
\definecolor{pgreen}{rgb}{0.0, 0.5, 0.0}
\definecolor{pviolet}{rgb}{0.57, 0.36, 0.51}

\ifpdf%
  \DeclareGraphicsExtensions{.eps,.pdf,.png,.jpg}
\else
  \DeclareGraphicsExtensions{.eps}
\fi
\usepackage{amsopn}
\usepackage{booktabs}
\usepackage{lipsum}
\usepackage{amsfonts}
\usepackage{graphicx}
\usepackage{epstopdf}
\usepackage{algorithmic}
\usepackage{cite}

\usepackage{diagbox}
\usepackage{array}
\newlength\colwidth
\newlength\secpad
\newlength\rowheight
\newcolumntype{f}{>{\centering\arraybackslash}m{\colwidth}@{\hskip 0.3em}}

\usepackage{multirow}
\usepackage{soul}

\ifpdf
  \DeclareGraphicsExtensions{.eps,.pdf,.png,.jpg}
\else
  \DeclareGraphicsExtensions{.eps}
\fi

\DeclareMathOperator*{\divg}{\operatorname{div}}
\DeclareMathOperator*{\curl}{\operatorname{curl}}

\DeclareMathOperator*{\grad}{\operatorname{grad}}

\def\Xint#1{\mathchoice
{\XXint\displaystyle\textstyle{#1}}%
{\XXint\textstyle\scriptstyle{#1}}%
{\XXint\scriptstyle\scriptscriptstyle{#1}}%
{\XXint\scriptscriptstyle\scriptscriptstyle{#1}}%
\!\int}
\def\XXint#1#2#3{{\setbox0=\hbox{$#1{#2#3}{\int}$ }
\vcenter{\hbox{$#2#3$ }}\kern-.6\wd0}}

\def\dashint{\Xint-}

\newsiamremark{remark}{Remark}
\newsiamremark{hypothesis}{Hypothesis}
\crefname{hypothesis}{Hypothesis}{Hypotheses}
\newsiamthm{claim}{Claim}

\headers{MFD for Convection-Dominated PDEs}{J. H. Adler, C. Cavanaugh, X. Hu, A. Huang, N. Trask}

\title{A stable mimetic finite-difference method for convection-dominated diffusion equations
\thanks{Revision submitted to the editors May, 22 2023.
}}

\author{James H. Adler\thanks{Department of Mathematics, Tufts University, Medford, MA 02155
  (\email{james.adler@tufts.edu}, \email{xiaozhe.hu@tufts.edu}).} \and Casey Cavanaugh\thanks{Center for Computation and Technology, Louisiana State University, Baton Rouge, LA 70803 \hbox{(\email{Ccavanaugh@cct.lsu.edu})}} \and Xiaozhe Hu\footnotemark[2]
\and Andy Huang \thanks{Radiation and Electrical Science, Sandia National Laboratories
  \hbox{(\email{ahuang@sandia.gov}).}} \and Nathaniel Trask\thanks{Center for Computing Research, Sandia National Laboratories
  \hbox{(\email{natrask@sandia.gov}).}}}

\usepackage{amsopn}

\ifpdf
\hypersetup{
  pdftitle={A stable mimetic finite-difference method for convection-dominated diffusion equations},
  pdfauthor={J. H. Adler, C. Cavanaugh, X. Hu, A. Huang, N. Trask}
}
\fi

\begin{document}

\maketitle

% REQUIRED
\begin{abstract}
Convection-diffusion equations arise in a variety of applications such as particle transport, electromagnetics, and magnetohydrodynamics. Simulation of the convection-dominated regime for these problems, even with high-fidelity techniques, is particularly challenging due to the presence of sharp boundary layers and shocks causing jumps and discontinuities in the solution, and numerical issues such as loss of the maximum principle in the discretization. These complications cause instabilities, admitting large oscillations in the numerical solution when using traditional methods. Drawing connections to the simplex-averaged finite-element method (S. Wu and J. Xu, 2020), this paper develops a mimetic finite-difference (MFD) discretization using exponentially-averaged coefficients to
% \edit{ \st{guarantee monotonicity of the scheme and} overcome in}stability of the
overcome instability of the
numerical solution as the diffusion coefficient approaches zero.  The finite-element framework allows for transparent analysis of the MFD, such as proving well-posedness and deriving error estimates. Numerical tests are presented confirming the stability of the method and verifying the error estimates.
\end{abstract}

% REQUIRED

% REQUIRED

\section{Introduction}\label{sec:intro}

Convection-diffusion equations, particularly describing phenomena in the convection-dominated regime, come in both scalar and vector forms, and have many important applications. In particular, the \emph{scalar} convection-diffusion equation, for unknown $u$, 
\begin{align}\label{eq:convdiff_grad}
\begin{split}
-\divg \left( \alpha \grad u + \bm \beta u \right) + \gamma u &= f \quad \text{in } \Omega,  \\ 
 u &= 0 \quad \text{on } \partial \Omega ,
 \end{split}
\end{align}
is commonly seen in particle transport \cite{PT1, PT2, PT3, PT4, PT5, PT6, PT7}. Here, $\alpha$ and $\gamma$ are real-valued scalar functions and $\bs{\beta}$ is a vector function representing various physical parameters.  The \emph{vector} convection-diffusion equations, which come in two forms,
\begin{align} \label{eq:convdiff_curl}
\begin{split}
\curl \left( \alpha \curl \bm u + \bm \beta \times \bm u \right) + \gamma \bm{u} &= \bm f \quad \text{in } \Omega,  \\ 
\bm n \times  \bm u &= \bm 0 \quad \text{on } \partial \Omega ,
 \end{split}
\end{align}
\vskip -12pt
\begin{align}\label{eq:convdiff_div}
\begin{split}
-\grad \left( \alpha \divg \bm u + \bm \beta \cdot \bm u \right) + \gamma \bm{u} &= \bm f \quad \text{in } \Omega,  \\ 
 \bm u \cdot \bm n  &= 0 \quad \text{on } \partial \Omega ,
 \end{split}
\end{align}
arise in electromagnetic applications such as magnetohydrodynamics \cite{EM1, EM2, EM3, EM4, EM5}, with $\bm \beta$ representing an externally applied field. 

For this paper, we are interested in the convection-dominated case, where $|\alpha| \ll |\bm \beta |$, as it poses difficulties in numerical simulations (see \cite{CDchallenge} for a survey). From PDE theory \cite{EvansPDE}, the convection-diffusion equation is elliptic; however, in the limiting case of $\alpha \to 0$, this is no longer true. As outlined in \cite{CDchallenge}, this loss of ellipticity causes methods developed for the diffusion-dominated (elliptic) problem to fail in the limiting case, resulting in large numerical oscillations when solutions involve internal or boundary layers and shocks. One approach to restoring the stability is developing a \emph{monotone} scheme for which a discrete maximum principle holds. 

For the scalar convection-diffusion equation, (\ref{eq:convdiff_grad}), stabilized and monotone methods have been developed to avoid these common numerical difficulties. One choice is finite-difference or finite-volume schemes using upwinding \cite{Upwind1, Upwind2}. Additionally, there are a variety of finite-element (FE) schemes, like discontinuous Galerkin \cite{StableDG2, StableDG1, blanca1}, streamline-upwind Petrov--Galerkin (SUPG) \cite{Stream1, Stream2}, bubble function stabilized methods \cite{bubble1, bubble2, bubble3, bubble4}, penalty methods \cite{penalty1, penalty2}, and Petrov--Galerkin methods \cite{PG1, PG2} that all provide stabilized discretizations for the scalar equation \eqref{eq:convdiff_grad}. In \cite{blanca2}, an $\bm H(\curl)$-elliptic problem with discontinuous coefficients is considered using DG and auxiliary space preconditioners acheiving scalable results. 
% However, none of these methods provide easy extensions to solving the vector equations (\ref{eq:convdiff_curl}) or (\ref{eq:convdiff_div}). 
Additionally, while there is no explicit maximum principle for the convection-diffusion vector equations, there are a few approaches which provide a unifying framework for the general (both scalar and vector) convection-diffusion equation using Whitney forms \cite{Whit} and $k$-forms \cite{kform1, kform2}, which numerically demonstrate stability in the limiting regime.

The approach that we focus on here is exponential fitting \cite{expfit1, expfit2, expfit3, expfit4, expfit5, EAFE, SAFE, expfit6, expfit7, expfitLZ, expfitOG}. In particular, for \eqref{eq:convdiff_grad} we consider the edge-averaged finite-element method (EAFE) \cite{EAFE}, which uses exponential shifting to recast the convection-diffusion problem as a modified Poisson equation. This approach has many advantages, such as a bilinear form to provide straightforward analysis of the scheme, and a provably monotone matrix formulation. Using Whitney forms and a discrete de Rham complex, the EAFE scheme was extended to the vector convection-diffusion equations, (\ref{eq:convdiff_curl}) and (\ref{eq:convdiff_div}), using the simplex-averaged finite-element method (SAFE) \cite{SAFE}. Again, using exponential shifting with average integrals over the appropriate simplex, each convection-diffusion equation can be recast as a pure diffusion problem with a low-order term added to deal with the large nullspaces of the vector equations. This allows for implementation using exponentially-averaged coefficients in a standard FE approximation of a diffusion equation, and a clear avenue to develop the corresponding analysis and establish monotonicity. Furthermore, the SAFE method establishes a general framework for all three convection-diffusion equations, and ensures that the corresponding exponentially-shifted differential operators, called the flux operators, are structure-preserving in the continuous problem. However, while the SAFE flux operators satisfy a discrete de Rham complex before quadrature is introduced, the exact sequence does not hold in implementation unless the diffusion and convection coefficients are both constant. 

The goal of this paper is to develop a general framework for a stable mimetic finite-difference (MFD) method for the convection-diffusion equation in both scalar and vector forms. There are many advantages to the MFD approach, including the ability to work with general polyhedral meshes, fast assembly of matrix systems, and the natural inheritance of a de Rham complex. Specifically for the convection-diffusion scheme, we show that the MFD method guarantees that the flux differential operators also satisfy a discrete de Rham sequence regardless of how any integrals are computed.   Another major advantage of the stable MFD scheme, and of MFD more generally, is that it is related to discrete exterior calculus (DEC). As opposed to finite-element exterior calculus (FEEC), the DEC approach eliminates the need for mass matrices, which are tied to a specific mesh geometry and dimension. More specifically, MFD guarantees structure-preserving operators which are scaled mesh incidence matrices. The incidence matrices more naturally extend to higher-order problems since they are not tied to any particular coordinate system, paving the way for connections to novel methods in machine learning. For instance, the MFD and DEC approach is advantageous to structure-preserving machine learning methods where a graph incidence matrix of a graph in $\mathbb{R}^n$ is known and the appropriate metrics can be learned using a neural network \cite{DDEC}. Thus, the stable MFD scheme aids the development of structure-preserving machine learning methods for problems where standard approaches struggle due to the curse of dimensionality.  Additionally, there are known connections between the MFD and traditional FE methods \cite{FeMFD, MFDMax} which can be applied to the convection-diffusion scheme. These connections provide a transparent approach to establishing well-posedness and error estimates in the FE setting for the MFD method. Finally, the MFD method yields a linear system resembling a scaled Laplacian from which monotonicity of the scalar convection-diffusion scheme is proven. This validates the stability of the scheme in the convection-dominated regime.

This paper is organized as follows. Section \ref{sec:preliminaries} defines notation and briefly presents the standard FE and MFD methods. In Section \ref{sec:CD_MFD}, 
we introduce the key features of the SAFE method and the corresponding mimetic approach. We also present the full discretized scheme, suggest quadrature rules, and discuss numerical stability in the implementation. Connections to the finite-element method are given in Section \ref{sec:MFD-theory}, along with an analysis of the scheme including well-posedness results and error estimates.  Numerical results are presented in Section \ref{sec:num_CD} and, finally, Section \ref{sec:conc_CD} summarizes the conclusions and presents ideas for future work.

\section{Preliminaries} \label{sec:preliminaries}
In this section, we introduce notation used in the paper and briefly recall the MFD method.  For ease of the analysis, we assume $\Omega \in \mathbb{R}^3$ is a bounded polyhedron domain. However, all the results hold for the two-dimensional case with minor modifications, and the illustrations and numerical results are performed on $\mathbb{R}^2$ for simplicity.  In the rest of this paper, for an open subset $\omega \in \mathbb{R}^d$, we use the standard notation for Sobolev spaces, $W^{\ell,p}(\omega)$, $0 \leq \ell < \infty$, $1 \leq p \leq \infty$, and their associated norm, $\| \cdot \|_{\ell, p, \omega}$, and semi-norm, $|\cdot|_{\ell, p, \omega}$. For $p = 2$, we use the standard notation, $H^{\ell}(\omega):= W^{\ell, 2}(\omega)$ with norm $\| \cdot \|_{\ell, \omega}$ and semi-norm $| \cdot |_{\ell, \omega}$ . Furthermore, for $\ell = 0$, $H^{0}(\omega)$ coincides with the standard $L^2(\omega)$ space with inner product $\langle \cdot, \cdot \rangle_{\omega}$ and norm $\| \cdot \|_{\omega}$.  Finally, when the context is clear, we omit $\omega$ in the notation, especially when $\omega = \Omega$.

\subsection{The Finite-Element Method}\label{sec:FEM}
Let the Sobolev space on $\Omega$ associated with differential operator, $\mathfrak{D}^k$, $k=0,1,2$, be given by,
\begin{equation*}
H(\mathfrak{D}^k) := \{ u \, \big| \, u\in L^2 (\Omega), \, \mathfrak{D}^ku \in L^2(\Omega)\}, %\label{eq:sobolev_defn}
\end{equation*}
with norm $\| \cdot \|_{H(\mathfrak{D}^k)}$, semi-norm $| \cdot |_{H(\mathfrak{D}^k)}$, and inner product $\langle \cdot , \cdot \rangle_{H(\mathfrak{D}^k)}$. 
%Let $\| \cdot \|$ and $\langle \cdot , \cdot \rangle$ denote the $L^2(\Omega)$ norm and inner product, respectively.  
Here, $\mathfrak{D}^0 = \grad$, $\mathfrak{D}^1 =\curl$, $\mathfrak{D}^2 =\divg$, and we have the following de Rham sequence,  
\begin{equation}
H(\grad) \xrightarrow{\grad} \bm H(\curl) \xrightarrow{\curl} \bm H(\divg) \xrightarrow{\divg} L^2. \label{eq:DeRhamL2}
\end{equation}
Note that we use boldface letters for the vector-valued functions and their corresponding function spaces.  To make the notation consistent, we define $H(\mathfrak{D}^3) := L^2(\Omega)$, and as a result we have $\|\cdot\|_{H(\mathfrak{D}^0)} = \|\cdot\|_1$ and $\|\cdot\|_{H(\mathfrak{D}^3)} = \|\cdot\|$, for example.

%\casey{Do we really need $k=0,1,2,3$ or $k=0,1,2$ after each definition? It makes this really bulky. And $\mathcal{M}$ is used for mass matrices later on. How about $\mathcal{T}$?} \james{I think we do for a while here.  \mathcal{T}is a good idea} 
Given a simplicial mesh, $\mathcal{T}$, of the domain $\Omega$, we consider the lowest-order  conforming FE spaces for $H(\mathfrak{D}^k)$, denoted by $H_h(\mathfrak{D}^k)$, $k=0,1,2$.
On the mesh, we have $k$-simplices, denoted by $s^0$ (vertices), $s^1$ (edges), $s^2$ (faces), and $s^3$ (volumes), where $s^k_i$ denotes the $i$-th $k$-simplex in the mesh enumeration. Next, define the  corresponding degrees of freedom (DoFs) $\{ \eta_{i}^{k} \}_{s_i^k \in \mathcal{T}}$ and FE basis functions $\{ \lambda_{i}^{k} \}_{s_i^k \in \mathcal{T}}$, for each space $H(\mathfrak{D}^k)$, $k=0,1,2$.  In addition, a piecewise constant FE space is used for $L^2(\Omega)$ and is denoted as $L^2_h$ (or $H_h(\mathfrak{D}^3)$) with DoFs $\{\eta_{i}^{3}\}$ and basis functions $\{ \lambda_{i}^{3} \}$.  Specifically, for a vertex $i$ with coordinate $\bm{x}_i$, edge $e_i$ with unit tangent $\bm{t}_{i}$ and length $|e_i|$, face $f_i$ with unit normal $\bm{n}_{i}$ and area $|f_i|$, and simplex $T_i$ with volume $|T_i|$, we have the DoFs,
\begin{align*}
\eta_i^{0}(u) &= u(\bm x_i), %\label{eq:DoFgrad}	 
\\
\bm{\eta}_{i}^{1} (\bm{u}) & = \frac{1}{|e_i|}\int_{e_i} ( \bm{u} \cdot \bm{t}_{i}) \mathrm{d}s, %\label{eq:DoFcurl} 
\\
\bm{\eta}_{i}^{2} (\bm{u}) &= \frac{1}{|f_i|}\int_{f_i} ( \bm{u} \cdot \bm{n}_{i}) \mathrm{d}A, %\label{eq:DoFdiv} 
\\
\eta^{3}_{i} (u) &= \frac{1}{|T_i|} \int_{T_i} u \, \mathrm{d}V. %\label{eq:DoFL2}
\end{align*}
%\casey{Should $x_i$, $e$, $f$, $T$, be defined?} \james{done}
%In addition, the finite-element basis functions for $H_h(\grad)$ are the usual barycentric coordinates as follows,
%\begin{equation}
%	\lambda^{\grad}_i(\bm x) = 1 - \frac{|f_i|}{d |T|} \bm n_i \cdot (\bm x - \bm x_i ). \label{eq:BaryBasis}
%\end{equation}

Since \eqref{eq:convdiff_grad}--\eqref{eq:convdiff_div} are given with Dirichlet boundary conditions, the following FE spaces with boundary conditions are considered,
\begin{align*}
H_{h,0}(\grad) &:= \{ u \, \big| \, u \in H_h(\grad), u \big|_{\partial \Omega} = 0\}, %\label{eq:Hhgrad} 
\\
\bm H_{h,0}(\curl) &:= \{ \bm u \, \big| \, \bm u \in \bm{H}_h(\curl), \bm n \times \bm u \big|_{\partial \Omega} = 0\}, %\label{eq:Hhcurl} 
\\
\bm H_{h,0} (\divg) &:=\{ \bm u \, \big| \, \bm u \in \bm{H}_h(\divg), \bm u \cdot \bm n \big|_{\partial \Omega} = 0\}. %\label{eq:Hhdiv}
\end{align*}
It is well-known from finite-element exterior calculus (FEEC) \cite{ArnoldFEEC, FEextcalc, FEextcalc2} that the corresponding function spaces for these choices of elements satisfy a de Rham sequence,
\begin{equation}
	H_{h,0}(\grad) \xrightarrow{\grad} \bm{H}_{h,0}(\curl) \xrightarrow{\curl} \bm{H}_{h,0}(\divg) \xrightarrow{\divg} L^2_h.\label{eq:FEEC_exact}
\end{equation}

Next, we introduce the corresponding discrete differential operators. The FE gradient, $G_{\text{FE}}$, curl, $K_{\text{FE}}$, and divergence, $D_{\text{FE}}$, operators are expressed in terms of the basis functions and DoFs,
\begin{align*} 
\left( G_{\text{FE}} \right)_{ij}&:= \bm{\eta}_{i}^{1} \left( \grad \lambda_j^{0} \right) = \frac{1}{|e_i|} \int_{e_i} \grad \lambda^{0}_j \cdot \bm{t}_i \, \mathrm{d}s, 
%\label{eq:FEgradPrelim} 
\\ 
\left( K_{\text{FE}} \right)_{ij}&:= \bm{\eta}_{i}^{2} \left( \curl \bm \lambda_{j}^{1} \right) = \frac{1}{|f_i|} \int_{f_i} \curl \bm \lambda^{1}_{j} \cdot \bm{n}_i \, \mathrm{d}A, 
%\label{eq:FEcurlPrelim}
\\ 
\left( D_{\text{FE}} \right)_{ij}&:= \eta_{i}^{3} \left( \divg \bm \lambda_{j}^{2} \right) = \frac{1}{|T_i|} \int_{T_i} \divg \bm \lambda^{2}_{j}  \, \mathrm{d}V. 
%\label{eq:FEdivPrelim}
\end{align*}
%where $e_i$, $f_i$, and $T_i$ are the $i$th edge, face, and element in the mesh enumeration, respectively \casey{not consistent with the notation defining the DoFs earlier}. 
By direct calculation, we have $K_{\text{FE}} G_{\text{FE}} = 0$ and $D_{\text{FE}} K_{\text{FE}} = 0$ which verifies properties of \eqref{eq:FEEC_exact}.

Finally, we introduce the usual canonical interpolation operator, which plays an essential role in drawing the connections between FEM and MFD schemes, 
\begin{equation}\label{eq:can_interp}
	\Pi^{k} v = \sum_{_{s_i^k \in \mathcal{T}}} \eta_{i}^{k}(v) \lambda_{i}^{k}, \quad  \forall \ v \in H(\mathfrak{D}^k).
\end{equation}
%More precisely, 
%\begin{align}
%\Pi^{0} v &= \sum_i \eta_{i}^{0}(v) \lambda_i^{0}, \quad  \forall \ v \in H(\grad), \label{eq:Pi_grad} \\
%\bm{\Pi}^{1} \bm{v} &= \sum_{e_i} \bm{\eta}_{e_i}^{1}(\bm{v}) \bm{\lambda}_{e_i}^{1}, \quad  \forall \ \bm{v} \in \bm{H}(\curl), \label{eq:Pi_curl} \\
%\bm{\Pi}^{2} \bm{v} &= \sum_{f_i} \bm{\eta}_{f_i}^{2}(\bm{v}) \bm{\lambda}_{f_i}^{2}, \quad  \forall \ \bm{v} \in \bm{H}(\divg), \label{eq:Pi_div} \\	
%\Pi^{3} v & = \sum_{T_i}  \eta_{T_i}^{3}(v) \lambda_{T_i}^{3}, \quad  \forall \ v \in L^2. \label{eq:Pi_L2}
%\end{align}
We also define notation for the average integral,
\begin{equation*}
	\dashint_{\mathcal{S}} f(\bm{x}) \, \mathrm{d} \bm{x} := \frac{1}{| \mathcal{S}|} \int_{\mathcal{S}} f(\bm{x}) \, \mathrm{d}\bm{x}, 
	%\label{eq:avgint}
\end{equation*}
where $\mathcal{S}$ denotes the region of the integration. 

\subsection{The Mimetic Finite-Difference Method}\label{sec:MFD}
We next introduce the notation and differential operators for MFD. The foundation of the MFD method is a dual mesh configuration. The primal mesh is a Delaunay triangulation in 2 or 3 dimensions, with corresponding dual Voronoi polygonal or polyhedral mesh. Note that the finite-element mesh, $\mathcal{T}$ defined above, corresponds to the primal Delaunay mesh defined here, but we change the notation slightly to more clearly describe the connection between the primal and dual meshes in the finite-difference setting.  

Let $N_D$ denote the number of Delaunay nodes $\bm{x}_i^D$. An edge connecting nodes $\bm{x}_i^D$ and $\bm{x}_j^D$ is given by $\bm{e}_{ij}^D$ with unit tangent vector $\bm{t}_{ij}^D$ pointing from vertex of lower index to vertex of higher index. The edge tangent vector (non-unit) is given by $\bm{T}_{ij}^D =\bm{t}_{ij}^D |\bm{e}_{ij}^D|$. A Delaunay element is given by $D_k$ with boundary set $\partial D_k$. The set of indices of the elements that share a boundary with $D_k$ is given by $\mathcal{N}_k^D := \{ m \, \big| \, \partial D_k \cap \partial D_m \neq \emptyset, m = 1, \ldots, N_V \}$. The shared boundary element (a face in 3D or edge in 2D) is given by $\partial D_{km}$ which has an associated unit normal vector $\bm{n}_{km}^D$, which points outward from element $D_k$. Similarly, on the Voronoi mesh we define number of nodes $N_V$, nodes $\bm{x}_k^V$, edge $\bm{e}_{km}^V$ with unit tangent $\bm{t}_{km}^V$, element $V_i$, boundary set $\partial V_i$, neighbor set $\mathcal{N}_i^V$, and boundary element $\partial V_{ij}$ with outward unit normal $\bm{n}_{ij}^V$. 

\begin{figure}[H]
\begin{center}

%\begin{tikzpicture}[scale = 5, every node/.style={scale = .9}]
\begin{tikzpicture}[xscale=5.5,yscale=5, every node/.style={scale = .9}]
\draw[thick] (.5,0) to (.25, .433) to (-.25, .433) to (-.5, 0) to (-.25, -.433) to (.25, -.433) to cycle;
\draw[thick] (.25, .433) to (-.25, -.433);
\draw[thick] (-.25, .433) to (.25, -.433);
\draw[thick] (.5,0) to (-.5, 0);
\draw[thick] (.25, .433) to (.75, .433) to (1,0) to (.75, -.433) to (.5, 0) to cycle;
\draw[thick] (.75, -.433) to (.25, -.433);

\draw[line width=2] (.5,0) to (1,0);
\draw[line width=2] (.5,0) to (.75, .433);
\draw[line width=2] (1,0) to (.75, .433);

\draw[dashed, red, line width=2] (.25, .1443) to (.25, -.1433) to (0, -.2887) to (-.25, -.1443) to (-.25, .1443) to (0, .2887) to cycle;
\draw[dashed, red, thick] (0, .2887) to (0, .433);
\draw[dashed, red, thick] (0, -.2887) to (0, -.433);
\draw[dashed, red, thick] (-.25, .1433) to (-0.375, 0.2165);
\draw[dashed, red, thick] (.25, .1433) to (0.5, .2887);
\draw[dashed, red, thick] (-.25, -.1433) to (-0.375, -0.2165);

\draw[dashed, red, thick](.5, .2887) to (.5, .433);
\draw[dashed, red, thick](.5, -.2887) to (.5, -.433);
\draw[dashed, red, thick] (.75, .1433) to (.875  , 0.2165);
\draw[dashed, red, thick] (.75, -.1433) to (.875  , -0.2165);

%\draw[dashed, red, thick] (0.5, .2887) to  (.75, .1433) to (.75, -.1433) to (0.5, -.2887) to (.25, -.1443);
\draw[dashed, red, thick] (0.5, .2887) to  (.75, .1433);
\draw[dashed, red, thick] (.75, -.1433) to (0.5, -.2887) to (.25, -.1443);
%changes linewidth
%\draw[dashed, red, line width = 3.5] (.75, .1433) to (.75, -.1433);
\draw[dashed, red, thick] (.75, .1433) to (.75, -.1433);

\fill (0,0) circle (.5pt);
\fill (.5,0) circle (.5pt);
\fill (.25, .433) circle (.5pt);
\fill (-.25, .433) circle (.5pt);
\fill (-.5, 0) circle (.5pt);
\fill (-.25, -.433) circle (.5pt);
\fill (.25, -.433) circle (.5pt);

\fill (.75, .433) circle (.5pt);
\fill (.75, -.433) circle (.5pt);
\fill (1,0) circle (.5pt);

\fill[red] (.25, .1443) circle (.5pt);
\fill[red] (.25, -.1433) circle (.5pt);
\fill[red] (0, -.2887) circle (.5pt);
\fill[red] (-.25, -.1443) circle (.5pt);
\fill[red] (-.25, .1443) circle (.5pt);
\fill[red] (0, .2887) circle (.5pt);

\fill[red] (.75, .1433) circle (.5pt);
\fill[red] (.75, -.1433) circle (.5pt);
\fill[red] (.5, -0.2887) circle (.5pt);
\fill[red] (0.5, .2887) circle (.5pt);

%del nodes
\node at (0.01, .1) {$\bm{x}_1^D$};
\node at (-.54, .1) {$\bm{x}_2^D$};

%vor nodes
\node at (.75, -.22){$\bm{x}_{2}^V$};
\node at (.75, .22) {$\bm{x}_1^V$};

\node at (-.4, .06) {$\bm e^D_{12}$};
%\node at (-.11, -.05) {$\partial V_{12}$};
\draw[->, very thick, red](-.45, .4)  to[bend left = 27] (-.1, 0.05);
\node at (-.5, .4) {$V_1$};

%\node at (-.38, .06) {$\bm e^D_{12}$};

\node at (-.4, -.06) {$\partial V_{12}$};
\draw[->, thick, red](-.345, -.05)   to (-.26, -.05) ;
\draw[->, thick](-.36, .06)  to[bend left = 35] (-.3, 0.005);

\draw[->, very thick, black](1, .4)  to[bend right = 27] (.75, .32);
\node at (1.05, .4) {$D_1$};
%\node at (.86, .06) {$\partial D_{12}$};
%\node at (.69, -.06) {$\bm e^V_{12}$};

\node at (.9, .06) {$\partial D_{12}$};
\node at (.9, -.06) {$\bm e^V_{12}$};
\draw[->, thick, red](.85, -.06)   to (.76, -.06) ;
\draw[->, thick](.83, .06)  to[bend right = 40] (.785, .005);

\end{tikzpicture}
\caption{Two-dimensional MFD mesh with Voronoi element $V_1$ corresponding to Delaunay point $\bm{x}_1^D$, and Delaunay element $D_1$ corresponding to Voronoi point $\bm{x}_1^V$. Delaunay edge $\bm{e}_{12}^D$ is orthogonal to Voronoi face (edge) $\partial V_{12}$, and similarly, Voronoi edge $\bm{e}_{12}^V$ is orthogonal to Delaunay face (edge) $\partial D_{12}$. } \label{fig:MFD_mesh}
\end{center}
\end{figure}
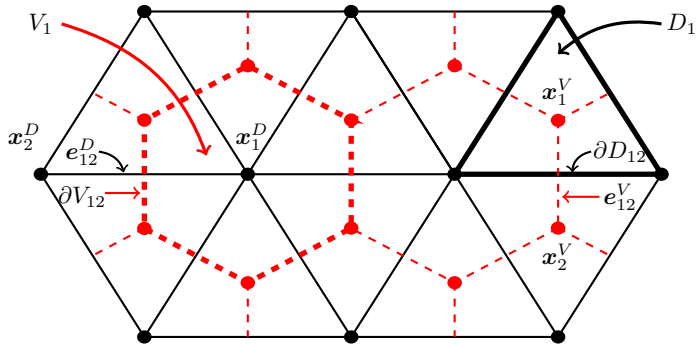

The MFD dual mesh configuration presents many useful geometric relationships. First, there are one-to-one relationships between nodes on one mesh and elements on the other. In particular, each Voronoi vertex $\bm{x}_k^V$ is the circumcenter of the Delaunay element $D_k$, and each Voronoi element $V_i$ is the set of points in the domain that lies closest to Delaunay point $\bm{x}_i^D$, $ V_i := \{ \bm{x} \in \Omega \, \big| \, | \bm{x} - \bm{x}_i^D| \leq | \bm{x} - \bm{x}_j^D|, j = 1, \ldots, N_D, j\neq i \}$. Another correspondence is between edges on one mesh and faces on the other. Delaunay edge $\bm{e}_{ij}^D$ and Voronoi face $\partial V_{ij}$ are orthogonal such that $\bm{t}_{ij}^D \cdot \bm{n}_{ij}^V = \pm 1$. Analogously, we have that $\bm{t}_{km}^V \cdot \bm{n}_{km}^D = \pm 1$ resulting in Voronoi edge $\bm{e}_{km}^V$ and Delaunay face $\partial D_{km}$ also being orthogonal. These relationships are demonstrated in Figure \ref{fig:MFD_mesh}. Note that the dual grid is not required to be a Voronoi diagram, though the construction of an alternate dual grid should maintain the same properties: the one-to-one relationships previously mentioned, and orthogonality of edges on one mesh with faces on the other.

%\subsection{Grid Functions and MFD Operators} 
%\casey{Should we be careful using $k$ as an MFD index here? Or do you think it's clear? $k$ is typically ised to denote Voronoi objects.}\james{I think it's OK here, it's just an index not a variable} 
As in \cite{MFD_Vabish, MFDBook_lipnikov}, we proceed by defining scalar and vector function spaces on the two meshes, and the usual discrete differential operators. The space of scalar functions will be approximated as piecewise constants on the nodes of the mesh as follows:
\begin{align*}
H_D &:= \{ u(\bm{x}) \, \big| \, u(\bm{x})  = u_i^D := u(\bm{x}_i^D), \,\, \forall \bm{x} \in V_i, \, i = 1,...,N_D \}, 
%\label{eq:scalarD}  
\\
H_V &:= \{ u(\bm{x}) \, \big| \, u(\bm{x}) = u_k^V := u(\bm{x}_k^V) , \,\, \forall \bm{x} \in D_k, \, k = 1,...,N_V \}. 
%\label{eq:scalarV}
\end{align*}

%Thus, scalar functions on the Delaunay nodes, $H_D$, are constant on Voronoi elements, and likewise, scalar functions on Voronoi nodes, $H_V$, are constant on Delaunay elements.

Vector functions, $\bm{u}(\bm{x})$, on the other hand, are approximated along the mesh edges. The Delaunay space of discrete vector functions, $\bm{H}_D$, is formed by taking projections onto edges and evaluating at the intersection of the Delaunay edge and corresponding Voronoi interface. The space of Voronoi vector functions, $\bm{H}_V$, is defined analogously:
\begin{align*}
\bm{H_D} &:= \{ \bm{u}(\bm{x}) \, \big| \, \bm{u}(\bm{x}) =  u_{ij}^D := (\bm{u} \cdot \bm{t}_{ij}^D) (\bm{x}_{ij}^D) , \,\, \bm{x}_{ij}^D = \bm{e}_{ij}^D \cap \partial V_{ij} \}, %\label{eq:vectorD} 
\\
\bm{H_V} &:= \{ \bm{u}(\bm{x}) \, \big| \, \bm{u}(\bm{x}) = u_{km}^V := (\bm{u} \cdot \bm{t}_{km}^V) (\bm{x}_{km}^V) , \,\, \bm{x}_{km}^V = \bm{e}_{km}^V \cap \partial D_{km} \}.
%\label{eq:vectorV}.
\end{align*}

The usual differential operators $\grad$, $\curl$, $\divg$ are well-known and easily derived from discrete versions of the calculus theorems \cite{MFD_Vabish, MFDBook_lipnikov, Manzini_MFD}. 
On the Delaunay mesh, we have ${\grad}_D u : H_D \to \bm{H_D}$ on edge $\bm{e}_{ij}^D$,  
\begin{equation}
({\grad}_D u )_{ij}^D = \frac{u_j^D - u_i^D}{|\bm{e}_{ij}^D|} \,\, \zeta(i,j), \label{eq:gradD_comp}
\end{equation}
where $\zeta$ is an orientation constant that $\zeta(i,j) = 1$ if $j > i$ and, otherwise, $\zeta = -1$.
%\begin{equation*}
%\zeta (i,j) = \begin{cases} \,\,\,\,1, & j>i \\
%						 -1, & \text{otherwise,}
%			\end{cases}
%\end{equation*}			  
We also have ${\curl}_D  u: \bs{H_D} \to \bs{H_V}$, measuring the rotation on Delaunay face $\partial D_{km}$, 
\begin{equation}
({\curl}_D  u)_{km}^V = \frac{(\bs t_{km}^V \cdot \bs n_{km}^D)}{|\partial D_{km}|} \sum_{\bm{e}_{ij}^D \in \partial D_{km}} u_{ij}^D \,\, |\bm{e}_{ij}^D| \,\, \chi(\bs n_{km}^D,\bs t_{ij}^D), \label{eq:curlD_comp}
\end{equation}
where the constant $\chi$ enforces the direction of rotation, i.e., $\chi(\bs n_{km}^D, \bs t_{ij}^D) = 1$ if $\bm t_{ij}^D$ is positively oriented and, otherwise, $\chi(\bs n_{km}^D, \bs t_{ij}^D) = -1$.
%\begin{equation*}
%\chi(\bs n_{km}^D, \bs t_{ij}^D) = \begin{cases}
%\,\,\,\,1, & \bs t_{ij}^D \text{ positively oriented,} \\
%-1, & \text{otherwise,}
%\end{cases}
%\end{equation*}
In addition, we have the Voronoi divergence  ${\divg}_V u : \bs{H_V} \to H_V$, corresponding to the outward flux of $D_k$,
\begin{equation}
({\divg}_V u)_k^V = \frac{1}{|D_k|} \sum_{m \in \mathcal{N}_k^D} |\partial D_{km}| u_{km}^V (\bs n_{km}^D \cdot \bm t_{km}^V). \label{eq:divV_comp}
\end{equation}
The Voronoi mesh operators $\grad_V$, $\curl_V$, and Delaunay $\divg_D$ are defined in a similar fashion.  To write the operators in matrix form, we define the Delaunay edge-vertex signed incidence matrix, $\mathcal{G} \in \mathbb{R}^{M_D \times N_D}$, the Voronoi edge-vertex signed incidence matrix, $\mathcal{G}_{V} \in \mathbb{R}^{M_V \times N_V}$, and the Delaunay face-edge signed incidence matrix, $\mathcal{K} \in \mathbb{R}^{M_V \times M_D}$ where $M_V$ is the number of Voronoi edges and $M_D$ the number of Delaunay edges. Assuming a pre-determined orientation of the edges and faces, the nonzero entries of $\mathcal{G}$, $\mathcal{G}_V$, and $\mathcal{K}$ are either $1$ or $-1$ and capture the action of the operators. To scale appropriately, define the following diagonal matrices encoded with metric information on the meshes: \\
\begin{center}
\begin{tabular}{ccc}
$\mathcal{D}_{\bm{e}^D} = \text{diag} \left(|\bm{e}_{ij}^D| \right)$, &  
$\mathcal{D}_{\partial D} = \text{diag} \left(|\partial D_{km}| \right)$, & $\mathcal{D}_{D} = \text{diag} \left(|D_k| \right)$, 
\\[3pt]
$\mathcal{D}_{\bm{e}^V} = \text{diag} \left(|\bm{e}_{km}^V| \right)$,
& $\mathcal{D}_{\partial V} = \text{diag} \left(|\partial V_{ij}| \right)$, 
&
$\mathcal{D}_{V} = \text{diag} \left(|V_i| \right)$. 
\end{tabular}
\end{center}
Thus,
\begin{align*}
 %\label{eq:graddivcurlD}
 {\grad}_D &:= \mathcal{D}_{\bm{e}^D}^{-1} \mathcal{G},& 
 {\divg}_D &:= \mathcal{D}_V^{-1} \mathcal{G}^T \mathcal{D}_{\partial V},& {\curl}_D &:= \mathcal{D}_{\partial D}^{-1} \mathcal{K} \mathcal{D}_{\bm{e}^D},\\
%\label{eq:graddivcurlV}
{\grad}_V &:= \mathcal{D}_{\bm{e}^V}^{-1} \mathcal{G}_V, &
{\divg}_V &:= \mathcal{D}_D^{-1} \mathcal{G}_V^T \mathcal{D}_{\partial D}, &
{\curl}_V &:= \mathcal{D}_{\partial V}^{-1} \mathcal{K}^T \mathcal{D}_{\bm{e}^V}.
\end{align*} 

A significant advantage of the mimetic differential operators is that they are structure-preserving, meaning that standard properties are preserved: $\text{curl}_D \text{grad}_D = 0$, $\text{curl}_V \text{grad}_V = 0$, $\text{div}_V \text{curl}_D = 0$, and $\text{div}_D \text{curl}_V = 0$ \cite{MFD_Vabish, MFDBook_lipnikov, Manzini_MFD}. In essence, this means that the operators satisfy two exact sequences,
\begin{align}
H_D &\xrightarrow{\text{grad}_D} \bm{H_D} \xrightarrow{\text{curl}_D} \bm{H_V} \xrightarrow{\text{div}_V} H_V, \label{eq:MFDderham1} \\
H_V &\xrightarrow{\text{grad}_V} \bm{H_V} \xrightarrow{\text{curl}_V} \bm{H_D} \xrightarrow{\text{div}_D} H_D.\label{eq:MFDderham2}
\end{align}

\section{The MFD Method for the Convection Diffusion Problems} \label{sec:CD_MFD}
In this section, we introduce the proposed MFD scheme, and present the full stable discretization for the convection-diffusion problems, along with some implementation details. First, though, we review the key definitions and features of the SAFE method \cite{SAFE}, which our approach is based on. 

 \subsection{Shifted Flux Operators}\label{sec:FluxOps}
 The SAFE framework relies on exponential shifting of the usual $\grad$, $\curl$, and $\divg$ operators to reformulate the convection-diffusion equation as a diffusion  problem with exponential coefficients. The goal is to preserve the de Rham complex \eqref{eq:DeRhamL2} with the new, shifted operators. 

To start, define $\bm{\theta} = \bm \beta / \alpha$ as the ratio of convection to diffusion parameters, and introduce a potential function, $\varphi$,  such that
\begin{equation}
\bm{\theta} = \grad \varphi. \label{eq:potential}
\end{equation}
Then, the shifted flux operators are obtained as follows:
\begin{align}
J^0 u &= \grad  u + \bm \theta u = e^{- \varphi} \grad \left( e^{\varphi} u \right),\label{eq:flux_grad}\\  
J^1 \bm u &= \curl \bm u + \bm \theta \times \bm u = e^{- \varphi} \curl \left( e^{\varphi} \bm u \right),\label{eq:flux_curl} \\ 
J^2 \bm u &= \divg \bm u + \bm \theta \cdot \bm u = e^{- \varphi} \divg \left( e^{\varphi} \bm u \right). \label{eq:flux_div} 
\end{align}
This shifting can be thought of as an exponential transform on the function $u$, applying the derivative, and transforming back to the original space. 
%The transform properties (\ref{eq:flux_grad})--(\ref{eq:flux_div}) are verified by applying the product rule. For example, taking (\ref{eq:flux_grad}),
%\begin{align*}
%J_0 u = e^{-\varphi} \grad \left( e^{\varphi} u \right) &= e^{-\varphi} \left[  e^{\varphi} \grad u + u \grad \left( e^{\varphi} \right)  \right] \\
%&= e^{-\varphi} \left[  e^{\varphi} \grad u +  u e^{\varphi} \grad \varphi  \right] \\
%&=  \grad u  + u \grad \varphi\\
%&=  \grad u + u \bm \theta .
%\end{align*}
Using this formulation, it is easily verified that,
\begin{align*}
J^1 \circ J^0 u &= e^{- \varphi} \curl \left( e^{\varphi} e^{- \varphi} \grad \left( e^{\varphi} u \right) \right) = e^{- \varphi} \curl \left( \grad \left( e^{\varphi} u \right) \right)= \bm{0}, \\
J^2 \circ J^1 \bm{u} &= e^{- \varphi} \divg \left( e^{\varphi} e^{- \varphi} \curl \left( e^{\varphi} \bm u \right) \right) = e^{- \varphi} \divg \left( \curl \left( e^{\varphi} \bm u \right) \right)=  0.
\end{align*}
Therefore, the flux operators $J^0$, $J^1$, and $J^2$ correspond to shifted $\grad$, $\curl$, and $\divg$ operators, respectively, and admit the desired de Rham exact sequence,
\begin{equation}
H(\grad) \xrightarrow{J^0} \bm H(\curl) \xrightarrow{J^1} \bm{H}(\divg) \xrightarrow{J^2} L^2. \label{eq:SAFEderham} 
\end{equation}
The benefit of this approach is that the convection term is built into the shifted differential operator. Thus, the convection-diffusion equations~\eqref{eq:convdiff_grad}--\eqref{eq:convdiff_div} are re-written as modified diffusion problems: 
\begin{align}\label{eq:convdiff_grad_SAFE}
\begin{split}
-\divg \left(  \alpha  J^0 u \right) + \gamma u &= f \quad \text{in } \Omega,  \\ 
 u &= 0 \quad \text{on } \partial \Omega; 
 \end{split}
\end{align}
\vskip -10pt
\begin{align}\label{eq:convdiff_curl_SAFE}
\begin{split}
\curl \left( \alpha J^1 \bm{u} \right) + \gamma \bm u &= \bm f \quad \text{in } \Omega,  \\ 
\bm n \times  \bm u &= \bm 0 \quad \text{on } \partial \Omega ;
 \end{split}
\end{align}
\vskip -10pt
\begin{align}\label{eq:convdiff_div_SAFE}
\begin{split}
-\grad \left( \alpha J^2 \bm u \right) + \gamma \bm u &= \bm f \quad \text{in } \Omega,  \\ 
 \bm u \cdot \bm n  &= 0 \quad \text{on } \partial \Omega.
 \end{split}
\end{align}

\begin{remark} \label{remark:helmholtz}
This framework relies on~\eqref{eq:potential}, the ability to re-write the ratio of the convection to diffusion parameters in terms of a potential function. Although this can be done for a wide range of functions,  for arbitrary function $\bm \theta(\bm{x})$, we should not expect this to always be possible.  One solution to handle this case would be to take a piecewise constant approximation of $\bm \theta$ on the mesh, as is done in \cite{SAFE}. However, this may destroy the de Rham complex~\eqref{eq:SAFEderham}. An alternative option would be to decompose $\bm \theta$ using a Helmholtz decomposition, i.e., $\bm{\theta} = \grad  \varphi + \curl \bm{\psi}$, and build the flux operators using the $\grad \varphi$ part.  Section \ref{subsec:helm_num} provides numerical examples for the scalar convection-diffusion equation, where the convection coefficient is written in terms of this Helmholtz decomposition.  In many practical applications, though, such as the drift-diffusion equations, the convection coefficient is naturally given as the gradient of a potential function. 
\end{remark}

\subsection{Discrete MFD flux operators} 
Next, we consider a discretization of the shifted flux operators~\eqref{eq:flux_grad}--\eqref{eq:flux_div} using MFD. In particular, we choose to discretize $J^0$ and $J^1$ on the Delaunay mesh, and $J^2$ on the Voronoi mesh to recover something resembling the  MFD de Rham complex \eqref{eq:MFDderham1}. 
%Similarly, one could place $J^0$ and $J^1$ on the Voronoi mesh, and $J^2$ on the Delaunay mesh to yield a relation like \eqref{eq:MFDderham2}, but we focus on the former in this paper. 

Starting with the $J^0$ operator, rearrange~\eqref{eq:flux_grad} slightly and integrate over a Delaunay edge $\bm{e}_{ij}^D$, to obtain 
\begin{equation*}
\dashint_{\bm e_{ij}^D} e^{\varphi} J^0 u \cdot \bm{t}_{ij}^D \, \mathrm{d}s = \dashint_{\bm e_{ij}^D} \grad \left( e^{\varphi} u \right) \cdot \bm{t}_{ij}^D \, \mathrm{d}s. %\label{eq:J0_deriv}
\end{equation*}
The right-hand side can be directly integrated using the fundamental theorem of calculus. On the left-hand side, we use the basic idea of an MFD method and approximate~$J^0 u \cdot \bm{t}_{ij}^D$ by a constant on the edge, denoted by~$(\mathcal{J}^0_D u)_{ij}^D$. Thus, denoting $\varphi_i = \varphi(\bm{x}_i^D)$, the discrete MFD flux operator,~$\mathcal{J}^0_D : H_D \to \bm{H_D}$, is defined as
\begin{equation}
\left( \mathcal{J}^0_D u \right)_{ij}^D :=  \frac{e^{\varphi_j} u_j^D - e^{\varphi_i} u_i^D}{|\bm{e}_{ij}^D| \, \dashint_{\bm e_{ij}^D} e^{\varphi} \mathrm{d}s} \zeta(i,j). \label{eq:J0_MFD}
\end{equation}
Next, for $J^1$, rearrange~\eqref{eq:flux_curl} and integrate over a Delaunay face to obtain
\[\dashint_{\partial D_{km}} e^{\varphi} J^1 \bm{u} \cdot \bm{n}_{km}^D \, \mathrm{d}A = \dashint_{\partial D_{km}} \curl \left( e^{\varphi} \bm{u} \right) \cdot \bm{n}_{km}^D \, \mathrm{d}A. \]
Approximating $J^1 u \cdot \bm{n}_{km}^D$ with constant $(\mathcal{J}^1_D u)_{km}^V$ on the face and directly evaluating the right-hand side using Stokes' Theorem, we define the following discrete MFD flux operator $\mathcal{J}_D^1  : \bm{H_D} \to \bm{H_V}$,
\begin{equation}
\left( \mathcal{J}^1_D u \right)_{km}^V := (\bm n_{km}^D, \bm t_{ij}^D) \sum_{\bm e_{ij}^D \in \partial D_{km}} \frac{|\bm{e}_{ij}^D| \, \dashint_{\bm e_{ij}^D} e^{\varphi} \mathrm{d}  s}{|\partial D_{km}| \, \dashint_{\partial D_{km}} e^{\varphi} \mathrm{d}A} \, u_{ij}^D \, \chi(\bm n_{km}^D, \bm t_{ij}^D). \label{eq:J1_MFD}
\end{equation}
Following the same procedure for approximating $J^2$, rearrange~\eqref{eq:flux_div} and integrate over a Delaunay element to obtain 
\[ \dashint_{D_k} e^{\varphi} J^2 \bm{u} \, \mathrm{d}V = \dashint_{D_k} \divg \left( e^{\varphi} \bm{u} \right) \mathrm{d}V. \]
Here, the right-hand side can be directly evaluated using the divergence theorem, and a constant approximation of $J^2 \bm u$ on the element $D_k$ yields the MFD flux operator $\mathcal{J}^2_V : \bm{H_V} \to H_V$, 
\begin{equation}
\left( \mathcal{J}^2_V u \right)_{k}^V := \sum_{m \in \mathcal{N}_k^D} \frac{ |\partial D_{km}| \, \dashint_{\partial D_{km}} e^{\varphi} \mathrm{d} A}{| D_{k}|\, \dashint_{D_k} e^{\varphi} \mathrm{d}V} \, u_{km}^V \left( \bm n_{km}^D \cdot \bm t_{km}^V \right). \label{eq:J2_MFD}
\end{equation}

All three discrete MFD flux operators, $\mathcal{J}^0_D$, $\mathcal{J}^1_D$, and $\mathcal{J}^2_V$,  have the same form as the usual MFD ${\grad}_D$~\eqref{eq:gradD_comp}, ${\curl}_D$~\eqref{eq:curlD_comp} and ${\divg}_V$~\eqref{eq:divV_comp}.  The only difference is the scaling, which reflects the exponential shifting seen as the average exponential integrals in the coefficients. Therefore, to write the discrete MFD flux operators in matrix form, we use the incidence matrices $\mathcal{G}$ and $\mathcal{K}$, along with the following diagonal matrices for the average exponential integrals,
\begin{align}
\begin{split}
E_{\bm{x}^D} = \text{diag} \left( e^{\varphi(\bm{x}_i^D)} \right), \,\,&E_{\bm{e}^D} = \text{diag} \left(\dashint_{\bm e_{ij}^D} e^{\varphi} \mathrm{d}s \right),\\ 
E_{\partial D} = \text{diag} \left(\dashint_{\partial D_{km}} e^{\varphi} \mathrm{d}A\right), \,\, &E_{D} = \text{diag} \left( \dashint_{D_k} e^{\varphi} \mathrm{d}V \right). \label{eq:exp_avg_mat}
\end{split}
\end{align}
Thus, the matrix form of the discrete flux operators is given by,
\begin{align}
 \mathcal{J}^0_D &:= E_{\bm{e}^D}^{-1} \mathcal{D}_{\bm{e}^D}^{-1} \mathcal{G}E_{\bm{x}^D} = E_{\bm{e}^D}^{-1} \, {\grad}_D\, E_{\bm{x}^D},\label{eq:J0Mat} \\ 
 \mathcal{J}^1_D &:= E_{\partial D}^{-1} \mathcal{D}_{\partial D}^{-1} \mathcal{K} \mathcal{D}_{\bm{e}^D} E_{\bm{e}^D} = E_{\partial D}^{-1}\, {\curl}_D \, E_{\bm{e}^D}, \label{eq:J1Mat}\\
 \mathcal{J}^2_V &:= E_D^{-1} \mathcal{D}_{D}^{-1} \mathcal{G}_V^T \mathcal{D}_{\partial D} E_{\partial D} = E_D^{-1} \, {\divg}_V \, E_{\partial D}. \label{eq:J2Mat} 
\end{align} 
%which is equivalent to,
%\begin{align}
% \mathcal{J}_0 &:= E_{\bm{e}^D}^{-1} \, {\grad}_D\, E_{\bm{x}^D},\label{eq:J0Mat} \\ 
% \mathcal{J}_1 &:= E_{\partial D}^{-1}\, {\curl}_D \, E_{\bm{e}^D}, \label{eq:J1Mat}\\  
% \mathcal{J}_2 &:= E_D^{-1} \, {\divg}_V \, E_{\partial D}. \label{eq:J2Mat}
%\end{align}
From this, it is easy to show that $\mathcal{J}^1_D \mathcal{J}^0_D = 0$ and $\mathcal{J}^2_V \mathcal{J}^1_D = 0$, and that the MFD exact sequence is preserved,
\[
H_D \xrightarrow{\mathcal{J}^0_D} \bm{H_D} \xrightarrow{\mathcal{J}^1_D} \bm{H_V} \xrightarrow{\mathcal{J}^2_V} H_V. \label{eq:SAFEMFDderham} \]
Note that this holds regardless of how the exponential integrals are computed, so long as the same approach (i.e., quadrature rule) is used consistently across the three discrete MFD flux operators.

\begin{remark}
The above sequence corresponds to the MFD sequence \eqref{eq:MFDderham1}, which relates to solving (3.6)--(3.8) in primal form. The derivation of the flux operators in Section 3.2 will give the flux operators corresponding to the second MFD exact sequence \eqref{eq:MFDderham2} by simply swapping primal mesh quantities ($\bm{x}_i^D$, $\bm{e}_{ij}^D$, $\partial D_{km}$, and $D_k$) for the corresponding dual mesh quantity ($\bm{x}_k^V$, $\bm{e}_{km}^V$, $\partial V_{ij}$, and $V_i$, respectively), and vice versa.
\end{remark}

\subsection{Full MFD Scheme and Implementation} \label{sec:full_disc}
With the flux operators defined in the MFD setting~\eqref{eq:J0Mat}--\eqref{eq:J2Mat}, the fully discretized convection-diffusion equations~\eqref{eq:convdiff_grad_SAFE}--\eqref{eq:convdiff_div_SAFE} can be written:
Find $u_D \in H_D$ such that,
\begin{equation}\label{eq:convdiff_grad_MFD}
\begin{split}
{\divg}_D \left( \mathcal{D}_{\alpha,\bm{e}^D} \, \mathcal{J}^0_D u_D \right) + \mathcal{D}_{\gamma, \bm{x}^D} \, u_D &= f_D \quad \text{in } \Omega,  \\ 
 u_D &= 0 \quad \text{on } \partial \Omega ;
 \end{split}
\end{equation}
find $\bm u_D \in \bm{H_D}$ such that,
\begin{equation}\label{eq:convdiff_curl_MFD}
\begin{split}
{\curl}_V \left( \mathcal{D}_{\alpha, \bm{e}^V} \mathcal{J}^1_D \bm u_D \right) + \mathcal{D}_{\gamma, \bm{e}^D} \bm{u}_D &= \bm f_D \quad \text{in } \Omega,  \\ 
\bm u_D &=  0 \quad \text{on } \partial \Omega ;
 \end{split}
\end{equation}
and find $\bm u_V \in \bm{H_V}$ such that
\begin{equation}\label{eq:convdiff_div_MFD}
\begin{split}
{\grad}_V \left( \mathcal{D}_{\alpha, \bm{x}^V} \mathcal{J}^2_V \bm u_V \right) + \mathcal{D}_{\gamma, \bm{e}^V} \bm{u}_V &= \bm f_V \quad \text{in } \Omega,  \\ 
 \bm u_V  &= 0 \quad \text{on } \partial \Omega,
 \end{split}
\end{equation}
where $\mathcal{D}_{\alpha, \bm{e}^D}$, $\mathcal{D}_{\alpha, \bm{e}^V}$, and $\mathcal{D}_{\alpha, \bm{x}^V}$ are diagonal matrices of suitable piecewise constant approximations of the average value of the diffusion coefficient, $\alpha (\bm{x})$, on the Delaunay edges, Voronoi edges, and Voronoi nodes, respectively.  Similarly, $\mathcal{D}_{\gamma, \bm{x}^D}$, $\mathcal{D}_{\gamma, \bm{e}^D}$, and $\mathcal{D}_{\gamma, \bm{e}^V}$ are diagonal matrices of suitable piecewise constant approximations of the average values of the coefficient $\gamma (\bm{x})$ on the Delaunay nodes, Delaunay edges, and Voronoi edges, respectively. Additionally, $f_D$, $\bm{f}_D$, and $\bm{f}_V$ are all suitable discretizations of the continuous right-hand sides of~\eqref{eq:convdiff_grad}--\eqref{eq:convdiff_div}. 
In particular,
\begin{align}
\mathcal{D}_{\alpha, \bm{e}^D} &= \operatorname{diag}\left( \dashint_{\bm{e}_{ij}^D} \hskip -2pt \alpha \mathrm{d}s \right), & \mathcal{D}_{\alpha, \bm{e}^V} &= \operatorname{diag}\left(\dashint_{\partial D_{km}} \hskip -8pt \alpha \mathrm{d}A  \right), & \mathcal{D}_{\alpha, \bm{x}^V} &= \operatorname{diag}\left( \dashint_{D_k} \hskip -2pt  \alpha \mathrm{d}V \right), \label{eq:Dalpha}\\
\mathcal{D}_{\gamma, \bm{x}^D} &= \operatorname{diag}\left( \gamma(\bm{x}_i^D) \right), & \mathcal{D}_{\gamma, \bm{e}^D} &= \operatorname{diag}\left(\dashint_{\bm{e}_{ij}^D} \hskip -2pt  \gamma \mathrm{d}s  \right), & \mathcal{D}_{\gamma, \bm{e}^V} &= \operatorname{diag}\left( \dashint_{\partial D_{km}} \hskip -8pt \gamma \mathrm{d}A \right). \label{eq:Dgamma}\\
f_D &= \operatorname{vec}\left( f(\bm{x}_i^D) \right), & \bm{f}_D &= \operatorname{vec}\left( \dashint_{\bm{e}_{ij}^D} \bm{f} \mathrm{d}s \right), & \bm{f}_V &= \operatorname{vec} \left( \dashint_{\partial D_{km}} \hskip -8pt \bm{f} \mathrm{d}A \right).  \label{eq:discrete-f}
\end{align}

%\casey{ I think the notation is wrong here, too. Should be something like 
%$\mathcal{D}_{\alpha, \bm{e}^D} = \operatorname{diag}\left( \dashint_{\bm{e}_{ij}^D} \alpha \mathrm{d}s \right)$, $\mathcal{D}_{\gamma, \bm{x}^D} = \operatorname{diag}\left( \gamma(\bm{x}_i^D) \right)$. $\bm x^D$ and $\bm e^D$ alone are not things defined by MFD.} \ja{I think I fixed it but maybe we don't need any of this?
\subsubsection{Quadrature Rules} \label{subsec:quad} 
The discrete MFD flux operators~\eqref{eq:J0_MFD}--\eqref{eq:J2_MFD} are defined with the exponential integrals computed exactly. In practice, this is intractable and a quadrature rule is required. Since we have an exponential function in the integrand, a simple numerical integration rule might not guarantee the desired accuracy. Therefore, we consider quadrature rules that are exact when $\bm \theta$ is a constant on the integral domain, which is equivalent to assuming a linear potential function $\varphi$.  Full derivations of these schemes can be found in \cite{casey:thesis}, and we state the results here.

For $E_{\bm{e}^D} = \text{diag} \left(\dashint_{\bm e_{ij}^D} e^{\varphi} \mathrm{d}s \right)$, we use the following edge quadrature rule,
\begin{equation}
	\dashint_{\bm{e}_{ij}^D} e^{\varphi}  \mathrm{d}s \approx \frac{e^{\varphi_j} - e^{\varphi_i}}{(\bm \theta_{ij} \cdot \bm T^D_{ij}  )}, \label{eq:edgequad}
\end{equation}
where $\bm{\theta}_{ij}$ is a constant approximation of $\bm{\theta}$ on the edge $\bm{e}^D_{ij}$ and $\bm T_{ij}^D = \bm t^D_{ij} | \bm{e}_{ij}^D|$.  It is easy to verify that this quadrature rule is exact if $\bm{\theta}$ is constant on the edge $\bm{e}^D_{ij}$.
For $E_{\partial D} = \text{diag} \left(\dashint_{\partial D_{pm}} e^{\varphi} \mathrm{d}A\right)$, let $\bm x_i^D, \bm x_j^D, \bm x_k^D$ be the indices of face $\partial D_{pm}$ where it is assumed that $i < j < k$.  Then, we obtain the following face quadrature rule,
\begin{align}
	\begin{split}
		\dashint_{\partial D_{pm}} \hskip -14pt e^{\varphi} \, \mathrm{d}A \approx 2 \Bigg(  \frac{  e^{\varphi_i}}{( \bm \theta_{ij} \cdot \bm{T}^D_{ij} )(\bm \theta_{ik}\cdot \bm{T}^D_{ik} ) }  -  \frac{ e^{\varphi_j}}{(\bm \theta_{ij}\cdot \bm{T}^D_{ij})(\bm \theta_{jk}\cdot \bm{T}^D_{jk})}
		+  \frac{ e^{\varphi_k}}{(\bm \theta_{jk}\cdot \bm{T}^D_{jk})(\bm \theta_{ik}\cdot \bm{T}^D_{ik})}  \Bigg), \label{eq:facequad}
	\end{split} 
\end{align}
which is exact if $\bm{\theta}$ is a constant on the face $\partial D_{pm}$. 
Lastly, for $E_{D} = \text{diag} \left( \dashint_{D_p} e^{\varphi} \mathrm{d}V \right)$, let $\bm x_i^D, \bm x_j^D, \bm x_k^D, \bm{x}_{\ell}^D$ be the indices of element $\partial D_p$ where it is assumed that $i < j < k < \ell$. Then, we have the following element quadrature rule,
\begin{align}
	\begin{split}
		\dashint_{D_p} e^{\varphi}  \, \mathrm{d}V & \approx \frac{6 }{\bm \theta_{k\ell} \cdot  \bm{T}^D_{k\ell} } \Bigg(  \frac{  e^{\varphi_i}}{( \bm \theta_{ij}\cdot \bm{T}^D_{ij} )(\bm \theta_{i\ell}\cdot \bm{T}^D_{i\ell} ) }  -  \frac{ e^{\varphi_j}}{(\bm \theta_{ij}\cdot \bm{T}_{ij}^D)(\bm \theta_{j\ell}\cdot \bm T^D_{j\ell})} \\
		&+  \frac{ e^{\varphi_\ell}}{(\bm \theta_{j\ell}\cdot \bm T^D_{j\ell})(\bm \theta_{i\ell}\cdot \bm T^D_{i\ell})}  \Bigg)  - \frac{6 }{\bm \theta_{k\ell} \cdot  \bm T^D_{k\ell}} \Bigg(  \frac{  e^{\varphi_i}}{( \bm \theta_{ij}\cdot \bm{T}^D_{ij} )(\bm \theta_{ik}\cdot \bm{T}^D_{ik} ) }  \\
		&-  \frac{ e^{\varphi_j}}{(\bm \theta_{ij}\cdot \bm T^D_{ij})(\bm \theta_{jk}\cdot \bm T^D_{jk})} +  \frac{ e^{\varphi_k}}{(\bm \theta_{jk}\cdot \bm T^D_{jk})(\bm \theta_{ik}\cdot \bm T^D_{ik})}  \Bigg),  \label{eq:elmquad}
	\end{split}
\end{align}
which is exact if $\bm{\theta}$ is a constant on the element.

\begin{remark}
Recall that the goal of this scheme is to yield a stable method for the convection-diffusion equation in the convection-dominated case. In other words, the mimetic scheme should handle the case when the diffusion coefficient tends to zero, $\alpha \to 0$. As $\varphi$ scales with $\alpha^{-1}$, this means that $\varphi \to \pm \infty$. Thus, for proper implementation of $\mathcal{J}^0_D$, $\mathcal{J}^1_D$, and $\mathcal{J}^2_V$, we must consider numerical stability with floating point arithmetic to avoid overflow.  
For example, substituting \eqref{eq:edgequad} into the denominator of~\eqref{eq:J0_MFD} yields
\begin{equation*}
\left( \mathcal{J}^0_D u \right)_{ij}^D = \left[ \left( \frac{\bm \theta_{ij} \cdot \bm T^D_{ij}}{1 - e^{\varphi_i - \varphi_j} }\right) u^D_j - \left( \frac{\bm \theta_{ij} \cdot \bm T^D_{ij}}{e^{\varphi_j - \varphi_i}-1}  \right) u^D_i \right] \frac{\zeta(i,j)}{|\bm{e}_{ij}^D|}.
\end{equation*}
The problematic case is when $\varphi_i \to \varphi_j$, which makes the denominator of both terms zero.  Taking the Taylor expansion of $\varphi$ at $\bm{\theta}_{ij}$, or  the constant approximation of $\bm{\theta} = \grad \varphi$ on the edge $\bm{e}_{ij}^D$, gives $\bm{\theta}_{ij} \cdot \bm{T}_{ij}^D \approx \varphi_j - \varphi_i$. Therefore,
\begin{align*}
&\lim_{\varphi_i \to \varphi_j} \frac{\bm \theta_{ij} \cdot \bm T^D_{ij}}{1 - e^{\varphi_i - \varphi_j}} \approx \lim_{\varphi_i \to \varphi_j} \frac{\varphi_j - \varphi_i}{1 - e^{\varphi_i - \varphi_j}} = 1, \\
& \lim_{\varphi_i \to \varphi_j} \frac{\bm \theta_{ij} \cdot \bm T^D_{ij}}{e^{\varphi_j - \varphi_i} -1} \approx \lim_{\varphi_i \to \varphi_j} \frac{\varphi_j - \varphi_i}{e^{\varphi_j - \varphi_i} -1} = 1.
\end{align*}
These limits imply that when $\varphi_i \to \varphi_j$, 
\begin{equation*}
\left( \mathcal{J}^0_D u \right)_{ij}^D  \approx  \frac{u^D_j - u^D_i}{|\bm{e}_{ij}^D|} \zeta(i,j),
\end{equation*}
which is the standard MFD gradient~\eqref{eq:gradD_comp}.  In the numerical results below, we identify when $|\varphi_i - \varphi_j| \leq 10^{-12}$ and replace the corresponding entries of $\mathcal{J}^0_D$ with the usual MFD gradient scheme. Similar approaches are performed for $\mathcal{J}^1_D$ and $\mathcal{J}^2_V$, and details can be found in~\cite{casey:thesis}. 
\end{remark}

\section{Analysis of the MFD method}\label{sec:MFD-theory}
By drawing connections between the proposed MFD scheme and SAFE, we provide a path to a concise and straightforward analysis of the MFD in the FE framework. In this section, we examine the connections to prove well-posedness and derive error estimates for MFD.  We also show the stability, i.e., the monotonicity, of the MFD discretization~\eqref{eq:convdiff_grad_MFD} for the scalar convection-diffusion equation~\eqref{eq:convdiff_grad}.

%In the $H(\grad)$ case, we consider the linear Lagrange (P1) finite element, in $\bm{H}(\curl)$ the N\'ed\'elec element, and in $\bm{H}(\divg)$ the Raviart--Thomas (RT) element. 

\subsection{Equivalence between MFD and SAFE} \label{subsec:equivalence}
Defining $\mathcal{E}_{\varphi} u :=  e^{\varphi} u$, the weak formulations of the convection-diffusion equations are:  for $k=0,1,2$, find $u \in H_{0}(\mathfrak{D}^k)$, such that
\begin{equation}
	a(u, v) := \langle \alpha J^k  u, \mathfrak{D}^k v \rangle + \langle \gamma u, v \rangle = \langle \alpha \mathcal{E}_{-\varphi} \mathfrak{D}^k  \mathcal{E}_{\varphi} u, \mathfrak{D}^k v \rangle + \langle \gamma u, v \rangle = \langle f, v \rangle,  \label{eq:CD_weak}
\end{equation}
for all $v \in H_{0}(\mathfrak{D}^k)$.  Note that $H_{0}(\mathfrak{D}^k)$ indicates the Sobolev space $H(\mathfrak{D}^k)$ with appropriate homogeneous Dirichlet boundary conditions applied. Then, the FE formulation is: for $k=0,1,2$, find $u_h \in H_{h,0}(\mathfrak{D}^k)$, such that
\begin{equation}
	a(u_h, v_h) =  \langle \alpha J^ku_h, \mathfrak{D}^k v_h \rangle + \langle \gamma u_h, v_h \rangle= \langle \alpha \mathcal{E}_{-\varphi} \mathfrak{D}^k  \mathcal{E}_{\varphi} u_h, \mathfrak{D}^k v_h \rangle + \langle \gamma u_h, v_h \rangle = \langle f, v_h \rangle,  \label{eq:CD_FE}
\end{equation}
for all $v_h \in H_{h,0}(\mathfrak{D}^k)$.
%Implementing~\eqref{eq:CD_FE} is challenging due to the presence of the $\mathcal{E}_{\varphi}$ and $\mathcal{E}_{-\varphi}$ terms. 
The SAFE scheme uses properties of the canonical interpolation, $\Pi^k$ in \eqref{eq:can_interp}, to show that
\begin{equation*}
\Pi^k \mathcal{E}_{\varphi} w_h =  \sum_{s_i^k \in \mathcal{T}} \left(\dashint_{s_i^k} e^{\varphi} \, \mathrm{d}\bm{x}\right) \eta_{i}^k(w_h) \lambda_{i}^k, \qquad \forall \ w_h \in H_{h}(\mathfrak{D}^k),
\end{equation*}
and then defines the so-called \emph{simplex-averaged} operator $\mathcal{H}_{\varphi}^k: H_h(\mathfrak{D}^k) \mapsto H_h(\mathfrak{D}^k)$,
\begin{equation*}
\mathcal{H}^k_{\varphi} w_h = \sum_{s_i^k \in \mathcal{T}} \left(\dashint_{s_i^k} e^{\varphi} \,\mathrm{d}\bm{x}\right)^{-1} \eta_{i}^k(w_h) \lambda_{i}^k, \qquad \forall \ w_h \in H_{h}(\mathfrak{D}^k).
\end{equation*}
%\casey{I know the $H^k$ keeps things consistent with the SAFE paper, but we use $H$ to denote Sobolev spaces. Maybe $\mathcal{H}^k$ would be better? Maybe same with those exponential operators, use $\mathcal{E}_{\varphi}$ since $E$'s are those diagonal matrix scalings? Because $\mathcal{E}_{\varphi}$ is not just a diagonal scaling--it's a continuous level operator. Just a suggestion.} \ja{OK}
With these definitions, $\mathcal{H}^k_{\varphi} = \left( \Pi^k \mathcal{E}_{\varphi} \right)^{-1}$ on $H_h(\mathfrak{D}^k)$, implying that $\mathcal{H}_{\varphi}^k \Pi^k \mathcal{E}_{\varphi} = I$ on $H_h(\mathfrak{D}^k)$ (or $\mathcal{H}_{\varphi}^k \Pi^k \mathcal{E}_{\varphi} \approx I$, more generally).  Then, the SAFE scheme yields the following approximation:
\begin{equation*}
\mathcal{E}_{-\varphi} \mathfrak{D}^k  \mathcal{E}_{\varphi} \approx \left(\mathcal{H}_{\varphi}^{k+1} \Pi^{k+1} \mathcal{E}_{\varphi}\right) \left( \mathcal{E}_{-\varphi} \mathfrak{D}^k  \mathcal{E}_{\varphi}\right) = \mathcal{H}_{\varphi}^{k+1} \Pi^{k+1} \mathfrak{D}^k \mathcal{E}_{\varphi} = \mathcal{H}_{\varphi}^{k+1} \mathfrak{D}^k \Pi^{k} \mathcal{E}_{\varphi},
\end{equation*}
where we use the commutative property, $\Pi^{k+1} \mathfrak{D}^k = \mathfrak{D}^k \Pi^{k}$, in the last step.  Then, the SAFE scheme is: for $k=0,1,2$, find $u_h \in H_{h,0}(\mathfrak{D}^k)$, such that
\begin{equation}
	a_{\text{SAFE}}(u_h, v_h) = \langle \alpha \mathcal{H}_{\varphi}^{k+1} \mathfrak{D}^k \Pi^{k} \mathcal{E}_{\varphi} u_h, \mathfrak{D}^k v_h \rangle + \langle \gamma u_h, v_h \rangle = \langle f, v_h \rangle,  \label{eq:CD_SAFE}
\end{equation}
for all $v_h \in H_{h,0}(\mathfrak{D}^k)$.
\begin{remark}
As outlined in \cite{SAFE}, \eqref{eq:CD_SAFE} requires the assembly of a weighted mass matrix in $H_h(\mathfrak{D}^{k+1})$. Based on the idea proposed in~\cite{EAFE}, a local constant interpolation is constructed to modify the bilinear form to avoid the explicit construction of the weighted mass matrix.  However, the presented MFD scheme more closely resembles \eqref{eq:CD_SAFE} via mass-lumping and left scaling as discussed later in this subsection, so we only consider \eqref{eq:CD_SAFE} here. 
%\casey{This remark is a really confusing way to say that  EAFE avoids assembling a mass matrix by using a local interpolation, but we don't do that because we can modify $a_{SAFE}$ a different way to get our MFD scheme that also does not require the assembly of a mass matrix. $a_{MFD}$ and mass lumping and scaling haven't ever been mentioned before, and the ``so we only consider (4.3)" when we haven't explicitly mentioned anything else to consider is confusing. How about: ``As outlined in \cite{SAFE}, \eqref{eq:CD_SAFE} requires the assembly of a weighted mass matrix in $H_h(\mathfrak{D}^{k+1})$. Based on the idea proposed in~\cite{EAFE}, a local constant interpolation is constructed to modify the bilinear form to avoid the explicit construction of the weighted mass matrix. However, we present an alternative approach for modifying \eqref{eq:CD_SAFE} which results in the MFD scheme and also avoids construction of a weighted mass matrix, making the local interpolation unnecessary.'' Or add this remark to remark 4.2 after the equivalence is made, which basically says the same thing as this with less detail.}
\end{remark}

As in \cite{MFDMax,FeMFD}, we use the method of mass lumping \cite{Nedlump, RTlump, Thomee_MassLump, Brezzi_MassLump} to draw connections to the MFD scheme.  This allows us to acquire a diagonal approximation of the two weighted mass matrices arising from the two 
weighted $L^2$ inner products on $H_h(\mathfrak{D}^k)$ and $H_h(\mathfrak{D}^{k+1})$ in \eqref{eq:CD_SAFE}. 
To handle the weights $\alpha$ and $\gamma$, we use the canonical interpolation $\Pi^k$ first and then apply the usual mass-lumping technique. This leads to the weak formulation of the MFD scheme:
%Use the low-order term $\langle \gamma u_h, v_h \rangle$ as an example, we approximate it by
%\begin{equation*}
%	\langle \gamma u_h, v_h \rangle \approx  \langle \Pi^k (\gamma u_h), v_h \rangle_h, \qquad \forall \ u_h, \ v_h \in H_h(\mathfrak{D}^k),
%\end{equation*}
%where $\langle \cdot, \cdot \rangle_h$ denotes the corresponding mass-lumping.  Now we are ready to define the weak formulation that leads to our proposed MFD scheme: 
for $k=0,1,2$, find $u_h \in H_{h,0}(\mathfrak{D}^k)$, such that
\begin{equation}
	a_{\text{MFD}}(u_h, v_h) = \langle \Pi^{k+1} \left( \alpha \mathcal{H}_{\varphi}^{k+1} \mathfrak{D}^k \Pi^{k} \mathcal{E}_{\varphi} u_h \right), \mathfrak{D}^k v_h \rangle_h + \langle \Pi^k \left( \gamma u_h \right), v_h  \rangle_h = \langle \Pi^k f, v_h \rangle_h,  \label{eq:CD_MFD}
\end{equation}
for all $v_h \in H_{h,0}(\mathfrak{D}^k)$.  Here, $\langle \cdot, \cdot \rangle_h$ denotes the corresponding mass-lumping.

To make the equivalence clear, consider the scalar convection-diffusion equation first ($k=0$) and ignore the boundary conditions for now. In this case, the matrix representation of $a_{\text{SAFE}}$ is
\begin{equation*}
a_{\text{SAFE}}(u_h, v_h) = v_D^T 
\left( G_{\text{FE}}^T \mathcal{M}_{\alpha, 1}E^{-1}_{\bm{e}^D} G_{\text{FE}} E_{\bm{x}^D}    +  \mathcal{M}_{\gamma, 0} \right) u_D
\end{equation*}
where $u_D$ and $v_D$ are vector representations of $u_h$ and $v_h$, respectively, and $\mathcal{M}_{\alpha, 1}$ and $\mathcal{M}_{\gamma, 0}$ are weighted mass matrices in $\bm{H}_h(\mathfrak{D}^1)=\bm{H}_h(\curl)$ and $\bm{H}_h(\mathfrak{D}^0)=H_h(\grad)$, respectively.  
%\casey{Should we use $\mathcal{M}_{\text{weight}, k}$ instead of $\mathcal{M}_{\text{weight}, \text{derivative}}$ to be consistent?} \ja{ok}

Defining the mass-lumped mass matrix of $H_h(\mathfrak{D}^k)$ as $\widetilde{\mathcal{M}}_k$ for $k = 0, 1, 2$, we obtain
\begin{eqnarray*}
\mathcal{M}_{\gamma, 0} \approx \widetilde{\mathcal{M}}_{0} \mathcal{D}_{\gamma, \bm{x}^D} = \mathcal{D}_V \mathcal{D}_{\gamma, \bm{x}^D}, \quad \mathcal{M}_{\alpha, 1} \approx \widetilde{\mathcal{M}}_{1} \mathcal{D}_{\alpha, \bm{e}^D} = \mathcal{D}_{\bm{e}^D} \mathcal{D}_{\partial V} \mathcal{D}_{\alpha, \bm{e}^D},
\end{eqnarray*}
where, $\mathcal{D}_{\gamma, \bm{x}^D}$ and $\mathcal{D}_{\alpha, \bm{e}^D}$ are defined in \eqref{eq:Dgamma} and \eqref{eq:Dalpha}, respectively.  Full derivations are found in \cite{MFDMax,casey:thesis}.  Thus,
\begin{align*}
	a_{\text{MFD}}(u_h, v_h) &= v_D^T
	\left( G_{\text{FE}}^T \mathcal{D}_{\bm{e}^D} \mathcal{D}_{\partial V} \mathcal{D}_{\alpha, \bm{e}^D}  E^{-1}_{\bm{e}^D} G_{\text{FE}} E_{\bm{x}^D}  +  \mathcal{D}_V \mathcal{D}_{\gamma, \bm{x}^D} \right) u_D  \\
	& =  v_D^T 
	\left(\mathcal{G}^T \mathcal{D}_{\partial V} \mathcal{D}_{\alpha, \bm{e}^D}  \mathcal{J}^0_D  +  \mathcal{D}_V \mathcal{D}_{\gamma, \bm{x}^D} \right) u_D.
\end{align*}
Left scaling by $\mathcal{D}_V^{-1}$ and using the fact $\divg_D = \mathcal{D}_V^{-1}\mathcal{G}^T \mathcal{D}_{\partial V}$, the left-hand side of the MFD system~\eqref{eq:convdiff_grad_MFD} is \textit{exactly} recovered.  Noting that $\langle \Pi^0 f, v_h \rangle_h = v_D^T \widetilde{\mathcal{M}}_{0} f_D$ with $f_D$ defined in \eqref{eq:discrete-f}, the same left scaling recovers the right-hand side as well.

Similar arguments extend to the vector convection-diffusion equations (case $k=1$ and $k=2$). The mass-lumped weighted mass matrices are defined as
\begin{eqnarray*}
	\mathcal{M}_{\gamma, 1} \approx \widetilde{\mathcal{M}}_{1} \mathcal{D}_{\gamma, \bm{e}^D} = \mathcal{D}_{\bm{e}^D} \mathcal{D}_{\partial V} \mathcal{D}_{\gamma, \bm{e}^D}, \quad \mathcal{M}_{\alpha, 2} \approx \widetilde{\mathcal{M}}_{2} \mathcal{D}_{\alpha, \bm{e}^V} = \mathcal{D}_{\partial D} \mathcal{D}_{\bm{e}^V} \mathcal{D}_{\alpha, \bm{e}^V},\\
		\mathcal{M}_{\gamma, 2} \approx \widetilde{\mathcal{M}}_{2} \mathcal{D}_{\gamma, \bm{e}^V} = \mathcal{D}_{\partial D} \mathcal{D}_{\bm{e}^V} \mathcal{D}_{\gamma, \bm{e}^V}, \quad \mathcal{M}_{\alpha, 3} \approx \widetilde{\mathcal{M}}_{3} \mathcal{D}_{\alpha, \bm{x}^V} = \mathcal{D}_{D} \mathcal{D}_{\alpha, \bm{x}^V},
\end{eqnarray*}
and lead to the weak formulation of the MFD scheme,~\eqref{eq:CD_MFD}. Then, for $k=1$, applying the left scaling $\mathcal{D}_{\partial V}^{-1} \mathcal{D}_{\bm{e}^D}^{-1}$,  the MFD matrix system~\eqref{eq:convdiff_curl_MFD} is recovered. For $k=2$, left scaling by $\mathcal{D}_{\bm e^V}^{-1} \mathcal{D}_{\partial D}^{-1}$ leads to the MFD matrix system \eqref{eq:convdiff_div_MFD}. 

\begin{remark}
%\casey{As outlined in \cite{SAFE}, \eqref{eq:CD_SAFE} requires the assembly of a weighted mass matrix in $H_h(\mathfrak{D}^{k+1})$. The approach taken by ~\cite{EAFE} is to construct a local constant interpolation to modify the bilinear form and avoid the explicit construction of the weighted mass matrix. However, this step is not necessary in the proposed MFD method. } 
Based on the equivalence between the SAFE and MFD scheme, we see that the MFD scheme~\eqref{eq:CD_MFD} provides another way to avoid assembly of the weighted mass matrix in $H_h(\mathfrak{D}^{k+1})$, since the mass-lumping only uses the geometric information of the meshes. 
\end{remark}

\begin{remark}
Although our discussion ignores the Dirichlet boundary conditions, the connection still holds after applying them properly on both the SAFE and MFD stiffness matrices. 
\end{remark}

\subsection{Well-posedness and Error Analysis}\label{sec:analysis}
%\ja{This needs to be double checked for norm and inner product notation consistency.}
Next, the well-posedness and error analysis of the MFD scheme for the convection-diffusion problems are considered.  Note that the only difference between the weak formulation of the MFD method \eqref{eq:CD_MFD} and the MFD systems \eqref{eq:convdiff_grad_MFD}--\eqref{eq:convdiff_div_MFD} is the left scaling, which does not affect the numerical solutions. 

We first give some useful lemmas that are crucial to our theoretical analysis. First, we define an interpolation $\widetilde{\Pi}^{k}_{\varphi}: H(\mathfrak{D}^k) \mapsto H_h(\mathfrak{D}^k)$ defined as 
%\begin{equation} \label{eq:tilde_Pi_interp}
$\widetilde{\Pi}^{k}_{\varphi} v := \mathcal{H}_{\varphi}^k \Pi^k \mathcal{E}_{\varphi} v$.
%\end{equation}
In \cite{SAFE}, it has been shown that $\widetilde{\Pi}^{k}_{\varphi} v_h = v_h$ for any $v_h \in H_h(\mathfrak{D}^k)$, and, in general, has the following property:
\begin{lemma}[\cite{SAFE}] \label{lem:widetilde-Pi}
For $T \in \mathcal{T}$, if $v \in W^{1,p}(T)$ and $ p > d$, where $d$ is the dimension, then,
\begin{equation*}
\| v - \widetilde{\Pi}^k_{\varphi} v\| \leq C h^{1 + d(\frac{1}{2} - \frac{1}{p})} \left( \sum_{T \in \mathcal{T}} | v |^2_{1,p,T} \right)^{\frac{1}{2}}.
\end{equation*}
For $v \in L^{\infty}(T)$, $T \in \mathcal{T}$, 
\begin{equation*}
\| \widetilde{\Pi}_{\varphi}^k v \| \leq C h^{\frac{d}{2}}  \left( \sum_{T \in \mathcal{T}} \| v \|_{0, \infty, T}^2 \right)^{\frac{1}{2}}. 
\end{equation*}
where the constant $C$ depends on $p$ and the Sobolev embedding number. %, $|\cdot|_{1,p,T}$ is the semi-norm associated with $W^{1,p}(T)$, and $\|\cdot\|_{0,\infty,T}$ is \ja{what is that??, just infinity norm?}
\end{lemma}
\begin{proof}
The proof can be found in Lemma 5.1 and Corollary 5.2 of \cite{SAFE}.	
\end{proof}

Since the connection between the SAFE and MFD methods depends on the mass-lumping, we assume the following mass-lumping error estimates holds \cite{Thomee_MassLump, Brezzi_MassLump}:
%\casey{Seems like this should be cited since we don't prove or justify, but someone commented it out?}
\begin{equation} \label{eq:mass-lump-error}
| \langle w_h, v_h \rangle - \langle w_h, v_h \rangle_h | \leq C h^2 \| w_h \|_{H(\mathfrak{D}^k)}  \| v_h \|_{H(\mathfrak{D}^k)}, \quad \forall \ w_h, \, v_h \in H_h(\mathfrak{D}^k),
\end{equation}
where $C > 0$ is a constant depending on the shape regularity of the mesh.  Based on the mass-lumping error \eqref{eq:mass-lump-error}, we have the following lemma which bounds the difference between the SAFE and MFD methods.
%\casey{before I make changes, all norms without a subscript are $L^2$ inner products and should have subscript 0?}
\begin{lemma} \label{lem:mass-lump-SAFE-MFD}
For any positive weight function, $0< \mu(\bm{x}) \in W^{1, \infty}(\Omega)$, assuming the mass-lumping error estimate, \eqref{eq:mass-lump-error}, holds, then, we have for any $ w_h, \, v_h \in H_h(\mathfrak{D}^k)$,
\begin{equation*}
|\langle \mu w_h, v_h \rangle - \langle  \Pi^k(\mu w_h), v_h \rangle_h | \leq C h \| \mu\|_{1, \infty} \| w_h \|  \| v_h \|_{H(\mathfrak{D}^k)},
\end{equation*}
where the constant $C>0$ does not depend on $\mu$.
%\casey{norms with two subscripts are not defined.}
\end{lemma}
\begin{proof}
For any $ w_h, \, v_h \in H_h(\mathfrak{D}^k)$,
\begin{align}\label{eq:mubound}
\langle \mu w_h, v_h \rangle - \langle  \Pi^k(\mu w_h), v_h \rangle_h &= \langle \mu w_h - \Pi^k(\mu w_h), v_h \rangle  \\\nonumber
& \quad + \left( \langle \Pi^k( \mu w_h), v_h \rangle - \langle  \Pi^k(\mu w_h), v_h \rangle_h \right).
\end{align}
Note that $w_h = \sum_{s_i^k} \eta^k_{i}(w_h) \lambda^k_{i} $ and $\Pi^k(\mu w_h) = \sum_{s_i^k} \left(\dashint_{s_i^k} \mu \right) \eta^k_{i}(w_h) \lambda^k_{i}$.  Then, 
\begin{align*}
& \quad |\langle \mu w_h - \Pi^k(\mu w_h), v_h \rangle | \\
& \leq  \| \mu w_h - \Pi^k (\mu w_h) \|   \| v_h \| = \|  \sum_{s_i^k} \left( \mu - \dashint_{s_i^k} \mu \, \mathrm{d}\bm{x}\right) \eta^k_{i}(w_h) \lambda^k_{i} \| \| v_h \| \\
& \leq \|\mu - \dashint_{s_i^k} \mu\,\mathrm{d}\bm{x}\|_{0,\infty}\|w_h\|\|v_h\| \leq  C h \| \mu \|_{1, \infty} \| w_h \| \| v_h \|
\end{align*}
For the second term on the right-hand side of \eqref{eq:mubound}, apply the mass-lumping error estimates \eqref{eq:mass-lump-error}, an inverse inequality, and the stability of $\Pi^k$, to get
\begin{align*}
| \langle \Pi^k( \mu w_h), v_h \rangle - \langle  \Pi^k(\mu w_h), v_h \rangle_h | &\leq C h^2  \| \Pi^k (\mu w_h) \|_{H(\mathfrak{D}^k)}   \| v_h \|_{H(\mathfrak{D}^k)} \\
&\leq C h \| \Pi^k(\mu w_h) \|  \| v_h \|_{H(\mathfrak{D}^k)} \\
& \leq C h \| \mu \|_{0, \infty} \| w_h \| \| v_h \|_{H(\mathfrak{D}^k)}
\end{align*}
Combining the two results completes the proof.
\end{proof}

Now we are ready to present an important lemma which measures the consistency between the bilinear forms $a(\cdot, \cdot)$ from \eqref{eq:CD_FE} and $a_{\text{MFD}}(\cdot, \cdot)$ from \eqref{eq:CD_MFD}, using the intermediary $a_{\text{SAFE}}(\cdot,\cdot)$ from \eqref{eq:CD_SAFE}.
\begin{lemma} \label{lem:diff-FE-MFD}
Assume $\alpha, \, \gamma \in W^{1, \infty}(\Omega)$, $\bm{\theta} \in W^{0,\infty}(\Omega)$ and the assumptions for the mass-lumping error estimate, \eqref{eq:mass-lump-error}, are satisfied. Then,  for any $w_h, \, v_h \in H_h(\mathfrak{D}^k)$, the following estimate holds,
\begin{equation*}
	|a(w_h, v_h) - a_{\rm{MFD}}(w_h, v_h)| \leq c_p \, h \|  w_h \|_{H(\mathfrak{D}^k)}  \| v_h \|_{H(\mathfrak{D}^k)},
\end{equation*}
where $c_p := C_1 \| \alpha \|_{0, \infty}  \| \bm{\theta} \|_{0, \infty} +  C_3 \| \alpha \|_{1, \infty}  (1 + \| \bm{\theta} \|_{0,\infty}) + C_2 \| \gamma \|_{1, \infty} >0$ with $C_1$, $C_2$, and $C_3$ being generic positive constants that do not depend on $\alpha$, $\bm{\theta}$, and $\gamma$.
\end{lemma}
\begin{proof}
For any $w_h, \, v_h \in H(\mathfrak{D}^k)$, we have
\begin{equation*}
a(w_h, v_h) - a_{\text{MFD}}(w_h, v_h) = a(w_h, v_h) - a_{\text{SAFE}}(w_h, v_h) + a_{\text{SAFE}}(w_h, v_h) - a_{\text{MFD}}(w_h, v_h).
\end{equation*}	
Since $  H^{k+1}_{\varphi} \mathfrak{D}^k \Pi^k \mathcal{E}_{\varphi} = \widetilde{\Pi}_{\varphi}^{k+1} J^{k}$, 
\begin{align*}
|a(w_h, v_h) - a_{\text{SAFE}}(w_h, v_h)| & = \left| \langle \alpha  (I - \widetilde{\Pi}_{\varphi}^{k+1}) J^k w_h, \mathfrak{D}^k v_h \rangle \right|	\\
& \leq \| \alpha \|_{0, \infty} \| (I - \widetilde{\Pi}_{\varphi}^{k+1}) J^k w_h \| \| \mathfrak{D}^k v_h \| \\
& \leq C_1 \| \alpha \|_{0, \infty} \, h^{1 + d(\frac{1}{2} - \frac{1}{p})} \left( \sum_{T \in \mathcal{T}} | J^k w_h |^2_{1,p,T} \right)^{\frac{1}{2}} \| \mathfrak{D}^k v_h \| \\
%& \leq C \| \alpha \|_{0, \infty} \, h^{1 + d(\frac{1}{2} - \frac{1}{p})} \left( \sum_{T \in \mathcal{T}} | \mathcal{D}^k w_h + \bm{\theta} *^k w_h |^2_{1,p,T} \right)^{\frac{1}{2}} \| \mathfrak{D}^k v_h \| \\
& \leq C_1 \| \alpha \|_{0, \infty}  \| \bm{\theta} \|_{0, \infty} \, h \| w_h \|_{H(\mathfrak{D}^k)} \| \mathfrak{D}^k v_h \| 
%\\
%& \leq C \| \alpha \|_{0, \infty}  | \varphi |_{1, \infty} \, h \| w_h \|_{H(\mathfrak{D}^k)} \| \mathfrak{D}^k v_h \|.
\end{align*}
Here, we use $|\mathfrak{D}^k w_h|_{1, p, T} = 0$ for any $w_h \in H_h(\mathfrak{D}^k)$ and an inverse inequality in the last step. Next we estimate the difference between $a_{\text{SAFE}}$ and $a_{\text{MFD}}$, 
\begin{align*}
& | a_{\text{SAFE}}(w_h, v_h) - a_{\text{MFD}}(w_h, v_h)| \\
& \leq \left| \langle \alpha \widetilde{\Pi}_{\varphi}^{k+1} J^k w_h, \mathfrak{D}^k v_h \rangle - \langle \Pi^{k+1} (  \alpha \widetilde{\Pi}_{\varphi}^{k+1} J^k w_h ), \mathfrak{D}^k v_h \rangle_h \right| \\
& \quad + \left|\langle \gamma w_h, v_h \rangle - \langle \Pi^{k}( \gamma w_h), v_h \rangle_h \right| 
\end{align*}
By \cref{lem:mass-lump-SAFE-MFD},  we have
%\begin{equation*}
$\left|\langle \gamma w_h, v_h \rangle - \langle \Pi^{k}( \gamma w_h), v_h \rangle_h \right| \leq C_2 \| \gamma \|_{1, \infty} \, h \| w_h \| \| v_h \|_{H(\mathfrak{D}^k)}$.
%\end{equation*}
Similarly, using \cref{lem:mass-lump-SAFE-MFD}, 
\begin{align*}
& \quad \left| \langle \alpha \widetilde{\Pi}_{\varphi}^{k+1} J^k w_h, \mathfrak{D}^k v_h \rangle - \langle \Pi^{k+1} (  \alpha \widetilde{\Pi}_{\varphi}^{k+1} J^k w_h ), \mathfrak{D}^k v_h \rangle_h \right|	\\
& \leq C \| \alpha \|_{1, \infty} \, h \|  \widetilde{\Pi}_{\varphi}^{k+1} J^k w_h \| \| v_h \|_{H(\mathfrak{D}^k)} 
\end{align*}
Based on \cref{lem:widetilde-Pi} and an inverse inequality, 
\begin{align*}
\|  \widetilde{\Pi}_{\varphi}^{k+1} J^k w_h \|^2 & \leq C h^{d}  \sum_{T \in \mathcal{T}} \| J^k w_h \|_{0,\infty, T}^2 \\
& \leq C \sum_{T \in \mathcal{T}} \left( h^{\frac{d}{2}} \| \mathfrak{D}^k w_h \|_{0,\infty, T} + \| \bm{\theta} \|_{0,\infty, T} \, h^{\frac{d}{2}} \|w_h \|_{0, \infty, T}  \right)^2  \\
& \leq C (1+ \| \bm{\theta} \|_{0, \infty}) \| w_h \|_{_{H(\mathfrak{D}^k)} }.
%& \leq C \sum_{T \in \mathcal{T}} \left( h^{\frac{d}{2}} \| \mathfrak{D}^k w_h \|_{0,\infty, T} + \| \varphi \|_{1,\infty, T} \, h^{\frac{d}{2}} \|w_h \|_{0, \infty, T}  \right)^2  \\
%& \leq C (1+ \| \varphi \|_{1, \infty}) \| w_h \|_{_{H(\mathfrak{D}^k)} }.
\end{align*}
Therefore, 
\begin{align*}
	& \quad \left| \langle \alpha \widetilde{\Pi}_{\varphi}^{k+1} J^k w_h, \mathfrak{D}^k v_h \rangle - \langle \Pi^{k+1} (  \alpha \widetilde{\Pi}_{\varphi}^{k+1} J^k w_h ), \mathfrak{D}^k v_h \rangle_h \right|	\\
	& \leq C_3 \| \alpha \|_{1, \infty}  (1 + \| \bm{\theta} \|_{0,\infty}) \, h \|  w_h \|_{H(\mathfrak{D}^k)}  \| v_h \|_{H(\mathfrak{D}^k)}.
\end{align*}
This implies that 
%\casey{next line doesn't make sense. Are we trying to say that $C$ is a function of those norms of the coeffs? If so, that's confusing.}
\begin{align*}
 & \quad \left|a_{\text{SAFE}}(w_h, v_h) - a_{\text{MFD}}(w_h, v_h) \right| \leq \tilde{C} \, h \|  w_h \|_{H(\mathfrak{D}^k)}  \| v_h \|_{H(\mathfrak{D}^k)},
\end{align*}
where the constant $\tilde{C} := C_3 \| \alpha \|_{1, \infty}  (1 + \| \bm{\theta} \|_{0,\infty}) + C_2 \| \gamma \|_{1, \infty} >0$. 
Overall, %\casey{same comment here.}
\begin{align*}
|a(w_h, v_h) - a_{\text{MFD}}(w_h, v_h)| \leq c_p \, h \|  w_h \|_{H(\mathfrak{D}^k)}  \| v_h \|_{H(\mathfrak{D}^k)},
\end{align*}
which completes the proof with $c_p := C_1 \| \alpha \|_{0, \infty}  \| \bm{\theta} \|_{0, \infty} +  C_3 \| \alpha \|_{1, \infty}  (1 + \| \bm{\theta} \|_{0,\infty}) + C_2 \| \gamma \|_{1, \infty} >0$.
\end{proof}

\subsubsection{Well-posedness}
We assume the continuous problem~\eqref{eq:CD_weak} is well-posed, i.e.,
\begin{equation} \label{eq:bounded-cont}
	\sup_{ w \in H_0(\mathfrak{D}^k)} \sup_{v \in H_0(\mathfrak{D}^k)} \frac{|a(w,v)|}{\| w\|_{H(\mathfrak{D}^k)} \| v \|_{H(\mathfrak{D}^k)} } = c_{b} > 0,
\end{equation}
and
\begin{equation} \label{eq:inf-sup-cont}
	\inf_{ w \in H_0(\mathfrak{D}^k)} \sup_{v \in H_0(\mathfrak{D}^k)} \frac{|a(w,v)|}{\| w\|_{H(\mathfrak{D}^k)} \| v \|_{H(\mathfrak{D}^k)} } = c_{i} > 0,
\end{equation}
where the constants $c_{b}$ and $c_{i}$ may depend on $\alpha$, $\bm{\beta}$, and $\gamma$.

The main tool to show the well-posedness of the weak-formulation of the MFD scheme, \eqref{eq:CD_MFD}, is \cref{lem:diff-FE-MFD}. First, we show that the MFD bilinear form, $a_{\text{MFD}}$, is bounded, and then that it satisfies an inf-sup condition.
\begin{lemma}[Boundedness of $a_{\text{MFD}}$] \label{lem:bound-MFD}
Under the same assumptions of \cref{lem:diff-FE-MFD} and assuming \eqref{eq:bounded-cont}, 
\begin{equation} \label{eq:bounded-MFD}
	\sup_{ w_h \in H_{h,0}(\mathfrak{D}^k)} \sup_{v_h \in H_{h,0}(\mathfrak{D}^k)} \frac{|a_{\rm{MFD}}(w_h,v_h)|}{\| w_h\|_{H(\mathfrak{D}^k)} \| v_h \|_{H(\mathfrak{D}^k)} } = \tilde{c}_{b} > 0,
\end{equation}
\end{lemma}
\begin{proof}
Using \eqref{eq:bounded-cont} and \cref{lem:diff-FE-MFD}, 
\begin{align*}
|a_{\text{MFD}}(w_h, v_h)| & \leq |a(w_h, v_h)| + |a_{\text{MFD}}(w_h, v_h) - a(w_h, v_h)|	 \\
& \leq c_b \| w_h\|_{H(\mathfrak{D}^k)} \| v_h \|_{H(\mathfrak{D}^k)} + c_p h \| w_h\|_{H(\mathfrak{D}^k)} \| v_h \|_{H(\mathfrak{D}^k)} \\
& = \tilde{c}_b \| w_h\|_{H(\mathfrak{D}^k)} \| v_h \|_{H(\mathfrak{D}^k)},
\end{align*}
which implies \eqref{eq:bounded-MFD} with $\tilde{c}_b = c_b + c_ph$.
\end{proof}
 
\begin{lemma}[Inf-sup condition of $a_{\text{MFD}}$] \label{lem:inf-sup-MFD}
Under the same assumptions of Lemma \ref{lem:diff-FE-MFD} and assuming \eqref{eq:inf-sup-cont}, for sufficiently small $h$, 
\begin{equation} \label{eq:inf-sup-MFD}
	\inf_{ w_h \in H_{h,0}(\mathfrak{D}^k)} \sup_{v_h \in H_{h,0}(\mathfrak{D}^k)} \frac{|a_{\rm{MFD}}(w_h,v_h)|}{\| w_h\|_{H(\mathfrak{D}^k)} \| v_h \|_{H(\mathfrak{D}^k)} } = \tilde{c}_{i} > 0.
\end{equation}	
\end{lemma}
\begin{proof}
According to the continuous inf-sup condition, \eqref{eq:inf-sup-cont}, for any given $v_h \in H_{h,0}(\mathfrak{D}^k)$, there exists $w_h \in H_{h,0}(\mathfrak{D}^k)$ such that
%\begin{equation*}
$	\frac{|a(w_h, v_h)|}{\| v_h \|_{H(\mathfrak{D}^k)}} \geq c_i \| w_h \|_{H(\mathfrak{D}^k)}$.
%\end{equation*}
Using \cref{lem:diff-FE-MFD}, 
\begin{align*}
\frac{|a_{\text{MFD}}(w_h, v_h)|}{\| v_h \|_{H(\mathfrak{D}^k)}} & \geq  	\frac{|a(w_h, v_h)|}{\| v_h \|_{H(\mathfrak{D}^k)}} -  	\frac{|a_{\text{MFD}}(w_h, v_h) - a(w_h, v_h)|}{\| v_h \|_{H(\mathfrak{D}^k)}} \\
& \geq   c_i \| w_h \|_{H(\mathfrak{D}^k)} - c_p h \| w_h \|_{H(\mathfrak{D}^k)}  = \tilde{c}_i \| w_h \|_{H(\mathfrak{D}^k)},
\end{align*}
where $\tilde{c}_i = c_i - c_ph > 0$ for sufficiently small $h$. This completes the proof.   
\end{proof}
Thus, we have the well-posedness of the MFD weak formulation \eqref{eq:CD_MFD} as follows.
\begin{theorem}[Well-posedness of the MFD weak formulation \eqref{eq:CD_MFD}] \label{thm:well-posed-MFD}
Under the same assumptions of \cref{lem:diff-FE-MFD} and assuming that the bilinear form $a(\cdot,\cdot)$ from \eqref{eq:CD_FE} satisfies the boundedness condition \eqref{eq:bounded-cont} and the inf-sup condition \eqref{eq:inf-sup-cont}, for sufficiently small $h$, the MFD weak formulation, \eqref{eq:CD_MFD}, is well-posed, i.e., it has a unique solution.  
\end{theorem}
\begin{proof}
The result follows from \cref{lem:bound-MFD}, 	~\cref{lem:inf-sup-MFD}, and Babu\v{s}ka's theory \cite{babuvska1971error}.
\end{proof}

\subsubsection{Error Analysis}
With the above well-posedness result, the error analysis of the MFD scheme, \eqref{eq:CD_MFD}, follows from Strang's first lemma \cite{strang1972variational}. 

\begin{theorem}\label{thm:MFD-error}
Let $u$ and $u_h$ be the solutions of \eqref{eq:CD_weak} and \eqref{eq:CD_MFD}, respectively, and $f \in H(\mathfrak{D}^k)$. Under the same assumptions of \cref{lem:diff-FE-MFD} and assuming that the bilinear form $a(\cdot, \cdot)$ satisfies the boundedness condition, \eqref{eq:bounded-cont}, and the inf-sup condition, \eqref{eq:inf-sup-cont}, for sufficiently small $h$,
\begin{align} \label{eq:MFD-error}
& \quad \| u - u_h \|_{H(\mathfrak{D}^k)} \leq\\\nonumber  
&\inf_{w_h \in H_{h,0}(\mathfrak{D}^k)} \left[ \left( 1+ \frac{c_b}{\tilde{c}_i} \right) \| u - w_h \|_{H(\mathfrak{D}^k)} +  \frac{c_p}{\tilde{c}_i} h \| w_h \|_{H(\mathfrak{D}^k)} \right] + Ch \| f \|_{H(\mathfrak{D}^k)},
\end{align}	
where constant $c_b$ is the continues upper bound from the assumption~\eqref{eq:bounded-cont}, $\tilde{c}_i$ is the constant from the discrete inf-sup condition~\eqref{eq:inf-sup-MFD}, $c_p$ is the constant defined in \cref{lem:diff-FE-MFD}, and the last constant $C$ is a genetic constant that does not depend on $\alpha$, $\bm{\theta}$, and $\gamma$.
\end{theorem}
\begin{proof}
For any $w_h \in H_{h,0}(\mathfrak{D}^k)$, we have
\begin{align*}
& \quad  a_{\text{MFD}}(u_h - w_h, v_h) =%&= a_{\text{MFD}} \pm a(u-w_h, v_h) 
\\
& a(u-w_h, v_h) + \left[ a(w_h, v_h) - a_{\text{MFD}}(w_h, v_h)  \right] + \left[ \langle \Pi^kf, v_h \rangle_h - \langle f, v_h \rangle \right] \\
& \leq c_b \| u -w_h \|_{H(\mathfrak{D}^k)} \| v_h \|_{H(\mathfrak{D}^k)} + c_p \, h \| w_h \|_{H(\mathfrak{D}^k)} \| v_h \|_{H(\mathfrak{D}^k)} + C h \| f \|_{H(\mathfrak{D}^k)} \| v_h \|_{H(\mathfrak{D}^k)}.
\end{align*}
By the discrete inf-sup condition \eqref{eq:inf-sup-MFD}, we arrive at
\begin{align*}
\| u_h - w_h \|_{H(\mathfrak{D}^k)} &\leq \tilde{c}_i^{-1} \sup_{v_h \in H_{h,0}(\mathfrak{D}^k)} \frac{|a_{\text{MFD}}(u_h - w_h, v_h)|}{\| v_h \|_{H(\mathfrak{D}^k)}} \\
& \leq \frac{c_b}{\tilde{c}_i} \| u - w_h \|_{H(\mathfrak{D}^k)} + \frac{c_p}{\tilde{c}_i}  h \| w_h \|_{H(\mathfrak{D}^k)} + Ch \| f \|_{H(\mathfrak{D}^k)}.
\end{align*}
Then the error estimates \eqref{eq:MFD-error} follows from the above inequality and the triangle inequality $\| u - u_h \|_{H(\mathfrak{D}^k)} \leq \| u - w_h \|_{H(\mathfrak{D}^k)}  + \| u_h - w_h \|_{H(\mathfrak{D}^k)} $. 
%\ja{I feel like we either show it or we don't...}
\end{proof}
By choosing $w_h = \Pi^k u$ in \cref{thm:MFD-error}, we have the following error estimates.
\begin{corollary} \label{coro:MFD-error}
Let $u$ and $u_h$ be the solutions of \eqref{eq:CD_weak} and \eqref{eq:CD_MFD}, respectively. Furthermore, let $u \in H^2(\Omega)$ %\casey{not defined, nor is $\| \cdot \|_2$. Confusing to define because $\mathfrak{D}^2 = \divg$. } 
and $f \in H(\mathfrak{D}^k)$. Under the same assumptions of \cref{thm:MFD-error},	
\begin{equation*}
\| u - u_h \|_{H(\mathfrak{D}^k)} \leq C \, h \left( \| u \|_{2} + \| f \|_{H(\mathfrak{D}^k)} \right)
\end{equation*}
where $C$ is a positive constant depending on $\alpha$, $\bm{\theta}$, and $\gamma$. 
\end{corollary}

\begin{remark}
If we assume $h$ is small enough so that $\tilde{c}_i \leq c_i/2$, then the constant $C$ in Corollary~\ref{coro:MFD-error} is proportional to $\max\{\dfrac{c_b}{c_i}, \dfrac{c_p}{c_i}\}$, which implies that the above error estimate depends on $\| \bm{\theta} \|_{0, \infty} = \| \frac{\bm{\beta}}{\alpha} \|_{0, \infty}$. As expected for convection-dominated problems, this implies that when the ratio $|\alpha|/|\bm{\beta}|$ gets smaller, the error gets bigger, and a sufficiently small $h$ is needed to observe a proper convergence rate. In~\Cref{sec:num_CD}, our numerical experiments confirm this result. 
\end{remark}

\begin{remark}
The assumption that $f \in H(\mathfrak{D}^k)$ makes the error estimate sub-optimal in terms of the regularity requirement.  This is due to the right-hand side, $\langle \Pi^kf, v_h \rangle_h$, used in the MFD scheme.  If $\langle f, v_h \rangle$ is used, the assumption on $f \in H(\mathfrak{D}^k)$ is no longer needed, but does require that the right-hand sides of \eqref{eq:convdiff_grad_MFD}--\eqref{eq:convdiff_div_MFD} be implemented accordingly. 
\end{remark}

\subsection{Monotonicity for $H(\grad)$}  \label{sec:stability}
Finally, we consider the monotonicity and stability of the scalar convection-diffusion problem, \eqref{eq:convdiff_grad}. Defining $\mathcal{L}u:= -\divg \left( \alpha J^0 u \right) + \gamma u$, we note that $\mathcal{L}$ satisfies a \emph{maximum principle} or \emph{monotonicity property}, which says the inverse of $\mathcal{L}$ is nonnegative. More precisely, if $(\mathcal{L}u)(\bm{x}) > 0$ for all $\bm{x} \in \Omega$, then $u(\bm{x}) \geq 0$ for all $\bm{x} \in \Omega$.  This property holds regardless of the ratio $|\bm{\beta}(\bm{x})|/|\alpha(\bm{x})|$ and, therefore, is important when convection is dominating diffusion, i.e.,  $|\alpha(\bm{x})| \ll |\bm{\beta}(\bm{x})|$.  

In practice, it is crucial to construct a numerical scheme that preserves the monotonicity property discretely. In particular, if $\mathcal{L}_h$ and $f_h$ are discretizations of $\mathcal{L}$ and $f$, respectively, then $\mathcal{L}_h$ has the monotonicity property, if $f_h \geq 0$, implies $\mathcal{L}_h^{-1} f_h \geq 0$. Here, the operand, $\geq$, is understood in a component-wise fashion.  It is also well-known from a linear algebra point of view, that if $\mathcal{L}_h$ is a non-singular M-matrix, then it is monotone (see e.g. \cite{MMatrix}).  

The MFD scheme for the scalar convection-diffusion problem, \eqref{eq:convdiff_grad_MFD} yields
\begin{equation*}
	\mathcal{L}_h =  {\divg}_D \mathcal{D}_{\alpha, \bm{e}^D} \mathcal{J}^0_D + \mathcal{D}_{\gamma, \bm{x}^D} := \mathcal{D}_V^{-1} \mathcal{G}^T  \mathcal{W} \mathcal{G} E_{\bm{x}^D} +  \mathcal{D}_{\gamma, \bm{x}^D},
\end{equation*}
where $\mathcal{W}$ is a diagonal matrix defined as $\mathcal{W} := \mathcal{D}_{\partial V} \mathcal{D}_{\alpha, \bm{e}^D} E_{\bm{e}^D}^{-1} \mathcal{D}_{\bm{e}^D}$. Here, we assume the Dirichlet boundary conditions have been applied so that the rows and columns corresponding to the boundary nodes have been eliminated. The next theorem states that $\mathcal{L}_h$ is a non-singular M-matrix, which implies the monotonicity of the MFD scheme, \eqref{eq:convdiff_grad_MFD}, for the scalar convection-diffusion problem \eqref{eq:convdiff_grad}.

\begin{theorem} \label{thm:M-matrix-MFD}
Let $\alpha(\bm{x}) > 0$ and $\gamma(\bm{x}) \geq 0$ for all $\bm{x} \in \Omega$. If a non-degenerate dual mesh configuration (e.g., non-degenerate Delaunay and Voronoi grids) is used, then the stiffness matrix $\mathcal{L}_h$ of the MFD scheme, \eqref{eq:convdiff_grad_MFD}, is a non-singular M-matrix, i.e., it is monotone. 
%\casey{We need to be careful about the mesh requirement. Delaunay is fine, but I don't think Voronoi is correct. We need a non-degenerate Voronoi mesh, or just a non-degenerate dual mesh. For example, the system is singular on a Delaunay right triangle mesh with the standard Voronoi mesh because the Voronoi mesh is degenerate. And if there are obtuse angles in the Delaunay triangulation, then the Voronoi point lies outside of its primal element and we'll have negative edge lengths, face areas, element volumes. So I think this needs to be changed to say a non-degenerate (every dual point lies in the interior of its corresponding Delaunay element) dual to make this true. I don't even really see why, for MFD specifically, the primal mesh has to be Delaunay. This is a requirement for the FEM/SAFE framework to hold, though. Alternatively, a Delaunay triangulation with all angles less than 90 degrees with the Voronoi diagram also works. But I think the first is better because we design a non-Voronoi dual mesh for the numerical experiments.}
\end{theorem}
\begin{proof}
Based on the non-degenerative assumption of the dual mesh configuration, the diagonal entries of the diagonal matrices, $\mathcal{D}_V$,  $\mathcal{D}_{\partial V}$, $\mathcal{D}_{\bm{e}^D}$, $E_{\bm{x}^D}$, and $E_{\bm{e}^D}$ are all positive. In addition, according to the assumptions on $\alpha$ and $\gamma$, the diagonal entries of $\mathcal{D}_{\alpha, \bm{e}^D}$ are all positive and the diagonal entries of $\mathcal{D}_{\gamma, \bm{x}^D}$ are all non-negative. Therefore, the diagonal entries of the diagonal matrix $\mathcal{W}$ are non-negative. Since $\mathcal{G}$ is the Delaunay edge-vertex signed incidence matrix (with the columns corresponding to the boundary nodes eliminated), the off-diagonal entries of $\mathcal{L}_h$ are
\begin{align*}
\left(\mathcal{L}_h\right)_{ij} = (\mathcal{D}_V)_{ii}^{-1} \, \left( \mathcal{G}^T \mathcal{W} \mathcal{G} \right)_{ij} (E_{\bm{x}^D})_{jj} = -(\mathcal{D}_V)_{ii}^{-1} \, \left( \mathcal{W}\right)_{\bm{e}^D_{ij}} (E_{\bm{x}^D})_{jj} \leq 0,
\end{align*}
where $\left(\mathcal{W}\right)_{\bm{e}^D_{ij}}$ denotes the diagonal entry of $\mathcal{W}$ corresponding to the edge $\bm{e}^D_{ij}$, which incidents with nodes $\bm{x}^D_i$ and $\bm{x}^D_j$. Therefore, $\mathcal{L}_h$ is a Z-matrix (a matrix with all non-positive off-diagonal entries).

Defining the inner product $(u_D, v_D)_{E_{\bm{x}^D} \mathcal{D}_V}:= v_D^T E_{\bm{x}^D} \mathcal{D}_V u_D$, we have 
%\casey{possibly not true for a degenerative mesh, $\mathcal{D}_V$ could perhaps have zeros or negatives, and this won't be an inner product.}
\begin{equation*}
(\mathcal{L}_h u_D, v_D)_{E_{\bm{x}^D} \mathcal{D}_V} = v_D^T E_{\bm{x}^D} \mathcal{G}^T \mathcal{W} \mathcal{G} E_{\bm{x}^D} u_D + v_D^T E_{\bm{x}^D} \mathcal{D}_V \mathcal{D}_{\gamma, \bm{x}^D} u_D.
\end{equation*}
Therefore, $\mathcal{L}_h$ is symmetric with respect to $(\cdot, \cdot)_{E_{\bm{x}^D} \mathcal{D}_V}$ and $(\mathcal{L}_h v_D, v_D)_{E_{\bm{x}^D} \mathcal{D}_V} \geq 0$ for any $v_D$. Furthermore, due to the well-posedness result \cref{thm:well-posed-MFD}, we have $(\mathcal{L}_h v_D, v_D)_{E_{\bm{x}^D} \mathcal{D}_V} > 0$ for any $v_D \neq 0$.  By the Courant minmax principle, all the eigenvalues of $\mathcal{L}_h$ are real and positive. Since $\mathcal{L}_h$ is also a Z-matrix, by one characterization given in \cite{MMatrix}, $\mathcal{L}_h$ is a non-singular M-matrix. 
\end{proof}

\begin{remark}
The non-degenerate Delaunay and Voronoi grid configuration is one where every dual point lies in the interior of the primal element, guaranteeing positivity of mesh metrics. However, the assumption on the dual mesh configuration, specifically requiring a Delaunay primal mesh and Voronoi dual mesh, can be relaxed. Given a primal mesh, which implies $\mathcal{D}_{\bm{e}^D}$ has positive diagonal entries, as long as the dual mesh is reasonable in the sense that diagonal entries of $\mathcal{D}_{\partial V}$ and $\mathcal{D}_V$ are all positive, the resulting MFD scheme is monotone.
\end{remark}

\begin{remark}
\cref{thm:M-matrix-MFD}	assumes that the integrals used to define $\mathcal{D}_{\alpha, \bm{e}^D}$ and $E_{\bm{e}^D}$ are computed exactly. In practice, numerical quadrature rules are used instead, but $\mathcal{L}_h$ is still a non-singular M-matrix if the numerical quadrature rules preserve positivity.  For computing $\mathcal{D}_{\alpha, \bm{e}^D}$, standard numerical quadrature rules, such as the midpoint rule, preserve the positivity.  For computing $E_{\bm{e}^D}$, special attention is needed. However, it can be shown that the quadrature rule \eqref{eq:edgequad} introduced in \cref{subsec:quad} does preserve the positivity.  Therefore, our implementation provides a monotone MFD discretization for the scalar convection-diffusion problem. 
\end{remark}

\begin{remark}
For the case of pure Neumann boundary condition, $\mathcal{L}_h$ is singular and its kernel is spanned by constants.  In this case, if we restrict ourselves to a subspace that is orthogonal to that kernel,  then $\mathcal{L}_h^{-1}$ is still well-defined and positive. Therefore, it is monotone. 
\end{remark}

\begin{remark}
The stability for the $\bm{H}(\curl)$ and $\bm{H}(\divg)$ cases are more subtle than the $H(\grad)$ case considered here. The stability cannot be simply explained or verified by the monotonicity of the corresponding stiffness matrices, and, in fact, the stiffness matrices are not monotone for both cases anymore. However, we still observe the stability of the MFD schemes numerically; see \Cref{sec:num_CD} for details.
\end{remark}

%\ja{I feel like we need some remark on the vector versions of the equations, no?}

\section{Numerical Results}\label{sec:num_CD}
\begin{figure}[h!]
  \centering
\begin{tikzpicture}[scale = 5., every node/.style={scale = 1}]
\draw[thick] (.5,0) to (.25, .433) to (-.25, .433) to (-.5, 0) to (-.25, -.433) to (.25, -.433) to cycle;
\draw[thick] (.25, .433) to (-.25, -.433);
\draw[thick] (-.25, .433) to (.25, -.433);
\draw[thick] (.5,0) to (-.5, 0);
\draw[dashed, red, thick] (.25, .1443) to (.25, -.1433) to (0, -.2887) to (-.25, -.1443) to (-.25, .1443) to (0, .2887) to cycle;
\draw[dashed, red, thick] (0, .2887) to (0, .433);
\draw[dashed, red, thick] (0, -.2887) to (0, -.433);
\draw[dashed, red, thick] (-.25, .1433) to (-0.375, 0.2165);
\draw[dashed, red, thick] (.25, .1433) to (0.375, 0.2165);
\draw[dashed, red, thick] (-.25, -.1433) to (-0.375, -0.2165);
\draw[dashed, red, thick] (.25, -.1433) to (0.375, -0.2165);

\fill (0,0) circle (.5pt);
\fill (.5,0) circle (.5pt);
\fill (.25, .433) circle (.5pt);
\fill (-.25, .433) circle (.5pt);
\fill (-.5, 0) circle (.5pt);
\fill (-.25, -.433) circle (.5pt);
\fill (.25, -.433) circle (.5pt);
\fill[red] (.25, .1443) circle (.5pt);
\fill[red] (.25, -.1433) circle (.5pt);
\fill[red] (0, -.2887) circle (.5pt);
\fill[red] (-.25, -.1443) circle (.5pt);
\fill[red] (-.25, .1443) circle (.5pt);
\fill[red] (0, .2887) circle (.5pt);

%del nodes
\node at (0.01, .11) {$\bm{x}_1^D$};
\node at (.57, 0) {$\bm{x}_2^D$};
\node at (-.57, 0) {$\bm{x}_5^D$};
\node at (.25, .5){$\bm{x}_7^D$};
\node at (-.25, .5){$\bm{x}_6^D$};
\node at (-.25, - .5){$\bm{x}_4^D$};
\node at (.25, -.5){$\bm{x}_3^D$};

%%vor nodes
%\node at (.25, -.2) {$\bm{x}_1^V$};
%\node at ( .07, -.3){$\bm{x}_2^V$};
%\node at (-.25, -.2) {$\bm{x}_3^V$};
%\node at (-.25, .2) {$\bm{x}_4^V$};
%\node at ( .07, .3){$\bm{x}_5^V$};
%\node at (.25, .2) {$\bm{x}_6^V$};

%every node/.style={scale = 2}
\end{tikzpicture}  
\caption{Hexagonal computational domain with $h=1$, centered at $(0,0)$. The primal mesh is drawn with solid black lines and the dual mesh is shown in red dashed lines.} \label{fig:hexagon_mesh}
\end{figure}
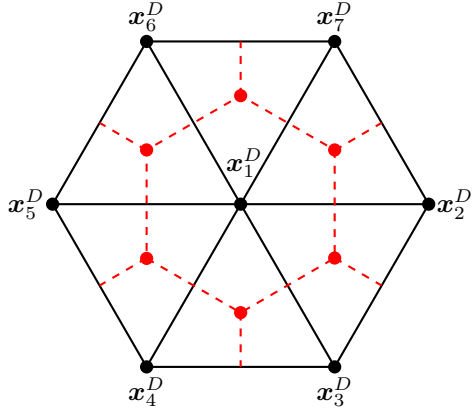

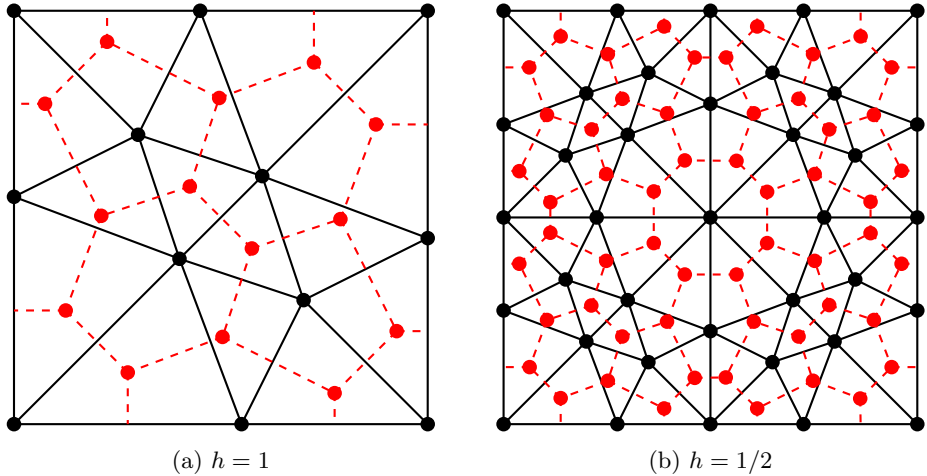
\begin{figure}
\centering
\begin{subfigure}{.5\textwidth}
  \centering
\begin{tikzpicture}[xscale=5.5,yscale=5.5, every node/.style={scale = .8}]
\fill (0.000, 0.000) circle (.5pt);
\fill (1.000, 0.000) circle (.5pt);
\fill (0.000, 1.000) circle (.5pt);
\fill (1.000, 1.000) circle (.5pt);
\fill (0.550, 0.000) circle (.5pt);
\fill (0.000, 0.550) circle (.5pt);
\fill (1.000, 0.450) circle (.5pt);
\fill (0.450, 1.000) circle (.5pt);
\fill (0.300, 0.700) circle (.5pt);
\fill (0.700, 0.300) circle (.5pt);
\fill (0.600, 0.600) circle (.5pt);
\fill (0.400, 0.400) circle (.5pt);
\draw[thick] (0.000,0.000) to (0.550, 0.000);
\draw[thick] (0.000,0.000) to (0.000, 0.550);
\draw[thick] (0.000,0.000) to (0.400, 0.400);
\draw[thick] (1.000,0.000) to (0.550, 0.000);
\draw[thick] (1.000,0.000) to (1.000, 0.450);
\draw[thick] (1.000,0.000) to (0.700, 0.300);
\draw[thick] (0.000,1.000) to (0.000, 0.550);
\draw[thick] (0.000,1.000) to (0.450, 1.000);
\draw[thick] (0.000,1.000) to (0.300, 0.700);
\draw[thick] (1.000,1.000) to (1.000, 0.450);
\draw[thick] (1.000,1.000) to (0.450, 1.000);
\draw[thick] (1.000,1.000) to (0.600, 0.600);
\draw[thick] (0.550,0.000) to (0.700, 0.300);
\draw[thick] (0.550,0.000) to (0.400, 0.400);
\draw[thick] (0.000,0.550) to (0.300, 0.700);
\draw[thick] (0.000,0.550) to (0.400, 0.400);
\draw[thick] (1.000,0.450) to (0.700, 0.300);
\draw[thick] (1.000,0.450) to (0.600, 0.600);
\draw[thick] (0.450,1.000) to (0.300, 0.700);
\draw[thick] (0.450,1.000) to (0.600, 0.600);
\draw[thick] (0.300,0.700) to (0.600, 0.600);
\draw[thick] (0.300,0.700) to (0.400, 0.400);
\draw[thick] (0.700,0.300) to (0.600, 0.600);
\draw[thick] (0.700,0.300) to (0.400, 0.400);
\draw[thick] (0.600,0.600) to (0.400, 0.400);
\fill[red] (0.211, 0.504) circle (.5pt);
\fill[red] (0.125, 0.275) circle (.5pt);
\fill[red] (0.075, 0.775) circle (.5pt);
\fill[red] (0.496, 0.789) circle (.5pt);
\fill[red] (0.225, 0.925) circle (.5pt);
\fill[red] (0.425, 0.575) circle (.5pt);
\fill[red] (0.504, 0.211) circle (.5pt);
\fill[red] (0.875, 0.725) circle (.5pt);
\fill[red] (0.275, 0.125) circle (.5pt);
\fill[red] (0.725, 0.875) circle (.5pt);
\fill[red] (0.775, 0.075) circle (.5pt);
\fill[red] (0.789, 0.496) circle (.5pt);
\fill[red] (0.575, 0.425) circle (.5pt);
\fill[red] (0.925, 0.225) circle (.5pt);
\draw[dashed, red, thick] (0.275, 0.125) to (0.275, 0.000);
\draw[dashed, red, thick] (0.125, 0.275) to (0.000, 0.275);
\draw[dashed, red, thick] (0.125,0.275) to (0.275, 0.125);
\draw[dashed, red, thick] (0.775, 0.075) to (0.775, 0.000);
\draw[dashed, red, thick] (0.925, 0.225) to (1.000, 0.225);
\draw[dashed, red, thick] (0.775,0.075) to (0.925, 0.225);
\draw[dashed, red, thick] (0.075, 0.775) to (0.000, 0.775);
\draw[dashed, red, thick] (0.225, 0.925) to (0.225, 1.000);
\draw[dashed, red, thick] (0.075,0.775) to (0.225, 0.925);
\draw[dashed, red, thick] (0.875, 0.725) to (1.000, 0.725);
\draw[dashed, red, thick] (0.725, 0.875) to (0.725, 1.000);
\draw[dashed, red, thick] (0.875,0.725) to (0.725, 0.875);
\draw[dashed, red, thick] (0.504,0.211) to (0.775, 0.075);
\draw[dashed, red, thick] (0.504,0.211) to (0.275, 0.125);
\draw[dashed, red, thick] (0.211,0.504) to (0.075, 0.775);
\draw[dashed, red, thick] (0.211,0.504) to (0.125, 0.275);
\draw[dashed, red, thick] (0.789,0.496) to (0.925, 0.225);
\draw[dashed, red, thick] (0.789,0.496) to (0.875, 0.725);
\draw[dashed, red, thick] (0.496,0.789) to (0.225, 0.925);
\draw[dashed, red, thick] (0.496,0.789) to (0.725, 0.875);
\draw[dashed, red, thick] (0.425,0.575) to (0.496, 0.789);
\draw[dashed, red, thick] (0.425,0.575) to (0.211, 0.504);
\draw[dashed, red, thick] (0.575,0.425) to (0.789, 0.496);
\draw[dashed, red, thick] (0.575,0.425) to (0.504, 0.211);
\draw[dashed, red, thick] (0.425,0.575) to (0.575, 0.425);
\end{tikzpicture}
  \caption{$h=1$}
  \label{fig:square_mesh1}
\end{subfigure}%
\begin{subfigure}{.5\textwidth}
  \centering
\begin{tikzpicture}[xscale=5.5,yscale=5.5, every node/.style={scale = .8}]
\fill (0.000, 0.000) circle (.5pt);
\fill (0.500, 0.000) circle (.5pt);
\fill (1.000, 0.000) circle (.5pt);
\fill (0.000, 0.500) circle (.5pt);
\fill (0.500, 0.500) circle (.5pt);
\fill (1.000, 0.500) circle (.5pt);
\fill (0.000, 1.000) circle (.5pt);
\fill (0.500, 1.000) circle (.5pt);
\fill (1.000, 1.000) circle (.5pt);
\fill (0.275, 0.000) circle (.5pt);
\fill (0.725, 0.000) circle (.5pt);
\fill (1.000, 0.275) circle (.5pt);
\fill (1.000, 0.725) circle (.5pt);
\fill (0.000, 0.275) circle (.5pt);
\fill (0.500, 0.225) circle (.5pt);
\fill (0.225, 0.500) circle (.5pt);
\fill (0.150, 0.350) circle (.5pt);
\fill (0.350, 0.150) circle (.5pt);
\fill (0.300, 0.300) circle (.5pt);
\fill (0.200, 0.200) circle (.5pt);
\fill (0.775, 0.500) circle (.5pt);
\fill (0.650, 0.150) circle (.5pt);
\fill (0.850, 0.350) circle (.5pt);
\fill (0.800, 0.200) circle (.5pt);
\fill (0.700, 0.300) circle (.5pt);
\fill (0.000, 0.725) circle (.5pt);
\fill (0.500, 0.775) circle (.5pt);
\fill (0.275, 1.000) circle (.5pt);
\fill (0.150, 0.650) circle (.5pt);
\fill (0.350, 0.850) circle (.5pt);
\fill (0.300, 0.700) circle (.5pt);
\fill (0.200, 0.800) circle (.5pt);
\fill (0.725, 1.000) circle (.5pt);
\fill (0.650, 0.850) circle (.5pt);
\fill (0.800, 0.800) circle (.5pt);
\fill (0.700, 0.700) circle (.5pt);
\fill (0.850, 0.650) circle (.5pt);
\draw[thick] (0.000,0.000) to (0.275, 0.000);
\draw[thick] (0.000,0.000) to (0.000, 0.275);
\draw[thick] (0.000,0.000) to (0.200, 0.200);
\draw[thick] (0.500,0.000) to (0.275, 0.000);
\draw[thick] (0.500,0.000) to (0.725, 0.000);
\draw[thick] (0.500,0.000) to (0.500, 0.225);
\draw[thick] (0.500,0.000) to (0.350, 0.150);
\draw[thick] (0.500,0.000) to (0.650, 0.150);
\draw[thick] (1.000,0.000) to (0.725, 0.000);
\draw[thick] (1.000,0.000) to (1.000, 0.275);
\draw[thick] (1.000,0.000) to (0.800, 0.200);
\draw[thick] (0.000,0.500) to (0.000, 0.275);
\draw[thick] (0.000,0.500) to (0.225, 0.500);
\draw[thick] (0.000,0.500) to (0.150, 0.350);
\draw[thick] (0.000,0.500) to (0.000, 0.725);
\draw[thick] (0.000,0.500) to (0.150, 0.650);
\draw[thick] (0.500,0.500) to (0.500, 0.225);
\draw[thick] (0.500,0.500) to (0.225, 0.500);
\draw[thick] (0.500,0.500) to (0.300, 0.300);
\draw[thick] (0.500,0.500) to (0.775, 0.500);
\draw[thick] (0.500,0.500) to (0.700, 0.300);
\draw[thick] (0.500,0.500) to (0.500, 0.775);
\draw[thick] (0.500,0.500) to (0.300, 0.700);
\draw[thick] (0.500,0.500) to (0.700, 0.700);
\draw[thick] (1.000,0.500) to (1.000, 0.275);
\draw[thick] (1.000,0.500) to (1.000, 0.725);
\draw[thick] (1.000,0.500) to (0.775, 0.500);
\draw[thick] (1.000,0.500) to (0.850, 0.350);
\draw[thick] (1.000,0.500) to (0.850, 0.650);
\draw[thick] (0.000,1.000) to (0.000, 0.725);
\draw[thick] (0.000,1.000) to (0.275, 1.000);
\draw[thick] (0.000,1.000) to (0.200, 0.800);
\draw[thick] (0.500,1.000) to (0.500, 0.775);
\draw[thick] (0.500,1.000) to (0.275, 1.000);
\draw[thick] (0.500,1.000) to (0.350, 0.850);
\draw[thick] (0.500,1.000) to (0.725, 1.000);
\draw[thick] (0.500,1.000) to (0.650, 0.850);
\draw[thick] (1.000,1.000) to (1.000, 0.725);
\draw[thick] (1.000,1.000) to (0.725, 1.000);
\draw[thick] (1.000,1.000) to (0.800, 0.800);
\draw[thick] (0.275,0.000) to (0.350, 0.150);
\draw[thick] (0.275,0.000) to (0.200, 0.200);
\draw[thick] (0.725,0.000) to (0.650, 0.150);
\draw[thick] (0.725,0.000) to (0.800, 0.200);
\draw[thick] (1.000,0.275) to (0.850, 0.350);
\draw[thick] (1.000,0.275) to (0.800, 0.200);
\draw[thick] (1.000,0.725) to (0.800, 0.800);
\draw[thick] (1.000,0.725) to (0.850, 0.650);
\draw[thick] (0.000,0.275) to (0.150, 0.350);
\draw[thick] (0.000,0.275) to (0.200, 0.200);
\draw[thick] (0.500,0.225) to (0.350, 0.150);
\draw[thick] (0.500,0.225) to (0.300, 0.300);
\draw[thick] (0.500,0.225) to (0.650, 0.150);
\draw[thick] (0.500,0.225) to (0.700, 0.300);
\draw[thick] (0.225,0.500) to (0.150, 0.350);
\draw[thick] (0.225,0.500) to (0.300, 0.300);
\draw[thick] (0.225,0.500) to (0.150, 0.650);
\draw[thick] (0.225,0.500) to (0.300, 0.700);
\draw[thick] (0.150,0.350) to (0.300, 0.300);
\draw[thick] (0.150,0.350) to (0.200, 0.200);
\draw[thick] (0.350,0.150) to (0.300, 0.300);
\draw[thick] (0.350,0.150) to (0.200, 0.200);
\draw[thick] (0.300,0.300) to (0.200, 0.200);
\draw[thick] (0.775,0.500) to (0.850, 0.350);
\draw[thick] (0.775,0.500) to (0.700, 0.300);
\draw[thick] (0.775,0.500) to (0.700, 0.700);
\draw[thick] (0.775,0.500) to (0.850, 0.650);
\draw[thick] (0.650,0.150) to (0.800, 0.200);
\draw[thick] (0.650,0.150) to (0.700, 0.300);
\draw[thick] (0.850,0.350) to (0.800, 0.200);
\draw[thick] (0.850,0.350) to (0.700, 0.300);
\draw[thick] (0.800,0.200) to (0.700, 0.300);
\draw[thick] (0.000,0.725) to (0.150, 0.650);
\draw[thick] (0.000,0.725) to (0.200, 0.800);
\draw[thick] (0.500,0.775) to (0.350, 0.850);
\draw[thick] (0.500,0.775) to (0.300, 0.700);
\draw[thick] (0.500,0.775) to (0.650, 0.850);
\draw[thick] (0.500,0.775) to (0.700, 0.700);
\draw[thick] (0.275,1.000) to (0.350, 0.850);
\draw[thick] (0.275,1.000) to (0.200, 0.800);
\draw[thick] (0.150,0.650) to (0.300, 0.700);
\draw[thick] (0.150,0.650) to (0.200, 0.800);
\draw[thick] (0.350,0.850) to (0.300, 0.700);
\draw[thick] (0.350,0.850) to (0.200, 0.800);
\draw[thick] (0.300,0.700) to (0.200, 0.800);
\draw[thick] (0.725,1.000) to (0.650, 0.850);
\draw[thick] (0.725,1.000) to (0.800, 0.800);
\draw[thick] (0.650,0.850) to (0.800, 0.800);
\draw[thick] (0.650,0.850) to (0.700, 0.700);
\draw[thick] (0.800,0.800) to (0.700, 0.700);
\draw[thick] (0.800,0.800) to (0.850, 0.650);
\draw[thick] (0.700,0.700) to (0.850, 0.650);
\fill[red] (0.287, 0.212) circle (.5pt);
\fill[red] (0.637, 0.562) circle (.5pt);
\fill[red] (0.713, 0.787) circle (.5pt);
\fill[red] (0.748, 0.895) circle (.5pt);
\fill[red] (0.038, 0.388) circle (.5pt);
\fill[red] (0.887, 0.462) circle (.5pt);
\fill[red] (0.863, 0.938) circle (.5pt);
\fill[red] (0.752, 0.605) circle (.5pt);
\fill[red] (0.887, 0.538) circle (.5pt);
\fill[red] (0.563, 0.362) circle (.5pt);
\fill[red] (0.138, 0.063) circle (.5pt);
\fill[red] (0.212, 0.287) circle (.5pt);
\fill[red] (0.388, 0.038) circle (.5pt);
\fill[red] (0.752, 0.395) circle (.5pt);
\fill[red] (0.113, 0.463) circle (.5pt);
\fill[red] (0.462, 0.112) circle (.5pt);
\fill[red] (0.395, 0.248) circle (.5pt);
\fill[red] (0.105, 0.252) circle (.5pt);
\fill[red] (0.248, 0.395) circle (.5pt);
\fill[red] (0.063, 0.138) circle (.5pt);
\fill[red] (0.252, 0.105) circle (.5pt);
\fill[red] (0.113, 0.537) circle (.5pt);
\fill[red] (0.038, 0.612) circle (.5pt);
\fill[red] (0.213, 0.713) circle (.5pt);
\fill[red] (0.062, 0.863) circle (.5pt);
\fill[red] (0.538, 0.887) circle (.5pt);
\fill[red] (0.252, 0.895) circle (.5pt);
\fill[red] (0.248, 0.605) circle (.5pt);
\fill[red] (0.138, 0.938) circle (.5pt);
\fill[red] (0.105, 0.748) circle (.5pt);
\fill[red] (0.287, 0.787) circle (.5pt);
\fill[red] (0.363, 0.563) circle (.5pt);
\fill[red] (0.563, 0.637) circle (.5pt);
\fill[red] (0.363, 0.438) circle (.5pt);
\fill[red] (0.395, 0.752) circle (.5pt);
\fill[red] (0.438, 0.638) circle (.5pt);
\fill[red] (0.388, 0.962) circle (.5pt);
\fill[red] (0.462, 0.887) circle (.5pt);
\fill[red] (0.613, 0.962) circle (.5pt);
\fill[red] (0.605, 0.752) circle (.5pt);
\fill[red] (0.938, 0.863) circle (.5pt);
\fill[red] (0.895, 0.748) circle (.5pt);
\fill[red] (0.787, 0.713) circle (.5pt);
\fill[red] (0.962, 0.613) circle (.5pt);
\fill[red] (0.962, 0.388) circle (.5pt);
\fill[red] (0.713, 0.212) circle (.5pt);
\fill[red] (0.438, 0.362) circle (.5pt);
\fill[red] (0.748, 0.105) circle (.5pt);
\fill[red] (0.637, 0.438) circle (.5pt);
\fill[red] (0.537, 0.113) circle (.5pt);
\fill[red] (0.605, 0.248) circle (.5pt);
\fill[red] (0.613, 0.038) circle (.5pt);
\fill[red] (0.862, 0.062) circle (.5pt);
\fill[red] (0.895, 0.252) circle (.5pt);
\fill[red] (0.787, 0.288) circle (.5pt);
\fill[red] (0.938, 0.138) circle (.5pt);
\draw[dashed, red, thick] (0.138, 0.063) to (0.138, 0.000);
\draw[dashed, red, thick] (0.063, 0.138) to (0.000, 0.138);
\draw[dashed, red, thick] (0.138,0.063) to (0.063, 0.138);
\draw[dashed, red, thick] (0.388, 0.038) to (0.388, 0.000);
\draw[dashed, red, thick] (0.613, 0.038) to (0.613, 0.000);
\draw[dashed, red, thick] (0.462,0.112) to (0.537, 0.113);
\draw[dashed, red, thick] (0.388,0.038) to (0.462, 0.112);
\draw[dashed, red, thick] (0.537,0.113) to (0.613, 0.038);
\draw[dashed, red, thick] (0.862, 0.062) to (0.863, 0.000);
\draw[dashed, red, thick] (0.938, 0.138) to (1.000, 0.138);
\draw[dashed, red, thick] (0.862,0.062) to (0.938, 0.138);
\draw[dashed, red, thick] (0.038, 0.388) to (0.000, 0.388);
\draw[dashed, red, thick] (0.113,0.463) to (0.113, 0.537);
\draw[dashed, red, thick] (0.038,0.388) to (0.113, 0.463);
\draw[dashed, red, thick] (0.038, 0.612) to (0.000, 0.613);
\draw[dashed, red, thick] (0.113,0.537) to (0.038, 0.612);
\draw[dashed, red, thick] (0.563,0.362) to (0.438, 0.362);
\draw[dashed, red, thick] (0.363,0.563) to (0.363, 0.438);
\draw[dashed, red, thick] (0.363,0.438) to (0.438, 0.362);
\draw[dashed, red, thick] (0.637,0.562) to (0.637, 0.438);
\draw[dashed, red, thick] (0.563,0.362) to (0.637, 0.438);
\draw[dashed, red, thick] (0.563,0.637) to (0.438, 0.638);
\draw[dashed, red, thick] (0.363,0.563) to (0.438, 0.638);
\draw[dashed, red, thick] (0.637,0.562) to (0.563, 0.637);
\draw[dashed, red, thick] (0.962, 0.388) to (1.000, 0.388);
\draw[dashed, red, thick] (0.962, 0.613) to (1.000, 0.613);
\draw[dashed, red, thick] (0.887,0.462) to (0.887, 0.538);
\draw[dashed, red, thick] (0.887,0.462) to (0.962, 0.388);
\draw[dashed, red, thick] (0.887,0.538) to (0.962, 0.613);
\draw[dashed, red, thick] (0.062, 0.863) to (0.000, 0.863);
\draw[dashed, red, thick] (0.138, 0.938) to (0.138, 1.000);
\draw[dashed, red, thick] (0.062,0.863) to (0.138, 0.938);
\draw[dashed, red, thick] (0.538,0.887) to (0.462, 0.887);
\draw[dashed, red, thick] (0.388, 0.962) to (0.388, 1.000);
\draw[dashed, red, thick] (0.388,0.962) to (0.462, 0.887);
\draw[dashed, red, thick] (0.613, 0.962) to (0.613, 1.000);
\draw[dashed, red, thick] (0.538,0.887) to (0.613, 0.962);
\draw[dashed, red, thick] (0.938, 0.863) to (1.000, 0.863);
\draw[dashed, red, thick] (0.863, 0.938) to (0.863, 1.000);
\draw[dashed, red, thick] (0.863,0.938) to (0.938, 0.863);
\draw[dashed, red, thick] (0.252,0.105) to (0.388, 0.038);
\draw[dashed, red, thick] (0.252,0.105) to (0.138, 0.063);
\draw[dashed, red, thick] (0.748,0.105) to (0.613, 0.038);
\draw[dashed, red, thick] (0.748,0.105) to (0.862, 0.062);
\draw[dashed, red, thick] (0.895,0.252) to (0.962, 0.388);
\draw[dashed, red, thick] (0.895,0.252) to (0.938, 0.138);
\draw[dashed, red, thick] (0.895,0.748) to (0.938, 0.863);
\draw[dashed, red, thick] (0.895,0.748) to (0.962, 0.613);
\draw[dashed, red, thick] (0.105,0.252) to (0.038, 0.388);
\draw[dashed, red, thick] (0.105,0.252) to (0.063, 0.138);
\draw[dashed, red, thick] (0.395,0.248) to (0.462, 0.112);
\draw[dashed, red, thick] (0.395,0.248) to (0.438, 0.362);
\draw[dashed, red, thick] (0.605,0.248) to (0.537, 0.113);
\draw[dashed, red, thick] (0.605,0.248) to (0.563, 0.362);
\draw[dashed, red, thick] (0.248,0.395) to (0.113, 0.463);
\draw[dashed, red, thick] (0.248,0.395) to (0.363, 0.438);
\draw[dashed, red, thick] (0.248,0.605) to (0.113, 0.537);
\draw[dashed, red, thick] (0.248,0.605) to (0.363, 0.563);
\draw[dashed, red, thick] (0.212,0.287) to (0.248, 0.395);
\draw[dashed, red, thick] (0.212,0.287) to (0.105, 0.252);
\draw[dashed, red, thick] (0.287,0.212) to (0.395, 0.248);
\draw[dashed, red, thick] (0.287,0.212) to (0.252, 0.105);
\draw[dashed, red, thick] (0.287,0.212) to (0.212, 0.287);
\draw[dashed, red, thick] (0.752,0.395) to (0.887, 0.462);
\draw[dashed, red, thick] (0.752,0.395) to (0.637, 0.438);
\draw[dashed, red, thick] (0.752,0.605) to (0.637, 0.562);
\draw[dashed, red, thick] (0.752,0.605) to (0.887, 0.538);
\draw[dashed, red, thick] (0.713,0.212) to (0.748, 0.105);
\draw[dashed, red, thick] (0.713,0.212) to (0.605, 0.248);
\draw[dashed, red, thick] (0.787,0.288) to (0.895, 0.252);
\draw[dashed, red, thick] (0.787,0.288) to (0.752, 0.395);
\draw[dashed, red, thick] (0.713,0.212) to (0.787, 0.288);
\draw[dashed, red, thick] (0.105,0.748) to (0.038, 0.612);
\draw[dashed, red, thick] (0.105,0.748) to (0.062, 0.863);
\draw[dashed, red, thick] (0.395,0.752) to (0.462, 0.887);
\draw[dashed, red, thick] (0.395,0.752) to (0.438, 0.638);
\draw[dashed, red, thick] (0.605,0.752) to (0.538, 0.887);
\draw[dashed, red, thick] (0.605,0.752) to (0.563, 0.637);
\draw[dashed, red, thick] (0.252,0.895) to (0.388, 0.962);
\draw[dashed, red, thick] (0.252,0.895) to (0.138, 0.938);
\draw[dashed, red, thick] (0.213,0.713) to (0.248, 0.605);
\draw[dashed, red, thick] (0.213,0.713) to (0.105, 0.748);
\draw[dashed, red, thick] (0.287,0.787) to (0.395, 0.752);
\draw[dashed, red, thick] (0.287,0.787) to (0.252, 0.895);
\draw[dashed, red, thick] (0.213,0.713) to (0.287, 0.787);
\draw[dashed, red, thick] (0.748,0.895) to (0.613, 0.962);
\draw[dashed, red, thick] (0.748,0.895) to (0.863, 0.938);
\draw[dashed, red, thick] (0.713,0.787) to (0.748, 0.895);
\draw[dashed, red, thick] (0.713,0.787) to (0.605, 0.752);
\draw[dashed, red, thick] (0.787,0.713) to (0.713, 0.787);
\draw[dashed, red, thick] (0.787,0.713) to (0.895, 0.748);
\draw[dashed, red, thick] (0.752,0.605) to (0.787, 0.713);
\end{tikzpicture}
  \caption{$h=1/2$}
  \label{fig:square_mesh2}
\end{subfigure}
\caption{Primal and dual mesh on the unit square, $[0,1] \times [0,1]$. The primal mesh is drawn with solid black lines and the dual mesh is shown in red dashed lines. A base, coarse mesh is shown in (a), and one refinement has been done in (b). }  \label{fig:square_mesh}
\end{figure}

To verify the theoretical results presented in this paper, the $H(\grad)$ and $\bm H(\curl)$ schemes are implemented in two dimensions in MATLAB.  
%Results for $\bm H(\divg)$ have been excluded since, in two dimensions, the divergence is a rotation of the curl. \xh{I am not sure about the last sentence...}
Two types of meshes are designed to test the MFD method. The first uses a hexagonal domain (see Figure \ref{fig:hexagon_mesh}), which has a Delaunay triangulation for the primal mesh and a non-degenerate dual Voronoi mesh, meaning the Voronoi points lie within their corresponding Delaunay element. More resolved meshes are obtained through uniform refinement and the error estimates are verified on these grids.  %where Table \ref{table:mesh} lists the geometric information.  
As many classical convection-diffusion problems with jumps and boundary layers are simulated on the unit square, we design a primal Delaunay triangulation and dual Voronoi mesh on $[0,1] \times [0,1]$. The refinement method, which preserves non-degeneracy of the meshes, is demonstrated in Figure \ref{fig:square_mesh}.

%However, the Delaunay triangulation is comprised entirely of right triangles which gives a degenerate Voronoi diagram. This is because each Voronoi point lies on the hypotenuse of the corresponding Delaunay triangle, which collapses two Voronoi points into one. Thus, the dual mesh we use is not a Voronoi diagram, and is designed to have one-to-one relationships between nodes on one mesh and elements on the other, as well as edges on one mesh to edges on the other. It also maintains the orthogonality of primal and dual edges, and is non-degenerate. 

%\begin{table}[H]
%\centering
%\begin{tabular}{|c||cccc|}
%\hline 
% & $h$ & Vertices & Edges & Elements \\ 
%\hline 
%\hline 
%Mesh 1 & $1/4$ & 61 & 156 & 96 \\ 
%\hline 
%Mesh 2 & $1/8$ & 217 & 600 & 384\\ 
%\hline 
%Mesh 3 & $1/16$ & 817 & 2,352 & 1,536 \\ 
%\hline 
%Mesh 4 & $1/32$ & 3,169  & 9,312 & 6,144  \\ 
%\hline 
%Mesh 5 & $1/64$ & 12,481  & 37,056 & 24,576 \\ 
%\hline 
%Mesh 6 & $1/128$ & 49,537  & 147,840 & 98,304\\ 
%\hline 
%Mesh 7 & $1/256$ & 197,377  & 590,592  & 393,216 \\ 
%\hline 
%%Mesh 8 & $1/512$ & 787,969  & 2,360,832  & 1,572,864   \\ 
%%\hline 
%\end{tabular} 
%\caption{Geometric information for the hexgonal Delaunay meshes (see Fig. \ref{fig:hexagon_mesh}).} \label{table:mesh}
%\end{table}

\subsection{Convergence results}\label{subsec:convergence}
First, the error estimates given by \cref{coro:MFD-error} are verified. Here, the domain is the hexagonal domain as in Figure \ref{fig:hexagon_mesh}. 
%where finer meshes are generated through uniform refinement, with mesh statistics listed in Table \ref{table:mesh}. 
For the $H(\grad)$ case, \eqref{eq:convdiff_grad}, the convergence is verified with the exact solution,
\begin{equation}
u(\bm x) = \sin(\pi x_1) \sin( \pi x_2), \label{eq:test4}
\end{equation}
where the right-hand side can be analytically computed. Similarly, the $\bm{H}(\curl)$ convergence results are verified using exact solution,
\begin{equation}
\bm u(\bm x) = \begin{bmatrix}
\sin(\pi x_1) \sin( \pi x_2) \\ \cos(\pi x_1) \cos( \pi x_2) \label{eq:test5}
\end{bmatrix},
\end{equation}
where the right-hand side can again be computed analytically via (\ref{eq:convdiff_curl}). In both cases, the convection coefficient is given by
\begin{equation*}
\bm \beta(\bm x) = \begin{bmatrix} \cos(x_1) + 4 \\ -\sin(x_2)+4 \end{bmatrix},
\end{equation*}
and we use $\gamma = 0$ and $\gamma = 1$ for the $H(\grad)$ and $\bm{H}(\curl)$ cases, respectively. To show the robustness of our schemes, particularly in the convection-dominated regime, we vary the diffusion coefficient $\alpha$ in our numerical experiments. 

\begin{figure}[h!]
	\begin{center}
		\includegraphics[scale=0.43]{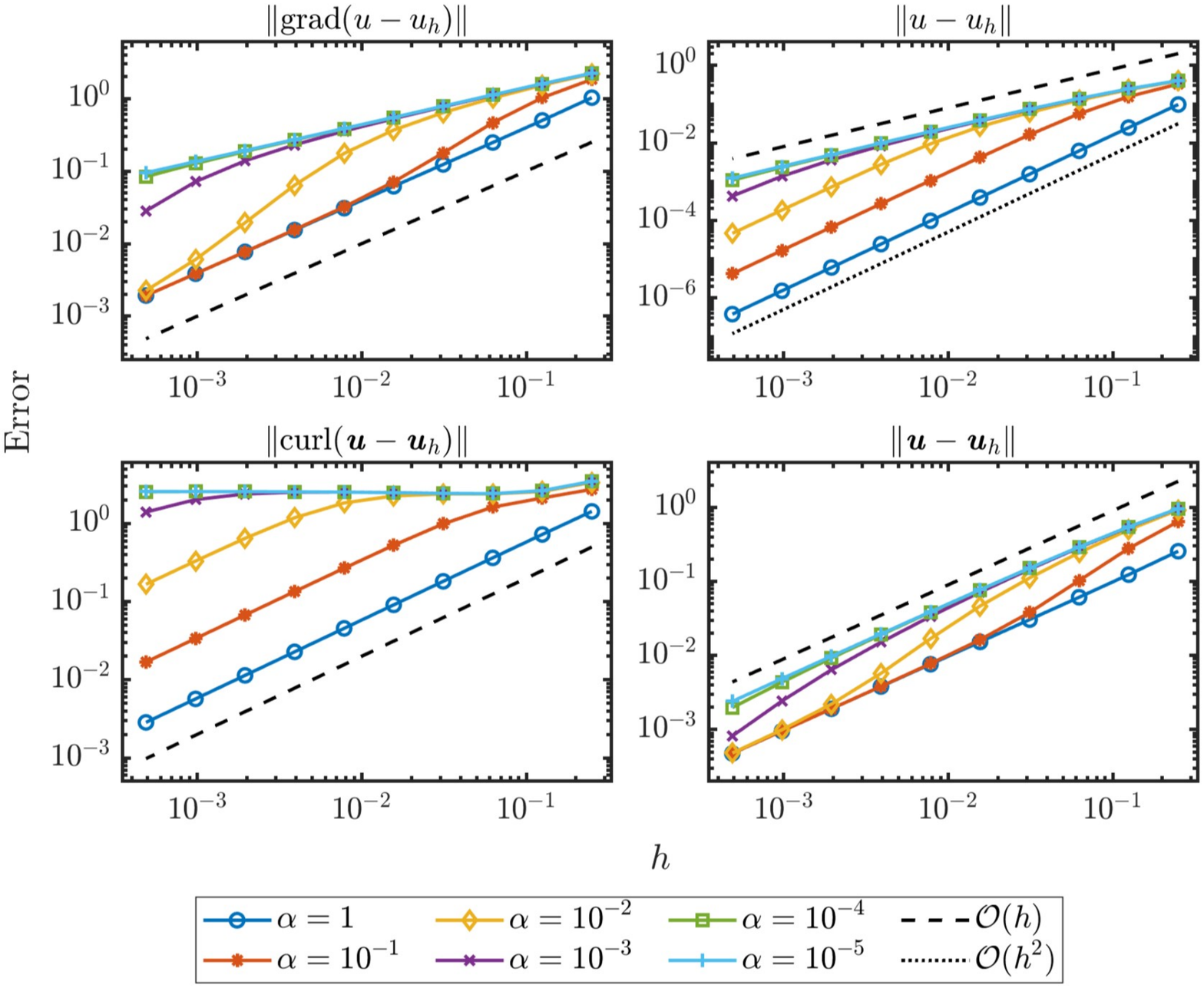}
		\caption{MFD results for solving the convection-diffusion equations with various values of $\alpha$ on a hexagonal computational domain. (Top) Error versus $h$ for scalar equation \eqref{eq:convdiff_grad_MFD} with exact solution \eqref{eq:test4} with error measured in the (left) $H(\grad)$ semi-norm and (right) $L^2$ norm. (Bottom) Error versus $h$ for vector equation \eqref{eq:convdiff_curl_MFD} with exact solution \eqref{eq:test5} with error measured in (left) $\bm H(\curl)$ semi-norm and (right) $L^2$ norm.}
		\label{fig:test_conv}
	\end{center}
\end{figure}

\cref{fig:test_conv} verifies the error estimates given by \cref{thm:MFD-error} and \cref{coro:MFD-error} for the $H(\grad)$ convection-diffusion problem.  The top-left graph of \cref{fig:test_conv} reports the error in the $H(\grad)$ semi-norm, i.e., $|u-u_h|_1 := \| \grad (u - u_h) \|$. We see that for $\alpha = 1$, $10^{-1}$, and $10^{-2}$ the $\mathcal{O}(h)$ rate is eventually recovered as the mesh is refined. In general, we expect to see the trend that the rate approaches $\mathcal{O}(h)$ once $h < | \alpha |/ |\bm{\beta} |$.  For example, for $\alpha = 10^{-3}$ the convergence rate is heading toward $\mathcal{O}(h)$ as $h$ gets closer to $10^{-3}$ (since $\bm{\beta} \sim \mathcal{O}(1)$). More refinement is needed to see a similar result for $\alpha = 10^{-4}$ and $10^{-5}$.  This is consistent with our theoretical results, which require $h$ to be sufficiently small. 

Since the MFD scheme can be considered as a mass-lumped SAFE scheme, in the $L^2$ norm, we would expect the error to converge as $\mathcal{O}(h^2)$ for sufficiently small $h$.  
This is numerically verified in the top-right graph of \cref{fig:test_conv}. For small $\alpha$, we see closer to $\mathcal{O}(h)$ convergence, which asymptotically approaches $\mathcal{O}(h^2)$ with finer meshes, whereas the large values of $\alpha$ clearly demonstrate the $\mathcal{O}(h^2)$ convergence. Once again, the turn from $\mathcal{O}(h)$ to $\mathcal{O}(h^2)$ happens roughly when $h$ becomes smaller than $| \alpha | / |\bm \beta |$. 
%This is clearly seen for $\alpha = 10^{-1}, 10^{-2}$, and more refinement is needed to see the trend for the smaller values of $\alpha$.

Next, we examine the convergence of the $\bm H(\curl)$ convection-diffusion problem. The $\bm H(\curl)$ semi-norm error, i.e., $|\bm{u} - \bm{u}_h |_{\bm H(\curl)} := \|\curl(\bm{u} - \bm{u}_h) \|  $, is plotted in the bottom-left graph of \cref{fig:test_conv}. When the mesh spacing, $h$, is larger than $| \alpha | / |\bm{\beta}|$, poor convergence is observed. However, the expected $\mathcal{O}(h)$ convergence is observed when $h \leq |\alpha|/|\bm{\beta}|$. In particular, we see the $\alpha = 10^{-3}$ case begin to trend toward $\mathcal{O}(h)$ around $h = 10^{-3}$. As in the $H(\grad)$ case, more refinement is needed to verify the convergence for smaller values of $\alpha$.  The bottom-right graph of \cref{fig:test_conv} also demonstrates $\mathcal{O}(h)$ convergence in the $L^2$ norm, which is expected due to the connection between our MFD scheme and the SAFE method (for the $\bm{H}(\curl)$ case, the SAFE discretization uses the lowest-order N\'ed\'elec elements).

\subsection{Boundary and internal layer in \boldmath{$H(\grad)$}} \label{subsec:grad_layer}
Next, we examine the stability of the $H(\grad)$ scheme on the domain $\Omega = [0,1] \times [0,1]$ using the meshes given by Figure \ref{fig:square_mesh}, with $h=1/64$. First, we consider a simple example which admits a jump at the boundary, using parameters,
\begin{equation}\label{eq:test1}
\alpha = 10^{-6}, \quad \bm \beta = \begin{bmatrix}
2+x_1 \\1+x_2 
\end{bmatrix}, \quad \gamma = 0, \quad  \quad f= 1, 
\end{equation}

\begin{figure}[h!]
	\centering
	\begin{subfigure}{.48\textwidth}
		\centering
		\includegraphics[scale=.4]{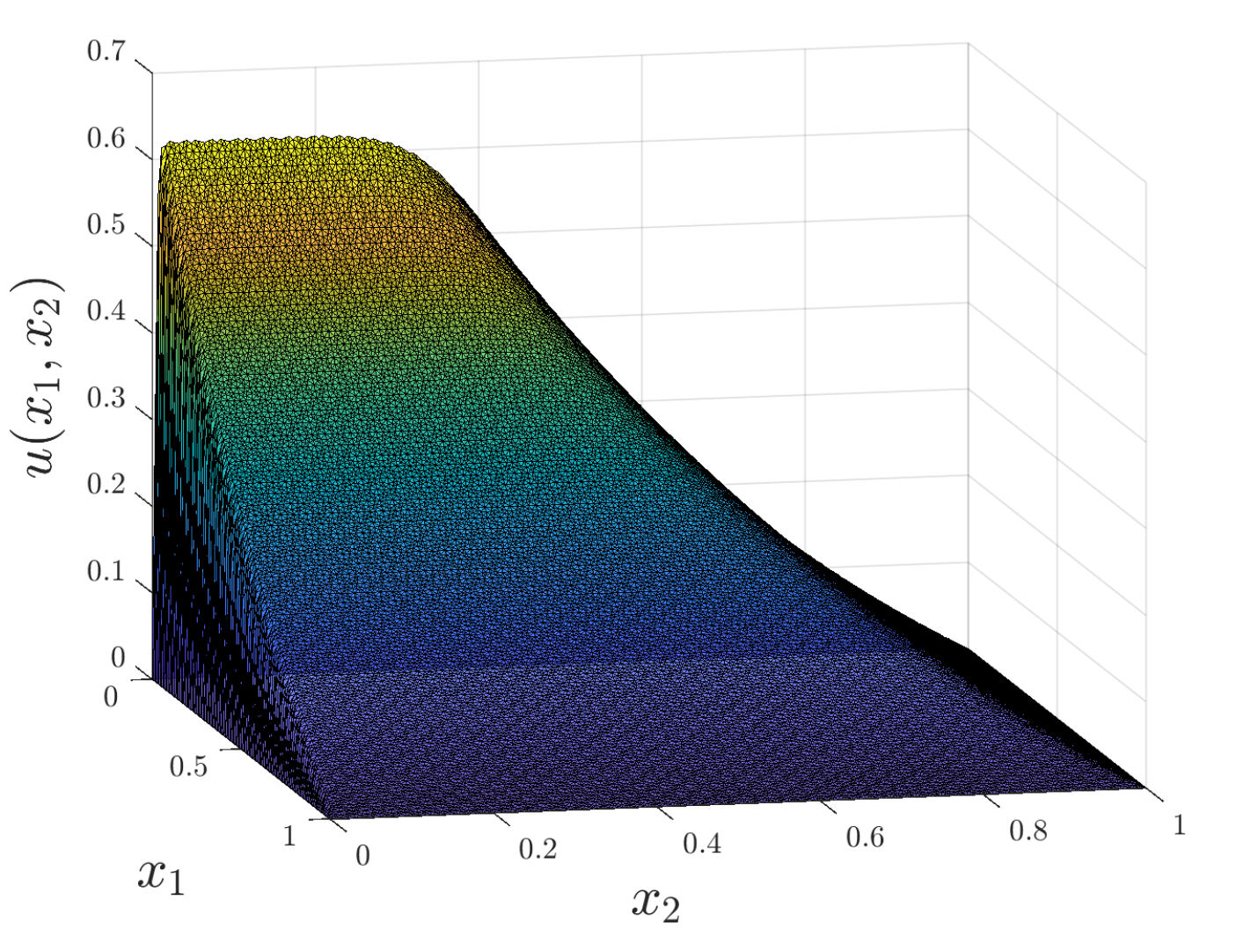}
		\caption{MFD solution of \eqref{eq:convdiff_grad_MFD}.}
		\label{fig:grad_MFD_test}
	\end{subfigure}%
	\begin{subfigure}{.49\textwidth}
		\centering
		\includegraphics[scale=.4]{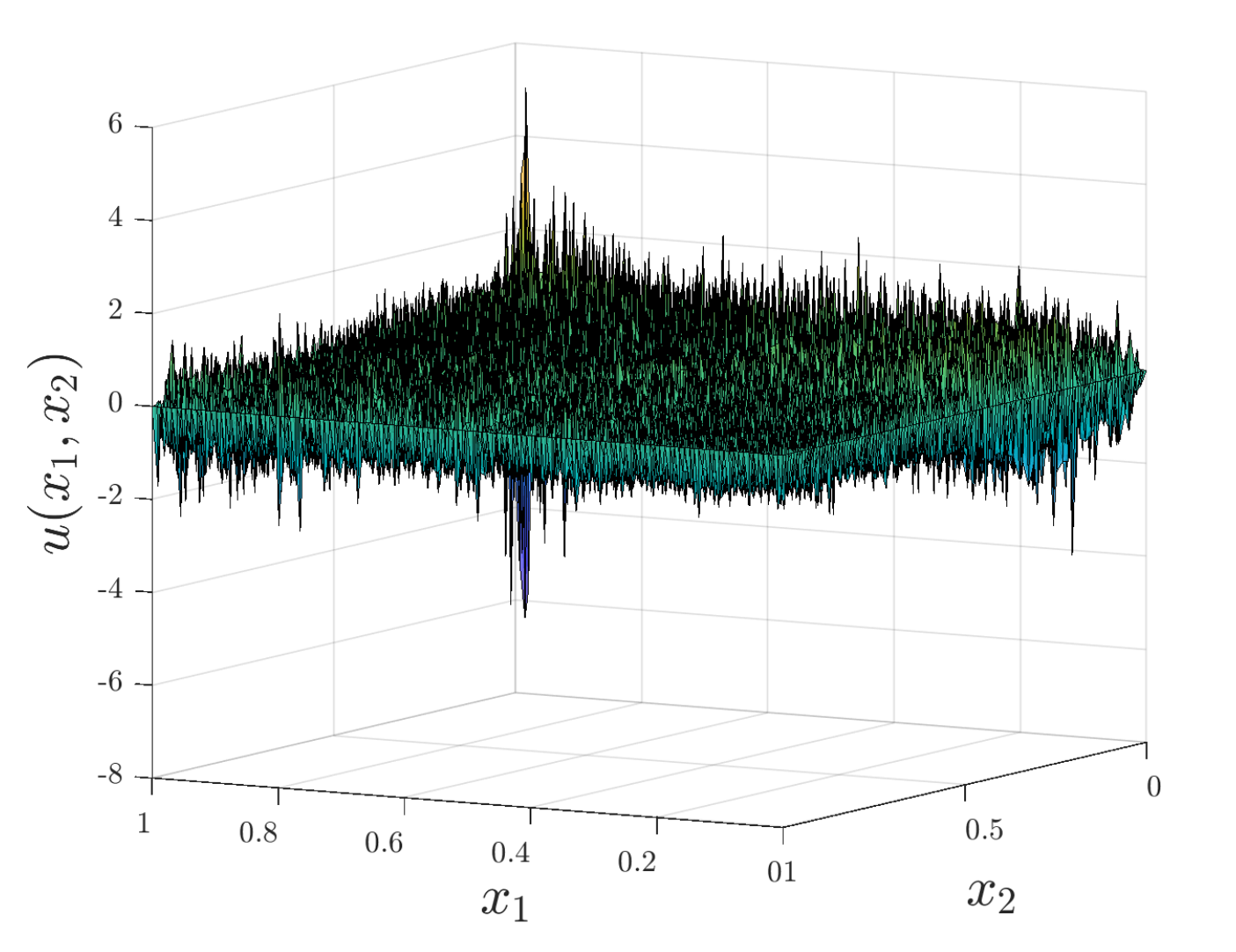}
		\caption{Linear Lagrange FEM solution of \eqref{eq:convdiff_grad}.}
		\label{fig:grad_FEM1}
	\end{subfigure}
		%\caption{Comparison of MFD and FEM solutions for $H$(grad) problem with parameters in \eqref{eq:test1} yielding a boundary layer.} --caption below doesn't compile on my computer -casey
	\caption{Comparison of MFD and FEM solutions for $H$(grad) problem with parameters in \eqref{eq:test1} (e.g. $\alpha = 10^{-6}$) yielding a boundary layer.}
	%: $\alpha = 10^{-6}, \quad \bm \beta = \begin{bmatrix}
	%		1 \\2 
	%	\end{bmatrix}, \quad \gamma = 0, \quad  \quad f= 1$.}
	\label{fig:grad_bound_test}
\end{figure}

Next, we consider a similar problem on the same domain, with parameters,
\begin{equation}\label{eq:test2}
\alpha(\bm{x}) = \begin{cases} 1, & x_1 < 0.5, \\
						10^{-6}, & x_1 \geq 0.5,						
						\end{cases}, \quad \bm \beta = \begin{bmatrix}
2+x_1 \\1+x_2
\end{bmatrix}, \quad \gamma =0,  \quad f= 1,
\end{equation}
where the jump in the diffusion coefficient causes an internal layer.

\begin{figure}[h!]
\centering
\begin{subfigure}{.48\textwidth}
 \centering
 \includegraphics[scale=.4]{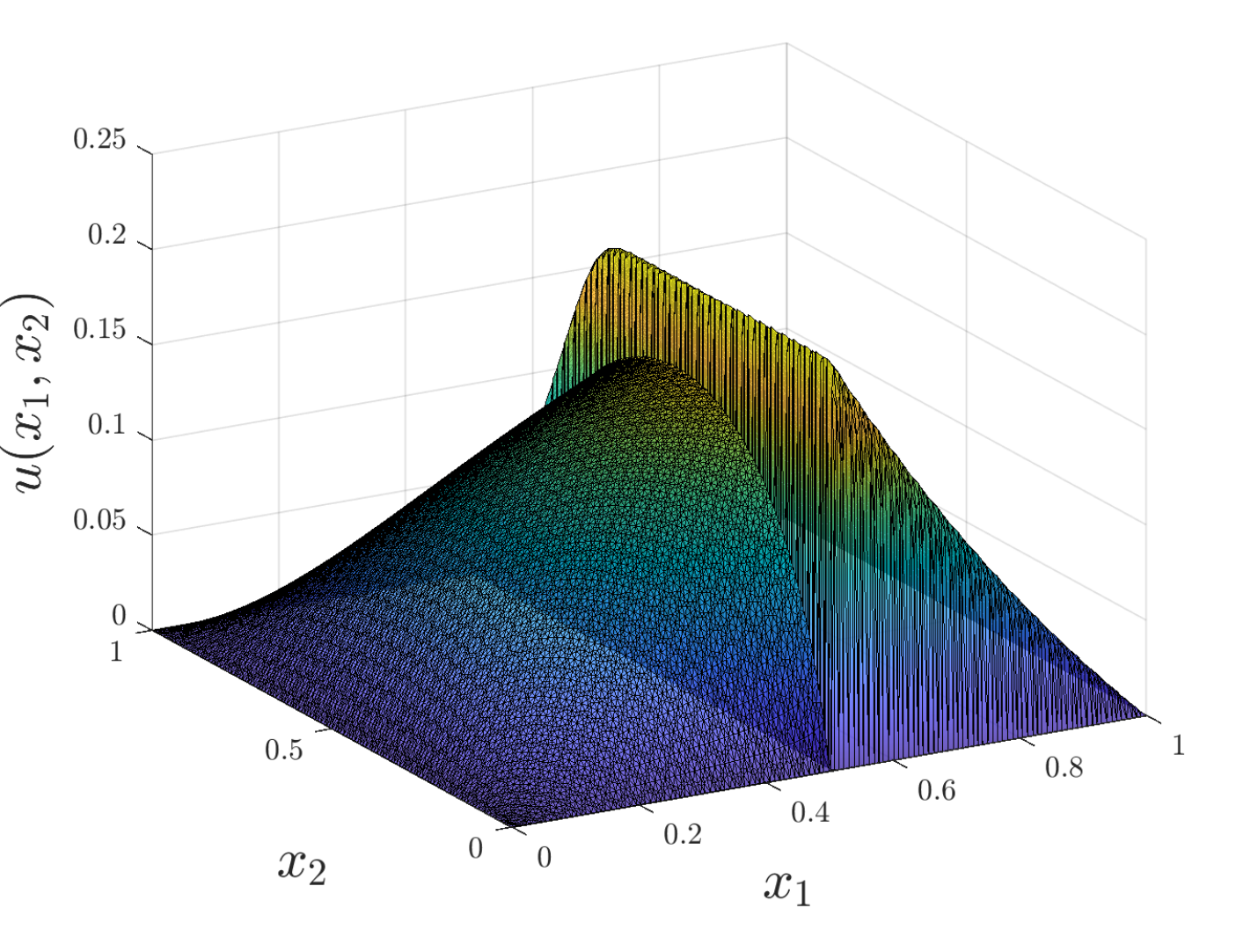}
 \caption{MFD solution of \eqref{eq:convdiff_grad_MFD}.}
 \label{fig:grad_MFD_testint}
\end{subfigure}%
\begin{subfigure}{.49\textwidth}
  \centering
\includegraphics[scale=.4]{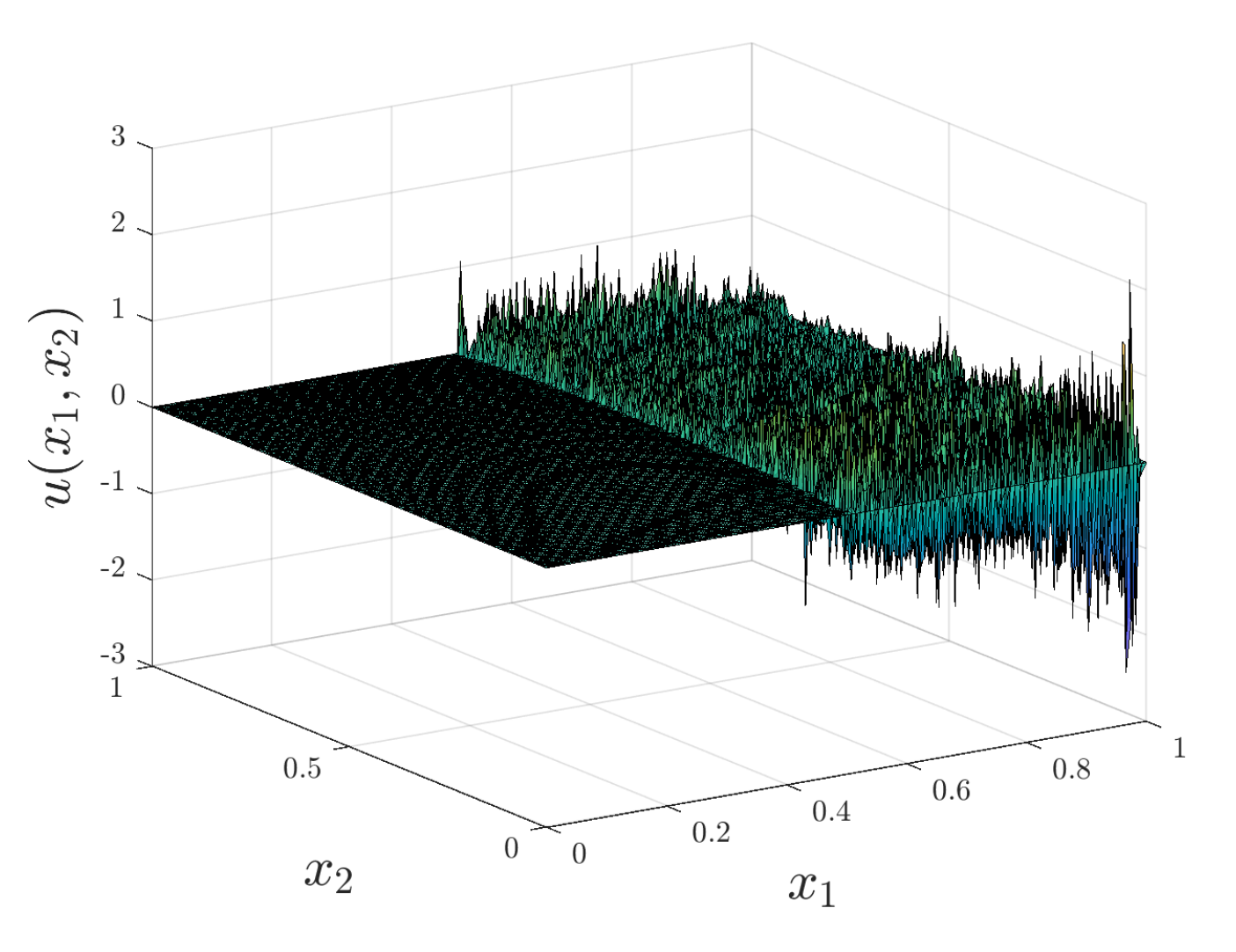}
  \caption{Linear Lagrange FEM solution of \eqref{eq:convdiff_grad}.}
  \label{fig:grad_FEM2}
\end{subfigure}
\caption{Comparison of MFD and FEM solutions for $H(\grad)$ problem with parameters in \eqref{eq:test2} (e.g. $\alpha = 1$ when $x_1 < 0.5$ and $\alpha = 10^{-6}$ when $x_1 \geq 0.5$) yielding an internal layer.}
\label{fig:grad_int_test}
\end{figure}

The solutions given by the MFD methods (solving \eqref{eq:convdiff_grad_MFD}) with parameters given by \eqref{eq:test1} and \eqref{eq:test2} are plotted in Figure \ref{fig:grad_MFD_test} and Figure \ref{fig:grad_MFD_testint}, respectively. It is clearly seen that there are no oscillations in the plots, and the jumps and layers are accurately captured. This demonstrates the monotonicity of the scheme. In contrast, both Figures \ref{fig:grad_FEM1} and \ref{fig:grad_FEM2} plot the solutions to the same boundary and internal layer problems but using a linear Lagrange FEM, which has large oscillations in the solution. 
\remark{Though omitted for brevity, when \eqref{eq:test2} is modified to let \[\alpha(\bm{x}) = \begin{cases} 10^6, & x_1 < 0.5, \\ 10^{-6}, & x_1 \geq 0.5,\end{cases},\] the stability result still holds and yields a plot similar to Figure \eqref{fig:grad_MFD_testint}. This is expected as the analysis is valid for the case where part of the domain is diffusion-dominated.}

%\xh{I think the captions of the plots are not clear, we should at least list $\alpha$ there. Moreover, I would put the MFD and FEM solutions next to each other to make the comparison clear.} \ja{OK, I fixed this up, but not sure having the parameters re-listed is overkill or not.}
%\xh{I change my mind, having the parameters is a overkill... sorry}

\subsection{Boundary layer in \boldmath{$H(\curl)$}}
\begin{figure}[h!]
\centering
\begin{subfigure}{.48\textwidth}
  \centering
\includegraphics[scale=.4]{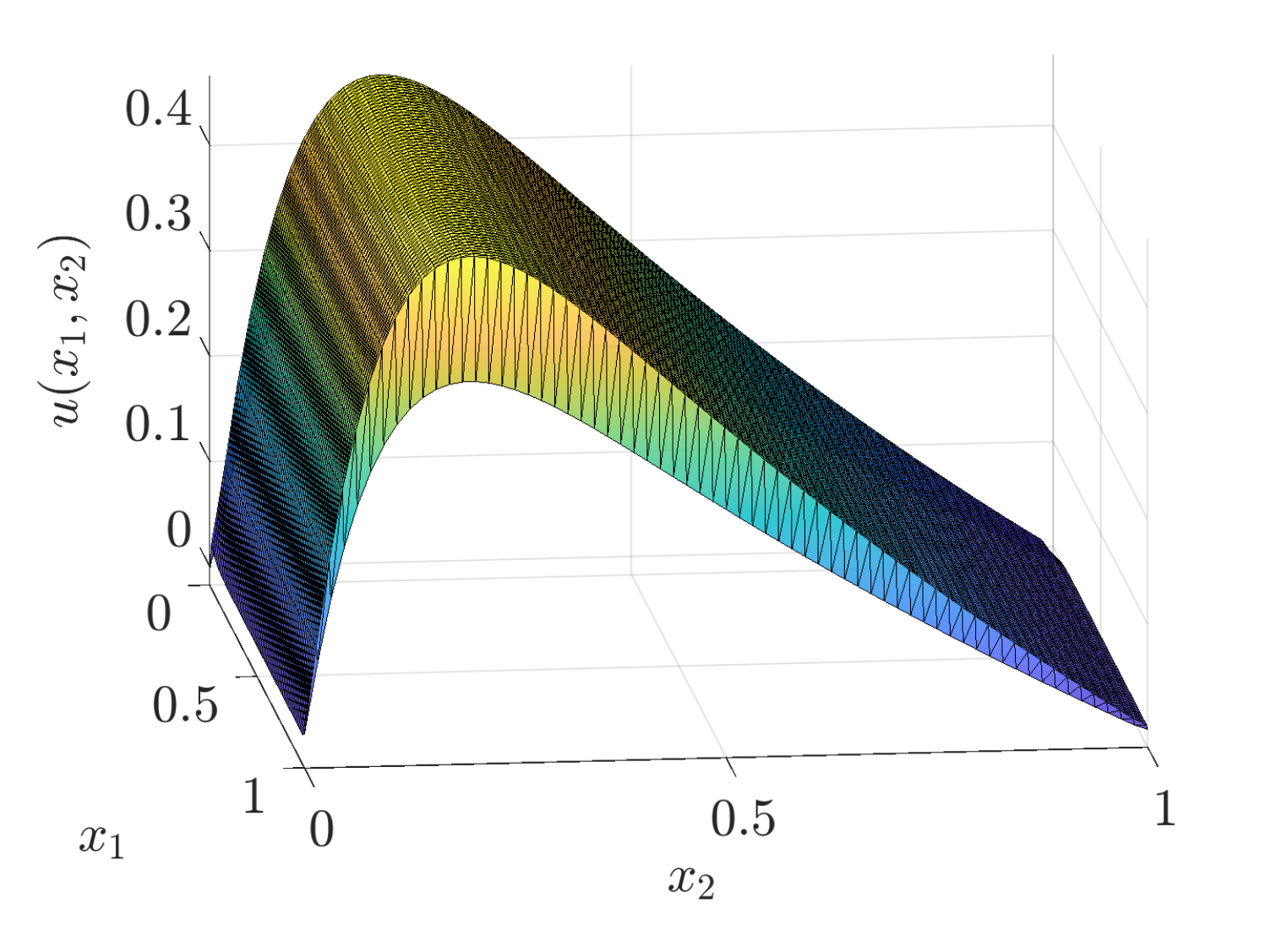}
 \caption{MFD solution of \eqref{eq:convdiff_curl_MFD}.}
  \label{fig:curl_MFD_test1}
\end{subfigure}%
\begin{subfigure}{.49\textwidth}
	\centering
	\includegraphics[scale=.4]{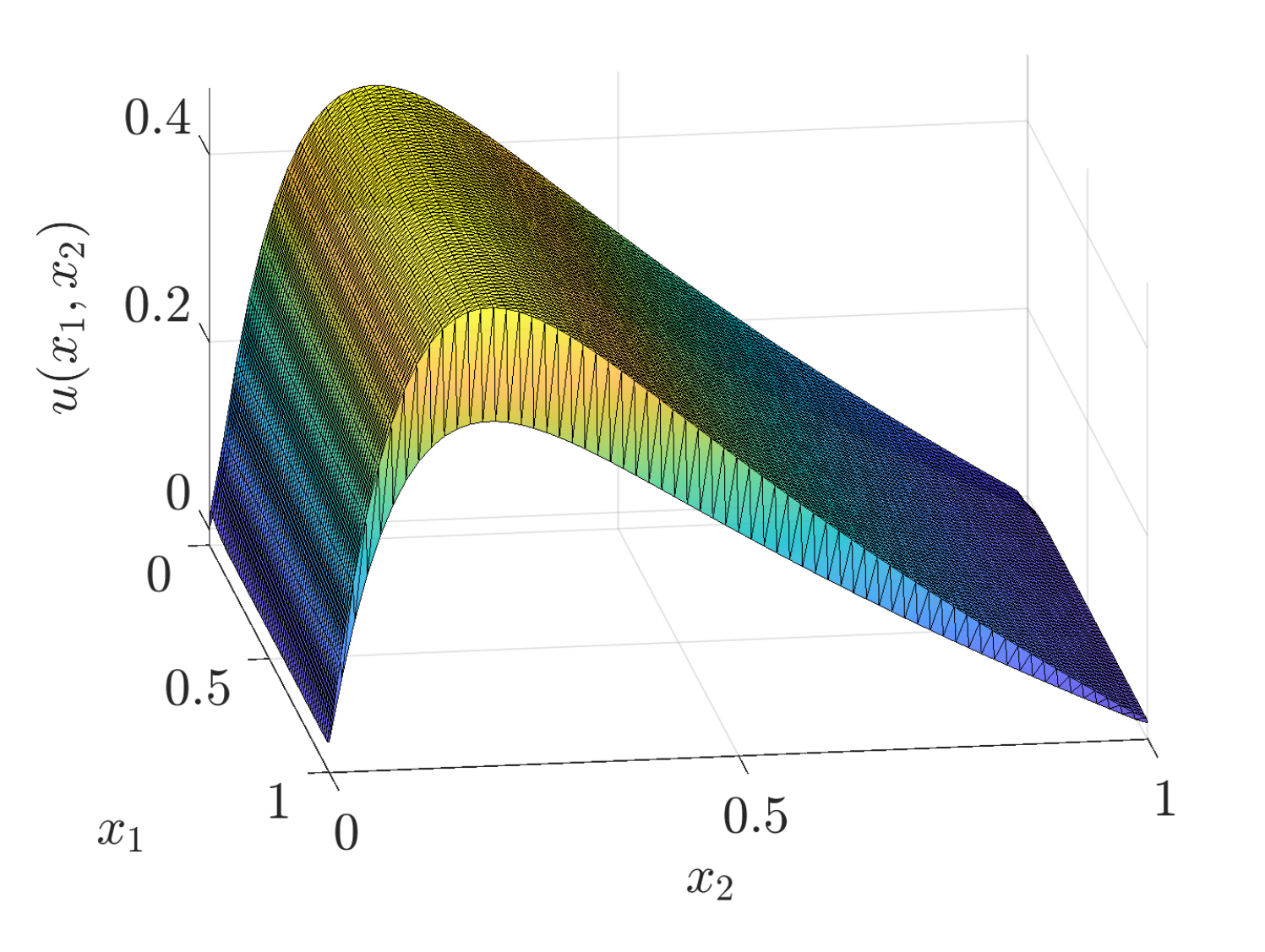}
	\caption{N\'ed\'elec FEM solution of \eqref{eq:convdiff_curl}.}
	\label{fig:curl_FEM_test1}
\end{subfigure}%
\caption{Comparison of MFD and FEM solutions for $\bm H$(curl) problem with parameters in \eqref{eq:test3} and $\alpha = 10^{-1}$}.
\label{fig:curl_test_alphapm2}
\end{figure}

\begin{figure}
\centering
\begin{subfigure}{.48\textwidth}
	\centering
	\includegraphics[scale=.4]{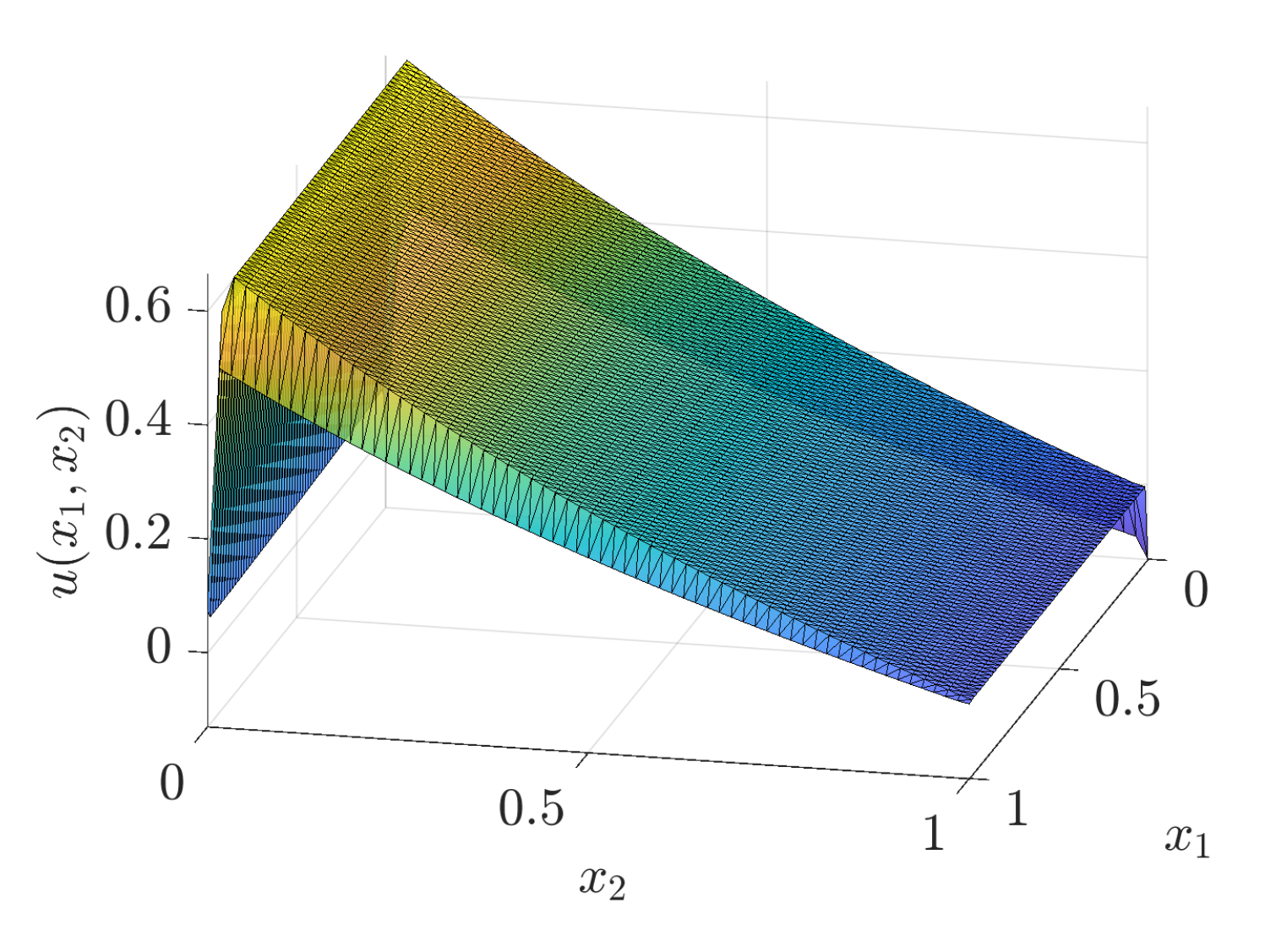}
 \caption{MFD solution of \eqref{eq:convdiff_curl_MFD}.}
	\label{fig:curl_MFD_test2}
\end{subfigure}
\begin{subfigure}{.49\textwidth}
  \centering
\includegraphics[scale=.4]{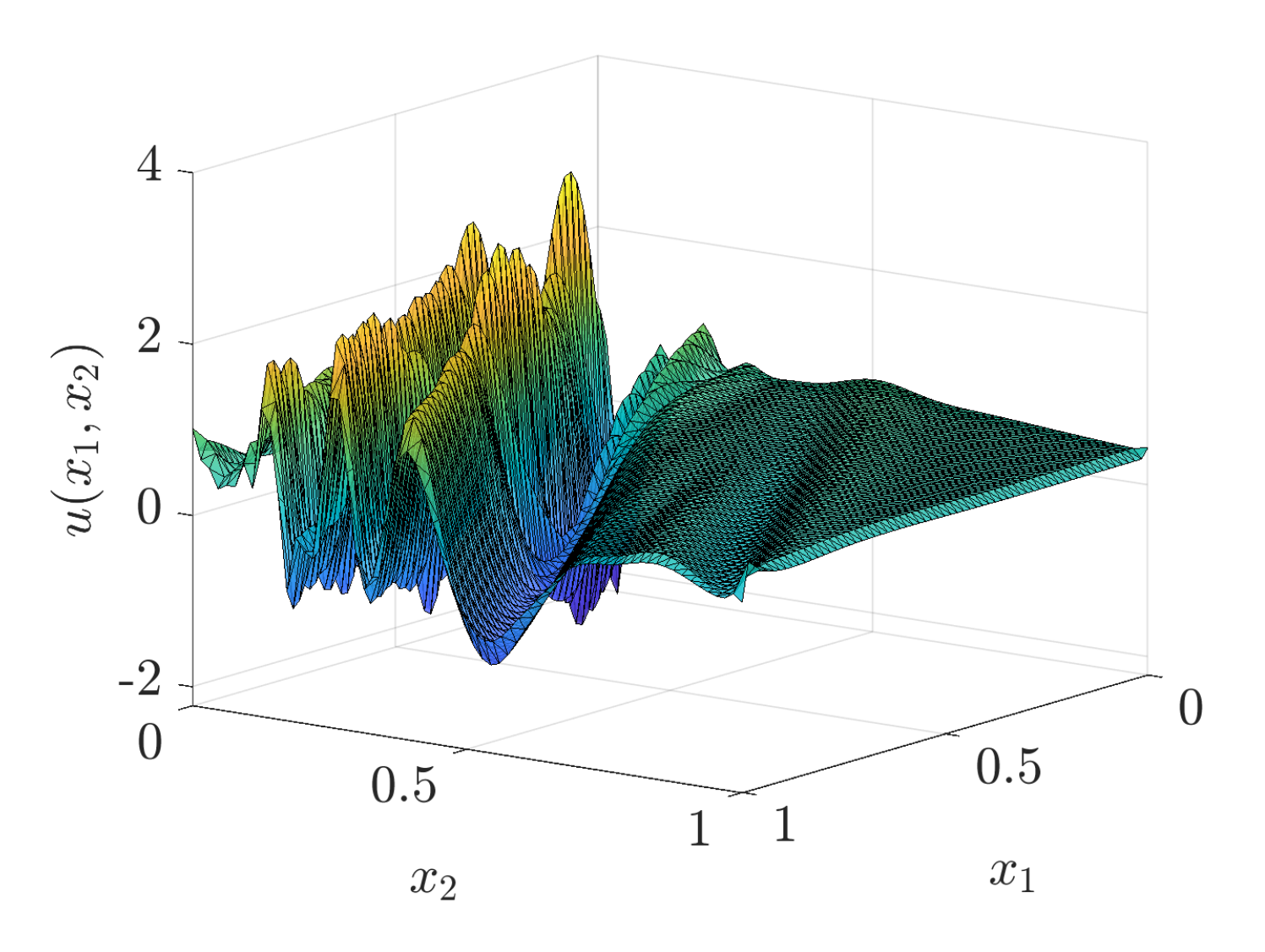}
	\caption{N\'ed\'elec FEM solution of \eqref{eq:convdiff_curl}.}
  \label{fig:curl_FEM_test2}
\end{subfigure}
\caption{Comparison of MFD and FEM solutions for $\bm H$(curl) problem with parameters in \eqref{eq:test3} and $\alpha = 10^{-6}$}.
\label{fig:curl_test_alphapm8}
\end{figure}

For the $\bm{H}(\curl)$ problem, \eqref{eq:convdiff_curl}, we construct an example where a boundary layer is expected with the parameters,
\begin{equation}\label{eq:test3}
\bm \beta = \begin{bmatrix}
2+x_1 \\1+x_2
\end{bmatrix}, \qquad \gamma = 1, \qquad \bm f=  \begin{bmatrix}
1 \\1
\end{bmatrix},
\end{equation}
and varying $\alpha$. The domain is again $\Omega = [0,1] \times [0,1]$ using a uniform mesh with $h=1/64$. The values $\alpha = 10^{-1}$ and $\alpha = 10^{-6}$ are tested, and results are plotted in \cref{fig:curl_test_alphapm2,fig:curl_test_alphapm8}. Again, the figure demonstrates that there are no spurious oscillations when using the MFD scheme. The boundary layer is accurately captured and the solution converges even as $\alpha$ varies by many orders of magnitude. In contrast, for $\alpha = 10^{-1}$ we see almost exactly the same solution as MFD when using the lowest-order N\'ed\'elec finite elements. When $\alpha$ is reduced to $10^{-6}$, the FEM solution is polluted with spurious oscillations.

\subsection{Helmholtz Decomposition in \boldmath{$H(\grad)$}} \label{subsec:helm_num}
\begin{figure}
	\begin{center}
		\includegraphics[scale = .55]{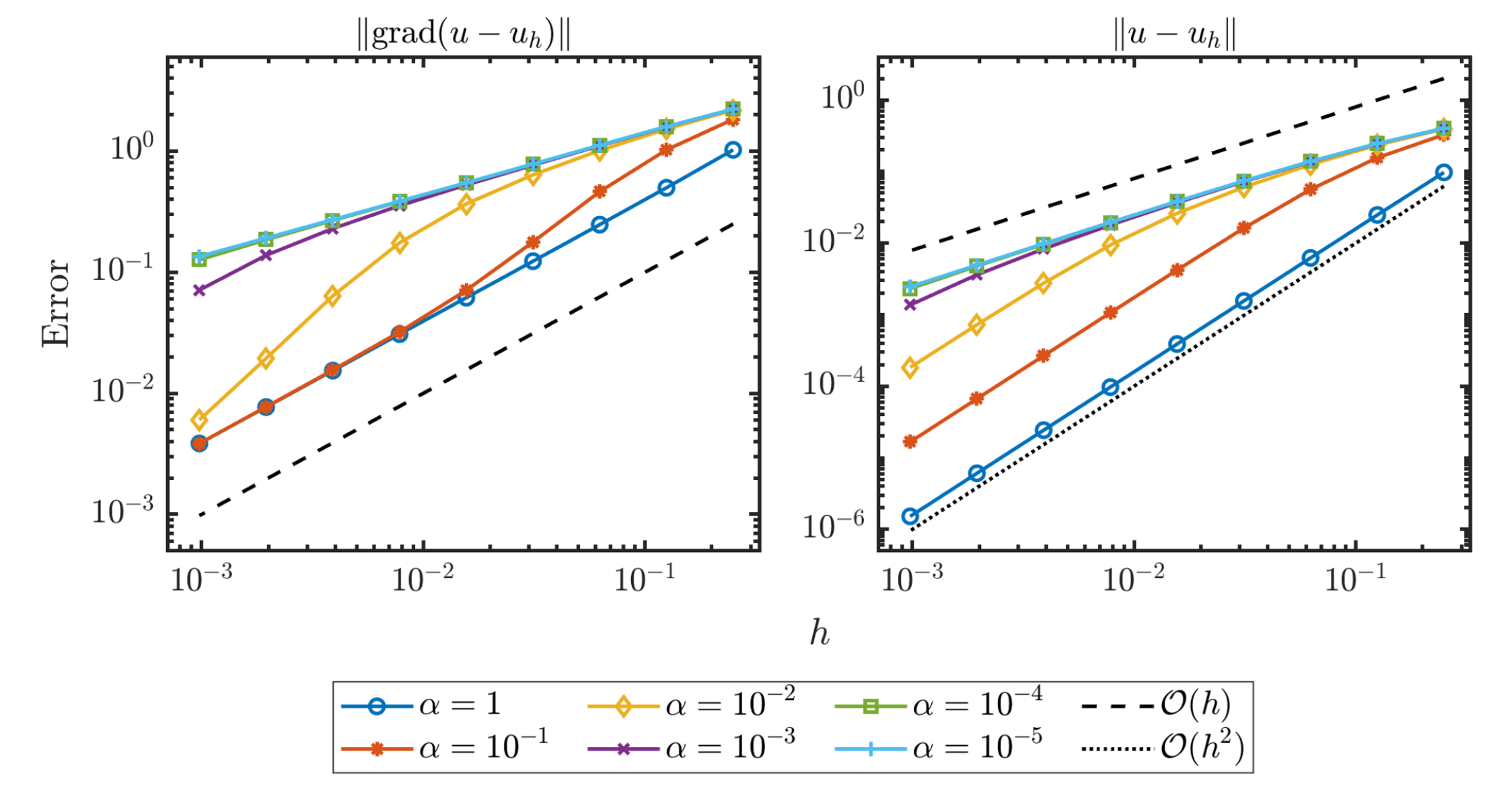}
		\caption{MFD convergence results for modified \eqref{eq:test4}, with error measured in $H(\grad)$ seminorm (left) and $L^2$ norm (right).}
		\label{fig:test_conv_HH}
	\end{center}
\end{figure}

\begin{figure}
\centering
\begin{subfigure}{.48\textwidth}
	\centering
	\includegraphics[scale=.4]{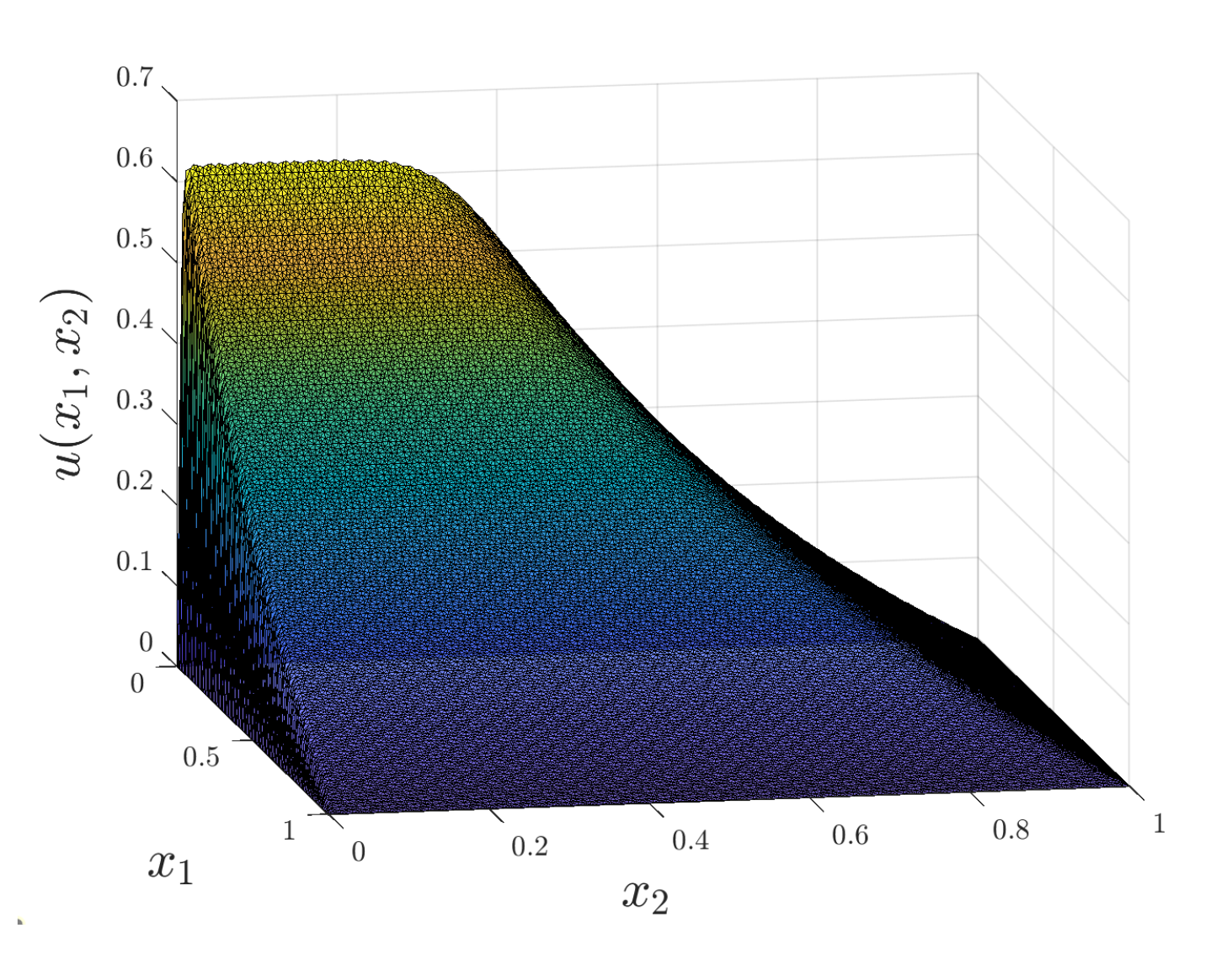}
 \caption{MFD solution to modified test \eqref{eq:test1}. }
	\label{fig:grad_layer_HH}
\end{subfigure}
\begin{subfigure}{.49\textwidth}
  \centering
\includegraphics[scale=.4]{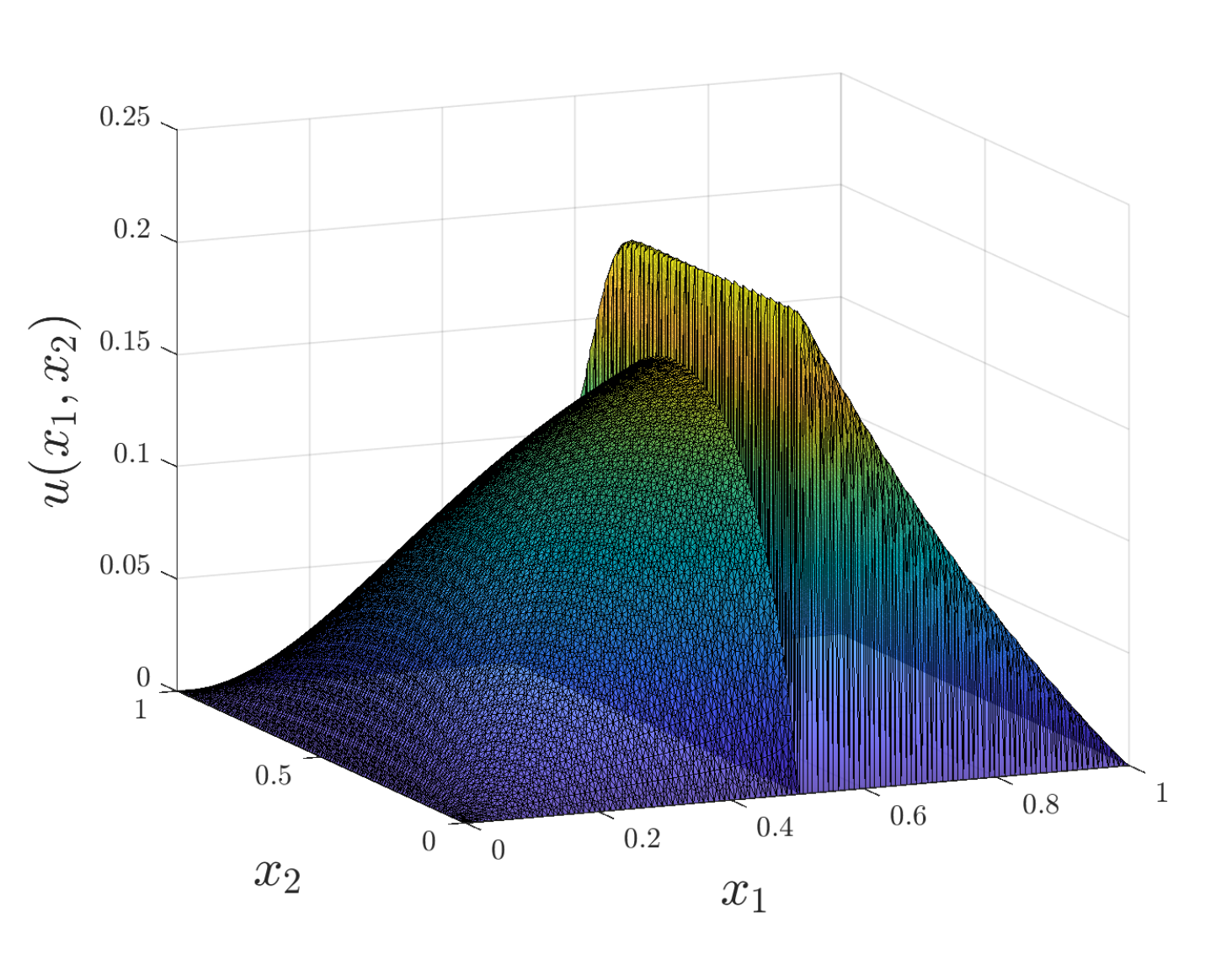}
	\caption{MFD solution to modified test \eqref{eq:test2}.}
  \label{fig:grad_jump_HH}
\end{subfigure}
\caption{MFD solutions to boundary layer and internal layer with Helmholtz decomposition for $\bm \beta$ in $H(\grad)$.}
\label{fig:grad_mono_HH}
\end{figure}

Finally, we note that while our method relies on being able to write $\bm \theta = \bm \beta / \alpha$ in terms of a potential function, we do consider a way to handle the case where no such potential function exists, by using a Helmholtz decomposition, $\bm \theta = \grad \varphi + \curl \bm \psi$ (see Remark \ref{remark:helmholtz}). Here, preliminary results are presented for the Helmholtz decomposition case for the scalar convection-diffusion equation. The tests given below are identical to the tests presented in \eqref{eq:test4} of Section \ref{subsec:convergence}, and \eqref{eq:test1} and \eqref{eq:test2} of Section \ref{subsec:grad_layer}, with the small modification of adding 
\begin{equation}
\curl \psi = \curl \left( .1 \sin(x_1 x_2) \right)
\end{equation}
to the convection parameter, $\bm \beta$.
Figure \ref{fig:test_conv_HH} demonstrates that the convergence result holds even with the $\curl \psi$ added to the convection. From Figures \ref{fig:grad_layer_HH} and \ref{fig:grad_jump_HH}, we can see that monotonicity is preserved even with the modification to the convection coefficient. 
%\casey{Should we mention that we expect that $\curl \psi$ needs to be ``small enough" for this to work?}

\section{Conclusions}\label{sec:conc_CD}
This manuscript presented a stable MFD method for the general $H(\mathfrak{D})$ convection-diffusion equation. The discrete flux operators allow for the convection-diffusion equations to be recast as a pure diffusion problem, and preserves the de Rham complex, guaranteeing that the method is structure-preserving regardless of the integration rule used to compute the exponential integrals. The MFD method can be thought of as a scaled, mass-lumped FE method. Using a FE framework allowed for a clear path for proving well-posedness via Babu\v{s}ka theory and deriving error estimates of the MFD scheme. Furthermore, the method is provably stable, meaning that as the diffusion coefficient tends to zero, $\alpha \to 0$, the discrete maximum principle or monotonicity property holds and the method will not suffer from numerical oscillations when the solution has shocks or boundary layers. 

While only the lowest-order FE and MFD methods are studied here, higher-order methods  \cite{Castillo_high_MFD, Lipkinov_high_MFD, MFDart3, WuZikatanov_SAFE} require further investigation to see if a similar framework holds. Additionally, since the flux operators satisfy a de Rham complex, iterative methods can be developed to efficiently solve the MFD system based on the FE theory. For instance, one could explore traditional multigrid methods for $H(\grad)$ and auxiliary space preconditioners for $\bm H(\curl)$ and $\bm H(\divg)$. Some preliminary work has been done on this for the SAFE method, and those ideas could be extended to apply to MFD.  Finally, we note that though the analysis presented only considers the case where $\bm \theta = \grad \varphi$, we show results validating the stability when the coefficients can be represented by a Helmholz decomposition, $\bm \theta = \grad \varphi + \curl \bm \psi$.  Completing the analysis for this case requires further exploration.

\section*{Acknowledgments}

Sandia National Laboratories is a multimission laboratory managed and operated by National Technology and Engineering Solutions of Sandia, LLC, a wholly owned subsidiary of Honeywell International, Inc., for the U.S. Department of Energy's National Nuclear Security Administration under contract DE-NA0003525. This paper describes objective technical results and analysis. Any subjective views or opinions that might be expressed in the paper do not necessarily represent the views of the U.S. Department of Energy or the United States Government.

N. Trask has been supported by the U.S. Department of Energy, Office of Advanced Scientific Computing Research under the U.S. Department of Energy, Office of Advanced Scientific Computing Research under the Early Career Research Program. A. Huang has been supported under the Sandia National Laboratories Laboratory Directed Research and Development (LDRD) program.

\bibliographystyle{siam}
\bibliography{references}
\end{document}